%% file: 1bitMatrixCompletion_2017_01_27.tex
\documentclass[11pt, a4paper]{article}

\input{packages_notes}

\input{commandes_guillaume}
\usepackage{multirow}

\begin{document}

\title{Estimation bounds and sharp oracle inequalities of regularized procedures with  Lipschitz loss functions}
\author{Pierre Alquier${}^{1,3,4}$ , Vincent Cottet${}^{1,3,4}$ , Guillaume Lecué${}^{2,3,4}$  
\\
\small{
(1) CREST, ENSAE, Université Paris Saclay}
\\
\small{
(2) CREST, CNRS, Université Paris Saclay
}
}

\footnotetext[3] {Email:
\{pierre.alquier\}, \{vincent.cottet\}, \{guillaume.lecue\}@ensae.fr}
\footnotetext[4] {
The authors gratefully acknowledge financial support
from Labex ECODEC (ANR - 11-LABEX-0047). Author n. 1 also acknowledge
financial support from the
research
programme {\it New Challenges for New Data} from LCL and GENES, hosted by
the {\it Fondation du Risque}. Author n. 3 acknowledge financial support from the "Chaire Economie et Gestion des Nouvelles Donn\'ees", under the auspices of Institut Louis Bachelier, Havas-Media and Paris-Dauphine.
}

\maketitle

\begin{abstract}
We obtain estimation error rates and sharp oracle inequalities for regularization procedures of the form
\begin{equation*}
  \hat f \in\argmin_{f\in
    F}\left(\frac{1}{N}\sum_{i=1}^N\ell(f(X_i), Y_i)+\lambda \norm{f}\right)
\end{equation*}
when $\norm{\cdot}$ is any norm, $F$ is a convex class of functions and $\ell$ is a Lipschitz loss function satisfying a Bernstein condition over $F$.  We explore both the bounded and subgaussian stochastic frameworks for the distribution of the $f(X_i)$'s, with no assumption on the distribution of the $Y_i$'s. The general results rely on two main objects: a complexity function, and a sparsity equation, that depend on the specific setting in hand (loss $\ell$ and norm $\norm{\cdot}$).

As a proof of concept, we obtain minimax rates of convergence in the following problems: 1) matrix completion with any Lipschitz loss function, including the hinge and logistic loss for the so-called 1-bit matrix completion instance of the problem, and quantile losses for the general case, which enables to estimate any quantile on the entries of the matrix; 2) logistic LASSO and variants such as the logistic SLOPE; 3) kernel methods, where the loss is the hinge loss, and the regularization function is the RKHS norm.  
\end{abstract}

\section{Introduction}
\label{sec:intro}

Many classification and prediction problems are solved in practice by regularized empirical risk minimizers (RERM). The risk is measured by a loss function and the quadratic loss function is the most popular function for regression. It has been extensively studied (cf. \cite{LM_sparsity,MR2829871} among others). Still many other loss functions are popular among practitioners and are indeed extremely useful in specific situations.

First, let us mention the quantile loss in regression problems. The $0.5$-quantile loss (also known as absolute or $L_1$ loss) is known to provide an indicator of conditional central tendency more robust to outliers than the quadratic loss. An alternative to the absolute loss for robustification is provided by the Huber loss. On the other hand, general quantile losses are used to estimate conditional quantile functions and are extremely useful to build confidence intervals and measures of risk, like {\it Values at Risk} (VaR) in finance.
 
Let us now turn to classification problems. The natural loss in this context, the so called $0/1$ loss, leads very often to computationally intractable estimators. Thus, it is usually replaced by a convex loss function, such as the hinge loss or the logistic loss. A thorough study of convex loss functions in classification can be found in~\cite{zhang2004statistical}.

All the aforementioned loss functions (quantile, Huber, hinge and logistic) share a common property: they are Lipschitz functions. This motivates a general study of RERM with any Lipschitz loss. Note that some examples were already studied in the literature: the $\norm{\cdot}_1$-penalty with a quantile loss was studied in~\cite{belloni2011l1} under the name ``quantile LASSO'' while the same penalty with the logistic loss was studied in~\cite{van2008high} under the name ``logistic LASSO'' (cf. \cite{MR3526202}). The ERM strategy with Lipschitz proxys of the $0/1$ loss are studied in~\cite{MR1892654}. The loss functions we will consider in the examples of this paper are reminded below:
\begin{enumerate}
  \item \textbf{hinge loss:} $\ell(y',y) = (1-yy')_+ = \max(0, 1-yy')$ for every $y\in\{-1,+1\}, y' \in \R$,
  \item \textbf{logistic loss:} $\ell(y',y) = \log(1 + \exp(-y y'))$ for every $y\in\{-1,+1\}, y' \in \R$;
  \item \textbf{quantile regression loss:} for some parameter $\tau\in(0,1)$, $\ell(y',y) = \rho_\tau(y-y')$  for every $y\in\R, y' \in \R$ where $\rho_\tau(z)= z (\tau - I(z\leq0))$ for all $z\in\R$.
  \end{enumerate}

The two main theoretical results of the paper, stated in Section~\ref{sec:theoretical_results}, are general in the sense that they do not rely on a specific loss function or a specific regularization norm. We develop two different settings that handle different assumptions on the design. In the first one, we assume that the family of predictors is subgaussian; in the second setting we assume that the predictors are uniformly bounded, this setting is well suited for classification tasks, including the 1-bit matrix completion problem. The rates of convergence rely on quantities that measure the complexity of the model and the size of the subdifferential of the norm.

To be more precise, the method works for any regularization function as long as it is a norm. If this norm has some sparsity inducing power, like the $\ell_1$ or nuclear norms, thus the statistical bounds depend on the underlying sparsity around the oracle because the subdifferential is large. We refer these bounds as \emph{sparsity dependent bounds}. If the norm does not induce sparsity, it is still possible to derive bounds that are now depending on the norm of the oracle because the subdifferential of the norm is very large in $0$. We call it \emph{norm dependent bounds} (aka ``complexity dependent bounds'' in \cite{LM_comp}).

We study many applications that give new insights on diverse problems: the first one is a classification problem with logistic loss and LASSO or SLOPE regularizations. We prove that the rate of the SLOPE estimator is minimax in this framework. The second one is about matrix completion. We derive new excess risk bounds for the 1-bit matrix completion issue with both logistic and hinge loss. We also study the quantile loss for matrix completion and  prove it reaches sharp bounds. We show several examples in order to assess the general methods as well as simulation studies. The last example involves the SVM and proves that  ``classic'' regularization method with no special sparsity inducing power can be analyzed in the same way as sparsity inducing regularization methods. 

A remarkable fact is that no assumption on the output $Y$ is needed (while most results for the quadratic loss rely on an assumption of the tails of the distribution of $Y$). Neither do we assume any statistical model relating the ``output variable'' $Y$ to the ``input variable'' $X$.

\paragraph{Mathematical background and notations.} The observations are $N$ i.i.d pairs $(X_i,Y_i)_{i=1}^N$ where $(X_i,Y_i) \in \cX \times \cY$ are distributed according to $P$. We consider the case where $\cY$ is a subset of $\bR$ and let $\mu$ denote the marginal distribution of $X_i$. Let $L_2$ be the set of real valued functions $f$ defined on $\cX$ such that $\bE f(X)^2 < +\infty$ where the distribution of $X$ is $\mu$. In this space, we define the $L_2$-norm as $\normL{f} = (\bE f(X)^2)^{1/2}$ and the $L_\infty$ norm such that $\norm{f}_{L_\infty} = {\rm esssup}(|f(X)|)$. We consider a set of predictors $F\subseteq E$, where $E$ is a subspace of $L_2$ and $\norm{\cdot}$ is a norm over $E$ (actually, in some situations we will simply have $F=E$, but in some natural examples we will consider bounded set of predictors, in the sense that $\sup_{f\in F} \norm{f}_{L_\infty}<\infty$, which implies that $F$ cannot be a subspace of $L_2$).

For every $f\in F$, the loss incurred when we predict $f(x)$, while the true output / label is actually $y$, is measured using a loss function $\ell$: $\ell(f(x),y)$. For short, we will also use the notation $\ell_f(x,y) = \ell(f(x),y) $ the loss function associated with $f$. In this work, we focus on loss functions that are nonnegative, and Lipschitz, in the following sense.
\begin{assumption}[Lipschitz loss function]\label{ass:lipschitz}
 For every $f_1, f_2\in F$, $x\in \cX$ and $y\in\bR$, we have
 \begin{equation*}
 \big|\ell(f_1(x),y) - \ell(f_2(x),y)\big|\leq |f_1(x)-f_2(x)|.
 \end{equation*}
 \end{assumption}
Note that we chose a Lipschitz constant equal to one in Assumption~\ref{ass:lipschitz}.  This can always be achieved by a proper normalization of the loss function. We define the oracle predictor as
\begin{equation*}
f^*\in \argmin_{f\in F}P \ell_f \mbox{ where}\footnote{Note that without any assumption on $Y$ it might be that $P \ell_f = \bE \ell_f(X,Y) =\infty$ for any $f\in F$. Our results remain valid in this case, but it is no longer possible to use the definition $f^*\in \argmin_{f\in F}P \ell_f$. A general definition is as follows: fix any $f_0\in F$. Note that for any $f\in F$, $\mathbb{E}[\ell_f(X,Y)-\ell_{f_0}(X,Y)]] \leq \mathbb{E}|(f-f_0)(X)|< \infty$ under the assumptions on $F$ that will be stated in Section~\ref{sec:theoretical_results}. It is then possible to define $f^*$ as any minimizer of $\mathbb{E}[\ell_f(X,Y)-\ell_{f_0}(X,Y)]]$. This definition obviously coincides with the defintion $f^*\in \argmin_{f\in F}P \ell_f$ when $P \ell_f$ is finite for some $f\in F$.} \: P \ell_f = \bE \ell_f(X,Y)
\end{equation*}and $(X,Y)$ is distributed like the $(X_i,Y_i)$'s.
The objective of machine learning is to provide an estimator $\hat{f}$ that predicts almost as well as $f^*$. We usually formalize this notion by introducing the excess risk $\cE(f)$ of $f\in F$ by
 \begin{equation*}
 \cL_f = \ell_f-\ell_{f^*} \mbox{ and } \cE(f) = P\cL_f.
 \end{equation*}
Thus we consider the estimator of the form
\begin{equation}\label{eq:estimator}
 \hM \in \argmin_{f\in F}\left\lbrace P_N \ell_f + \lambda \norm{f} \right\rbrace
 \end{equation}
where $P_N \ell_f = (1/N)\sum_{i=1}^N \ell_f(X_i, Y_i)$ and $\lambda$ is a regularization parameter to be chosen. Such estimators are usually called Regularized Empirical Risk Minimization procedure (RERM).
 
For the rest of the paper, we will use the following notations: let $rB$ and $rS$ denote the radius $r$ ball and sphere for the norm $\norm{\cdot}$, i.e. $rB = \{f\in E:\norm{f} \leq r\}$ and $rS = \{f\in E:\norm{f} = r\}$. For the $L_2$-norm, we write $rB_{L_2} = \{f\in L_2:\normL{f} \leq r\}$ and $rS_{L_2} = \{f\in L_2:\normL{f} = r\}$ and so on for the other norms. 

Even though our results are valid in the general setting introduced above, we will develop the examples mainly in two directions that we will refer to \emph{vector} and \emph{matrix}. The \emph{vector} case involves $\cX$ as a subset of $\R^p$; we then consider the class of linear predictors, i.e. $E=\{\inr{t,\cdot}, t \in \R^p\}$. In this case, we denote for $q\in [1,+\infty]$, the $l_q$-norm in $\bR^p$ as $\norm{\cdot}_{l_q}$. The \emph{matrix} case is also referred as the trace regression model: $X$ is a random matrix in $\R^{m\times T}$ and we consider the class of linear predictors $E = \{\inr{M,\cdot}, M \in \R^{m\times T}\}$ where $\inr{A,B} = \textrm{Trace}(A^\top B)$ for any matrices $A,B$ in $\R^{m\times T}$. The norms we consider are then, for $q \in [1,+\infty[$, the Schatten-$q$-norm for a matrix: $\forall M \in \bR^{m\times T}, \norm{M}_{S_q} = (\sum \sigma_i(M)^q)^{1/q}$ where $\sigma_1(M)\geq \sigma_2(M)\geq \cdots$ is the family of the singular values of $M$. The Schatten-$1$ norm is also called trace norm or nuclear norm. The Schatten-$2$ norm is also known as the Frobenius norm. The $S_\infty$ norm, defined as $\norm{M}_{S_\infty} = \sigma_1(M)$ is known as the operator norm. 

The notation ${\bf C}$ will be used to denote positive constants, that might change from one instance to the other. For any real numbers $a,b$, we write $a \lesssim b$ when there exists a positive constant ${\bf C}$ such that $a \leq{\bf C} b$. When $a \lesssim b$ and $b \lesssim a$, we write $a \sim b$.

\paragraph{Proof of Concept.} We now present briefly one of the outputs of our global approach: an oracle inequality for the $1$-bit matrix completion problem with hinge loss (we refer the reader to Section~\ref{sec:hinge_} for a detailed exposition of this example). While the general matrix completion problem has been extensively studied in the case of a quadratic loss, see~\cite{MR2906869,LM_sparsity} and the references therein, we believe that there is no satisfying solution to the so-called $1$-bit matrix completion problem, that is for binary observations $\cY = \{-1,+1\}$. Indeed, the attempts in~\cite{srebro2004maximum,Cottet2016} to use the hinge loss did not lead to rank dependent learning rates. On the other hand,~\cite{lafond2014probabilistic} studied RERM procedure using a statistical modeling approach and the logistic loss. While these authors prove optimal rates of convergence of their estimator with respect to the Frobenius norm,  the excess classification risk, is not studied in their paper. However we believe that the essence of machine learning is to focus on this quantity -- it is directly related to the average number of errors in prediction.

From now on we assume that $\cY = \{-1,+1\}$ and we consider the \emph{matrix} framework. In matrix completion, we write the observed location as a mask matrix $X$: it is an element of the canonical basis $(E_{1,1}, \cdots, E_{m,T})$ of $\R^{m\times T}$ where for any $(p,q)\in\{1,\ldots, m\}\times \{1, \ldots,T\}$ the entry of $E_{p,q}$ is $0$ everywhere except for the $(p,q)$-th entry where it equals to $1$. We assume that there are constants $0<\underbar{c}\leq \bar{c}<\infty$ such that, for any $(p,q)$, $\underbar{c} / (mT) \leq \mathbb{P}(X=E_{p,q}) \leq \bar{c} / (mT)$ (this extends the uniform sampling distribution for which $\underbar{c} = \bar{c} = 1$). These assumptions are encompassed in the following definition.

\begin{assumption}[Matrix completion design]\label{ass:matrix_design}
The sample size $N$ is in $\{\min(m,T),\dots,\max(m,T)^2\}$ and $X$ takes value in the canonical basis $(E_{1,1}, \cdots, E_{m,T})$ of $\R^{m\times T}$. There are positive constants $\underbar{c}, \bar{c}$ such that for any $(p,q)\in\{1,\dots,m\}\times\{1,\dots,T\}$,
\begin{equation*}
 \underbar{c} / (mT) \leq \mathbb{P}(X=E_{p,q}) \leq \bar{c} / (mT).
\end{equation*}
\end{assumption}

A predictor can be seen, for this problem, as the natural inner product with a real $m\times T$ matrix: $f(X) = \inr{M,X} = \textrm{Tr}(X^\top M)$. The class $F$ that we consider in Section~\ref{sec:hinge_} is the set of linear predictors where every entry of the matrix is bounded: $F = \{\inr{\cdot, M}: M\in b B_\infty\}$ where $b B_\infty = \{M=(M_{pq}):\max_{p,q}|M_{pq}|\leq b\}$ for a specific $b$. This set is very common in matrix completion studies. But it is especially natural in this setting: indeed, the Bayes classifier, defined by $\overline{M} = \argmin_{M\in \R^{m\times T}} \E \left(1-Y \inr{X, M}\right)_+$, has entries in $[-1,1]$. So, by taking $b=1$ in the definition of $F$, we ensure that the oracle $M^*=\argmin_{M\in \mathbb{R}^{m\times T}} \E \left(1-Y \inr{X, M}\right)_+$ satisfies $M^*=\bar{M}$, so there would be no point in taking $b>1$. We will therefore consider the following RERM (using the hinge loss)
\begin{equation}
\label{eq:RERM_mat_comp}
\widehat{M} \in\argmin_{M\in B_\infty}\left(\frac{1}{N}\sum_{i=1}^N \left(1-Y_i\inr{X_i,M}\right)_{+} + \lambda \norm{M}_{S_1}\right)
\end{equation}
where $\lambda>0$ is some parameter to be chosen. We prove in Section~\ref{sec:hinge_} the following result.

\begin{theo}\label{theo:main_mat_comp_intro}
Assume that Assumption \ref{ass:matrix_design} holds and there is $\tau>0$ such that, for any $(p,q) \in \{1,\dots,m\}\times\{1,\dots,T\}$, 
\begin{equation}
\left|\overline{M}_{p,q}-\frac{1}{2}\right| \geq \tau.
\end{equation}
There is a $c_0(\underbar{c},\bar{c})>0$, that depends only on $\underbar{c}$ and $\bar{c}$, and that is formally introduced in Section~\ref{sec:hinge_} below, such that if one chooses the regularization parameter
\begin{equation*}
\lambda = c_0(\underbar{c},\bar{c})\sqrt{\frac{\log(m+T)}{N \min(m,T)}}
\end{equation*}
then, with probability at least 
\begin{equation}\label{eq:proba_theo_main}
 1- {\bf C} \exp\left(-{\bf C} {\rm rank}(\overline{M}) \max(m,T) \log(m+T) \right),
 \end{equation} 
the RERM estimator $\widehat{M}$ defined in~\eqref{eq:RERM_mat_comp} satisfies for every $1\leq p\leq2$,
\begin{equation*}
\frac{1}{(mT)^{\frac{1}{p}}}\norm{\widehat{M}-\overline{M}}_{S_p}
\leq {\bf C} {\rm rank}(\overline{M})^{\frac{1}{p}} \sqrt{\frac{\log(m+T)}{N}}
\frac{\max(m,T)^{1-\frac{1}{p}}}{\min(m,T)^{\frac{1}{p}-\frac{1}{2}} }
\end{equation*}
and as a special case for $p=2$,
\begin{equation}
\frac{1}{\sqrt{mT}}\norm{\widehat{M}-\overline{M}}_{S_2}
\leq
{\bf C} \sqrt{\frac{{\rm rank}(\overline{M})\max(m,T) \log(m+T)}{N}}
\end{equation}
and its excess hinge risk is such that
    \begin{equation*}
 \cE_{hinge}(\widehat{M}) = \E(1-Y\inr{X, \widehat{M}})_+- \E(1-Y\inr{X, \overline{M}})_+\leq {\bf C} \frac{\rank(\overline{M})\max(m,T)\log(m+T)}{N}
\end{equation*} 
where the notation ${\bf C}$ is used for constants that might change from one instance to the other but depend only on $\underbar{c}$, $\bar{c}$ and $\tau$.
\end{theo}

The excess hinge risk bound from Theorem~\ref{theo:main_mat_comp_intro} is of special interest as it can be related to the classic excess $0/1$ risk. The excess $0/1$  risk of a procedure is really the quantity we want to control since it measures the difference between the average number of mistakes of a procedure with the best possible theoretical classification rule. Indeed, let us define the $0/1$ risk of $M$ by $R_{0/1}(M)= \mathbb{P}[Y \neq {\rm sign}(\left<M,X\right>)]$. It is clear that $\overline{M} \in \argmin_{M\in\R^{m\times T}} R_{0/1}(M)$. Then, it follows from Theorem 2.1 in~\cite{zhang2004statistical} that for some universal constant $c>0$, for every $M\in\R^{m\times T}$,
\begin{equation*}
R_{0/1}(M) - \inf_{M \in B_\infty} R_{0/1}(M) \leq c \cE_{hinge}(M).
\end{equation*} 
Therefore, the RERM from \eqref{eq:RERM_mat_comp} for the choice of regularization parameter $\lambda$ as in Theorem~\ref{theo:main_mat_comp_intro} satisfies with probability larger than in \eqref{eq:proba_theo_main},
\begin{equation}\label{equation-zhang}
 \cE_{0/1}(\widehat{M}) = R_{0/1}(\widehat{M}) - \inf R_{0/1}(M) \leq  {\bf C} \frac{\rank(\overline{M})\max(m,T)\log(m+T)}{N}
\end{equation}
where ${\bf C}$ depends on $c$, $\underbar{c}$, $\bar{c}$ and $\tau$. This yields a bound on the average of excess number of mistakes of $\widehat{M}$. To our knowledge such a prediction bound was not available in the literature on the $1$-bit matrix completion problem. Let us compare Theorem~\ref{theo:main_mat_comp_intro} to the main result in \cite{lafond2014probabilistic}. In~\cite{lafond2014probabilistic}, the authors focus on the estimation error $\|\widehat{M}-M^*\|_{S_2}$, which seems less relevant for practical applications. In order to connect such a result to the excess classification risk, one can use the results in~\cite{zhang2004statistical} and in this case, the best bound that can be derived is of the order of $\sqrt{\rank(M^*)\max(m,T)/N}$. Note that other authors focused on the classification error:~\cite{srebro2004maximum} proved an excess error bound, but the bound does not depend on the rank of the oracle. The rate $\rank(M^*)\max(m,T)/N$ derived from Theorem~\ref{theo:main_mat_comp_intro} for the $0/1$-classification excess risk was only reached in~\cite{Cottet2016}, but in the very restrictive noiseless setting, which is equivalent to $\inf_{M} R_{0/1}(M) = 0$.

We hope that this example convinced the reader of the practical interest of the general study of $\hat{f}$ in~\eqref{eq:estimator}. The rest of the paper is organized as follows. In Section~\ref{sec:theoretical_results} we introduce the concepts necessary to the general study of~\eqref{eq:estimator}: namely, a complexity parameter, and a sparsity parameter. Thanks to these parameters, we define the assumptions necessary to our general results: the Bernstein condition, which is classic in learning theory to obtain fast rates~\cite{LM_sparsity}, and a stochastic assumption on $F$ (subgaussian, or bounded). The general results themselves are eventually presented. The remaining sections are devoted to applications of our results to different estimation methods: the logistic LASSO and logistic SLOPE in Section \ref{sec:the_logistic_lasso}, matrix completion in Section~\ref{sec:hinge_} and Support Vector Machines (SVM) in Section \ref{sec:any_lipschitz_loss_rkhs_norm}. For matrix completion, the optimality of the rates for the logistic and the hinge loss, that were not known, is also derived. In Section~\ref{sec:bernstein_margin_condition} we discuss the Bernstein condition for the three main loss functions of interest: hinge, logistic and quantile. 
  
\section{Theoretical Results} 
\label{sec:theoretical_results}

\subsection{Applications of the main results: the strategy} 
\label{sec:applications_of_the_main_results_the_strategy}
The two main theorems in Sections \ref{sec:subgaussian} and \ref{sec:bounded} below are general in the sense that they allow the user to deal with any (nonnegative) Lipschitz loss function and any norm for regularization, but they involve quantities that depend on the loss and the norm. The aim of this Section is first to provide the definition of these objects and some hints on their interpretation, through examples. The theorems are then stated in both settings. Basically, the assumptions for the theorems are of three types:
\begin{enumerate}
 \item the so-called Bernstein condition, which is a quantification of the identifiability condition. It basically tells how the excess risk $\mathcal{E}(f) = P\cL_f=P(\ell_f-\ell_{f^*})$ is related to the $L_2$ norm $\|f-f^*\|_{L_2}$.
 \item a stochastic assumption on the distribution of the $f(X)$'s for $f\in F$. In this work, we consider both a subgaussian assumption and a uniform boundedness assumption. Analysis of the two setups differ only on the way the ``statistical complexity of $F$'' is measured (cf. below the functions $r(\cdot)$ in Definition~\ref{def:function_r} and Definition~\ref{def:function_r_boundedness}).
 \item finally, we introduce a sparsity parameter as in \cite{LM_sparsity}. It reflects how the norm $\|\cdot\|$ used as a regularizer can induce sparsity - for example, think of the ``sparsity inducing power'' of the $l_1$-norm used to construct the LASSO estimator.
\end{enumerate}

Given a scenario, that is a loss function $\ell$, a random design $X$, a convex class $F$ and a regularization norm, statistical results (exact oracle inequalities and estimation bounds w.r.t. the $L_2$ and regularization norms) for the associated regularized estimator together with the choice of the regularization parameter follow from the derivation of the three parameters $(\kappa, r, \rho^*)$ as explained in the next box together with Theorem~\ref{thm:main_subgaussian_simple} and Theorem~\ref{thm:main_bounded_simple}. 

\begin{center}
\begin{mybox}{Application of the main results}
\begin{enumerate}
	\item[1.] find the \textbf{Bernstein parameter} $\kappa\geq1$ and $A>0$ associated to the loss and the class $F$;
	\item[2.]  compute the \textbf{Complexity function}
\begin{equation*}
r(\rho) = \left[\frac{A \rho {\rm comp}(B)}{\sqrt{N}}\right]^{1/2\kappa}
\end{equation*}where ${\rm comp}(B)$ is defined either through the Gaussian mean width  $w(B)$, in the subgaussian case, or the Rademacher complexity ${\rm Rad}(B)$, in the bounded case;
\item[3.] Compute the sub-differential $\partial\norm{\cdot}(f^*)$ of $\norm{\cdot}$ at the oracle $f^*$ (or in the neighborhood $f^*+(\rho/20)B$ for approximately sparse oracles) and solve the \textbf{sparsity equation} ``find $\rho^*$ such that $\Delta(\rho^*)\geq 4\rho^*/5$''.
\item[4.] Apply  Theorem~\ref{thm:main_subgaussian_simple} in the subgaussian framework and Theorem~\ref{thm:main_bounded_simple} in the bounded framework. In each case, with large probability,
\begin{equation*}
 \norm{\hat{f}-f^*} \leq \rho^*\text{, } \norm{\hat{f}-f^*}_{L_2} \leq r(2\rho^*) \text{ and } \mathcal{E}(\hat{f}) \leq {\bf C} \left[r(2\rho^*)\right]^{2\kappa}.
\end{equation*}
\end{enumerate}
\end{mybox} 
\end{center}

For the sake of simplicity, we present the two settings in different subsections with both the exact definition of the complexity function and the theorem. As the sparsity equation is the same in both settings, we define it before even though it involves the complexity function.

\subsection{The Bernstein condition}
\label{subsec:bernstein}

The first assumption needed is called \emph{Bernstein} assumption and is very classic in order to deal with Lipschitz loss. 
 
\begin{assumption}[Bernstein condition]\label{ass:margin}
 There exists $\kappa\geq1$ and $A>0$ such that for every  $f\in F$, $\norm{f-f^*}_{L_2}^{2\kappa}\leq A P\cL_{f}$.
 \end{assumption}
  
The most important parameter is $\kappa$ and will be involved in the rate of convergence. As usual fast rates will be derived when $\kappa=1$. In many situations, this assumption is satisfied and we present various cases in Section~\ref{sec:bernstein_margin_condition}. In particular, we prove that it is satisfied with $\kappa=1$ for the logistic loss in both bounded and Gaussian framework, and we exhibit explicit conditions to ensure that Assumption~\ref{ass:margin} holds for the hinge and the quantile loss functions. 

We call Assumption~\ref{ass:margin} a {\it Bernstein condition} following \cite{MR2240689} and that it is different from the margin assumption from \cite{Mammen1999,Tsy04}: in the so-called margin assumption, the oracle $f^*$ in $F$ is replaced by the minimizer $\overline{f}$ of the risk function  $f \to P \ell_f $ over all measurable functions $f$, sometimes called the Bayes rules. We refer the reader to Section~\ref{sec:bernstein_margin_condition} and to the discussions in \cite{MR2933668} and Chapter~1.3 in \cite{HDR-lecue} for more details on the difference between the margin assumption and the Bernstein condition.

\begin{rmk}
The careful reader will actually realize that the proof of Theorem~\ref{thm:main_subgaussian_simple} and Theorem~\ref{thm:main_bounded_simple} requires only a weaker version of this assumption, that is: there exists $\kappa\geq1$ and $A>0$ such that for every  $f\in \mathcal{C}$, $\norm{f-f^*}_{L_2}^{2\kappa}\leq A P\cL_{f}$,
where $\mathcal{C}$ is defined in terms of the complexity function $r(\cdot)$ and the sparsity parameter $\rho^*$ to be defined in the next subsections,
\begin{equation}\label{eq:set_cC}
\cC :=\left\{f\in F: \norm{f-f^*}_{L_2}\geq r(2\norm{f-f^*}) \mbox{ and } \norm{f-f^*}\geq \rho^*\right\}.
\end{equation}
Note that the set $\cC$ appears to play a central role in the analysis of regularization methods, cf. \cite{LM_sparsity}. However, in all the examples presented in this paper, we prove that the Bernstein condition holds on the entire set $F$.
\end{rmk}

\subsection{The complexity function $r(\cdot)$}
\label{sec:cplxity_fn}

The complexity function $r(\cdot)$ is defined by
\begin{equation*}
 \forall \rho>0,\quad r(\rho) = \left[\frac{A \rho {\rm comp}(B)}{\sqrt{N}}\right]^{1/2\kappa}
\end{equation*}
where $A$ is the constant in Assumption~\ref{ass:margin} and where ${\rm comp}(B)$ is a measure of the complexity of the unit ball $B$ associated to the regularization norm. Note that this complexity measure will depend on the stochastic assumption of $F$. In the bounded setting, ${\rm comp}(B) = C {\rm Rad}(B)$ where $C$ is an absolute constant and ${\rm Rad}(B)$ is the Rademacher complexity of $B$ (whose definition will be reminded in Subsection~\ref{sec:bounded}). In the subgaussian setting, ${\rm comp}(B) = C L w(B)$ where $C$ is an absolute constant, $L$ is the subgaussian parameter of the class $F-F$ and $w(B)$ is the Gaussian mean-width of $B$ (here again, exact definitions of $L$ and $w(B)$ will be reminded in Subsection~\ref{sec:subgaussian}).

Note that sharper (localized) versions of $r(\cdot)$ are provided in Section~\ref{sec:proof_of_theorem_ref_thm_main}. However, as it is the simplest version that is used in most examples, we only introduce this version for now.

\subsection{The sparsity parameter $\rho^*$}
\label{sec:sparsity_param}
The size of the sub-differential of the regularization function $\norm{\cdot}$ in a neighborhood of the oracle $f^*$ will play as well a central role in our analysis. We recall now its definition: for every $f\in F$
\begin{equation*}
\partial\norm{\cdot}(f) = \left\{g \in E : \norm{f+h} - \norm{f}\geq \inr{g, h} \mbox{ for all } h\in E\right\}.
\end{equation*}
It is well-known that $\partial \norm{\cdot}(f)$ is a subset of the unit sphere of the dual norm of $\norm{\cdot}$ when $f\neq 0$. Note also that when $f=0$, $\partial\norm{\cdot}(f)$ is the entire unit dual ball, a fact we will also use in two situations, either when  the regularization norm has no ``sparsity inducing power'' -- in particular, when it is a smooth function as in the RKHS case treated in Section~\ref{sec:any_lipschitz_loss_rkhs_norm}; or when one wants extra \emph{norm dependent} upper bounds (cf. \cite{LM_comp} for more details where these bounds are called \emph{complexity dependent}) in addition to \emph{sparsity dependent} upper bounds. In the latter, the statistical bounds that we get are the minimum between an error rate that depends on the notion of sparsity naturally associated to the regularization norm (when it exists) and an error rate that depends on $\norm{f^*}$. 

\begin{definition}[From \cite{LM_sparsity}]\label{def:sparsity_param}
The \textbf{sparsity parameter} is the function $\Delta(\cdot)$ defined by 
\begin{equation*}
\Delta(\rho) = \inf_{h \in \rho S \cap r(2\rho) B_{L_2}} \sup_{g \in \Gamma_{f^*}(\rho)}\inr{h, g}
\end{equation*}where $\Gamma_{f^*}(\rho) = \bigcup_{f\in f^*+(\rho/20)B} \partial \norm{\cdot}(f)$.
\end{definition}

Note that there is a slight difference with the definition of the \textit{sparsity parameter} from \cite{LM_sparsity} where there $\Delta(\rho)$ is defined taking the infimum over the sphere $\rho S$ intersected with a $L_2$-ball of radius $r(\rho)$ whereas in Definition~\ref{def:sparsity_param}, $\rho S$ is intersected with a $L_2$-ball of radius $r(2\rho)$. Up to absolute constants this has no effect on the behavior of $\Delta(\rho)$ and the difference comes from technical detains in our analysis (a peeling argument that we use below whereas a direct homogeneity argument was enough in \cite{LM_sparsity}).  

In the following, estimation rates with respect to the regularization norm $\norm{\cdot}$, the norm $\normL{\cdot}$ as well as sharp oracle inequalities  are given. All the convergence rates depend on a single radius $\rho^*$ that satisfies the \textit{sparsity equation} as introduced in \cite{LM_sparsity}.

\begin{definition}\label{def:sparsity_equation}
The radius $\rho^*$ is any solution of the sparsity equation:
\begin{equation}\label{eq:sparsity_equation}
\Delta(\rho^*) \geq (4/5)\rho^*.
\end{equation}
\end{definition}
Since $\rho^*$ is central in the results and drives the convergence rates, finding a solution to the sparsity equation will play an important role in all the examples that we worked out in the following. Roughly speaking, if the regularization norm induces sparsity, a sparse element in $f^*+(\rho/20) B$ (that is an element $f$ for which $\partial\norm{\cdot}(f)$ is almost extremal -- that is almost as large as the dual sphere) yields the existence of a small $\rho^*$. In this case, $\rho^*$ satisfies the sparsity equation.

In addition, if one takes $\rho =  20 \norm{f^*}$ then $0\in \Gamma_{f^*}(\rho)$ and since $\partial\norm{\cdot}(0)$ is the entire dual ball associate to $\norm{\cdot}$, one has directly that $\Delta(\rho)=\rho$ and so $\rho$ satisfies the sparsity Equation~\eqref{eq:sparsity_equation}. We will use this observation to obtain \emph{norm dependent} upper bounds, i.e. rates of convergence depending on $\norm{f^*}$ and that do not depend on any sparsity parameter. Such a bound holds for any norm; in particular, for norms with no sparsity inducing power as in Section~\ref{sec:any_lipschitz_loss_rkhs_norm}.
 
\subsection{Theorem in the subgaussian setting}
\label{sec:subgaussian}
 
First, we introduce the subgaussian framework (then we will turn to the bounded case in the next section).

\begin{definition}[Subgaussian class]\label{def:subgaussian}
We say that a class of functions $\mathcal{F}$ is $L$-subgaussian (w.r.t. $X$) for some constant  $L\geq1$ when for all $f\in \mathcal{F}$ and all $\lambda\geq1$,
\begin{equation}\label{eq:subgaussian_assum}
\E \exp\left(\lambda|f(X)|/\norm{f}_{L_2}^2\right)\leq \exp\left(\lambda^2 L^2 \right)
\end{equation}
where $\norm{f}_{L_2} = \left(\E f(X)^2\right)^{1/2}$.
\end{definition}

We will use the following operations on sets: for any $F^\prime\subset E$ and $f\in E$, 
\begin{equation*}
F^\prime + f = \{f^\prime+f: f^\prime\in F^\prime\}, \quad F^\prime - F^\prime = \{f^\prime_1-f^\prime_2: f^\prime_1,f^\prime_2\in F^\prime\} \mbox{ and } d_{L_2}(F^\prime) = \sup\left(\norm{f^\prime_1-f^\prime_2}_{L_2}: f^\prime_1,f^\prime_2\in F^\prime\right).
\end{equation*}
 
\begin{assumption}
\label{ass:sub-gauss} The class $F-F$ is $L$-subgaussian.
\end{assumption}

Note that there are many equivalent formulations of the subgaussian property of a random variable based on $\psi_2$-Orlicz norms, deviations inequalities, exponential moments, moments growth characterization,  etc. (cf., for instance Theorem~1.1.5 in \cite{MR3113826}). The one we should use later is as follows: there exists some absolute constant ${\bf C}$ such that $F-F$ is $L$-subgaussian if and only if for all $f, g\in F$ and $t\geq 1$, 
\begin{equation}\label{eq:subgaussian_deviation}
\bP[|f(X)-g(X)|\geq {\bf C} t L  \norm{f-g}_{L_2}]\leq 2 \exp(-t^2).
\end{equation}

 There are several examples of subgaussian classes. For instance,  when $F$ is a class of linear functionals $F=\{\inr{\cdot, t}:t\in T\}$ for $T\subset \R^p$ and $X$ is a random variable in $\R^p$ then $F-F$ is $L$-subgaussian in the following cases:
 \begin{enumerate}
   \item $X$ is a Gaussian vector in $\R^p$,
   \item $X=(x_j)_{j=1}^p$ has independent coordinates that are subgaussian, that is, there are constants $c_0>0$ and $c_1>0$ such that $\forall j$, $\forall t>c_0, \bP[|x_j|\geq t (\E x_j^2)^{1/2}]\leq 2 \exp(-c_1t^2)$,
   \item for $2\leq q < \infty$, $X$ is uniformly distributed over $p^{1/q} B_{l_q}$ (cf. \cite{MR2123199}), 
   \item $X=(x_j)_{j=1}^p$ is an unconditional vector (meaning that for every signs $(\eps_j)_j\in\{-1,+1\}^p$, $(\eps_j x_j)_{j=1}^p$ has the same distribution as $(x_j)_{j=1}^p$), $\E x_j^2\geq c^2$ for some $c>0$ and $\norm{X}_{l_\infty}\leq R$ almost surely then one can choose $L\leq {\bf C} R/c$ (cf. \cite{LM13}).
 \end{enumerate}


In the \emph{subgaussian framework}, a natural way to measure the \emph{statistical complexity} of the problem is via Gaussian mean-width that we introduce now. 

 \begin{definition}\label{def:gauss_mean_width}
 Let $H$ be a subset of $L_2$ and denote by $d$ the natural metric in $L_2$. Let $(G_h)_{h\in H}$ be the canonical centered Gaussian process indexed by $H$ (in particular, the covariance structure of $(G_h)_{h\in H}$  is given by $d$:  $\left(\E (G_{h_1}- G_{h_2})^2\right)^{1/2} = \left(\E(h_1(X)-h_2(X))^2\right)^{1/2}$ for all $h_1,h_2\in H$). The \textbf{Gaussian mean-width} of $H$ (as a subset of $L_2$) is
 \begin{equation*}
 w(H) = \E \sup_{h\in H} G_h.
 \end{equation*}
 \end{definition}

We refer the reader to Section~12 in \cite{MR1932358} for the construction of Gaussian processes in $L_2$. There are many natural situations where Gaussian mean-widths can be computed. To familiarize with this quantity let us consider an example in the \emph{matrix} framework. Let $H = \{\inr{M, \cdot}: \norm{M}_{S_1}\leq1\}$ be the class of linear functionals indexed by the unit ball of the $S_1$-norm and $d$ be the distance associated with the Frobenius norm (i.e. $d(\inr{\cdot, M_1}, \inr{\cdot, M_2}) = d(M_1, M_2) = \norm{M_1-M_2}_{S_2}$) then
 \begin{equation*}
 w(H) = w(B_{S_1}) = \E \sup_{\norm{M}_{S_1}\leq1 }\inr{\bG, M} = \E \norm{\bG}_{S_1}^* = \E \norm{\bG}_{S_\infty} \sim \sqrt{m+T}
 \end{equation*}where $\bG$ is a standard Gaussian matrix in $\R^{m\times T}$, $\norm{\cdot}_{S_1}^*$ is the dual norm of the nuclear norm which is the operator norm $\norm{\cdot}_{S_\infty}$.

 We are now in position to define the complexity parameter as announced previously.

\begin{definition}\label{def:function_r_simple}
The \textbf{complexity parameter} is the non-decreasing function $r(\cdot)$ defined for every $\rho\geq0$, 
\begin{equation*}
r(\rho) = \left(\frac{A C L w(B) \rho}{ \sqrt{N}}\right)^\frac{1}{2\kappa}
\end{equation*}
where $\kappa,A$ are the Bernstein parameters from Assumption~\ref{ass:margin}, $L$ is the subgaussian parameter from Assumption~\ref{ass:sub-gauss} and $C>0$ is an absolute constant (the exact value of $C$ can be deduced from the proof of Proposition~\ref{prop:stochastic}). The Gaussian mean-width $w(B)$ of $B$ is computed with respect to the the metric associated with the covariance structure of $X$, i.e. $d(f_1, f_2) = \normL{f_1-f_2}$ for every $f_1,f_2\in F$.
\end{definition}

After the computation of the Bernstein parameter $\kappa$, the complexity function $r(\cdot)$ and the radius $\rho^*$, it is now possible to explicit our main result in the sub-Gaussian framework.

\begin{theo}\label{thm:main_subgaussian_simple}
Assume that Assumption~\ref{ass:lipschitz}, Assumption~\ref{ass:margin} and Assumption~\ref{ass:sub-gauss} hold and let $C>0$ from the definition of $r(\cdot)$ in Definition \ref{def:function_r_simple}. Let the regularization parameter $\lambda$ be
\begin{equation*}
\lambda = \frac{5}{8} \frac{CLw(B)}{ \sqrt{N}}
\end{equation*}
and $\rho^*$ satisfying \eqref{eq:sparsity_equation}. Then, with probability larger than
\begin{equation}\label{eq:proba_2}
 1- {\bf C} \exp\left(- {\bf C} N^{1/{2\kappa}} (\rho^* w(B))^{(2\kappa-1)/{\kappa}}\right)
 \end{equation}
 we have
 \begin{align*}
 \norm{\hat f - f^*} \leq \rho^*,\quad
 \norm{\hat f - f^*}_{L_2} &\leq r(2\rho^*) = \left[\frac{A C L w(B)2 \rho^*}{ \sqrt{N}}\right]^{1/{2\kappa}} \mbox{ and }
 \cE(\hat f) \leq  \frac{r(2\rho^*)^{2\kappa}}{A} = \frac{C L w(B)2 \rho^*}{\sqrt{N}}
 \end{align*}
  where ${\bf C}$ denotes positive constants that might change from one instance to the other and depend only on $A$, $\kappa$, $L$ and $C$.
\end{theo}

 \begin{rmk}[Deviation parameter]
 Replacing $w(B)$ by any upper bound does not affect the validity of the result. As a special case, it is possible to increase the confidence level of the bound by replacing $w(B)$ by $w(B)+x$: then, with probability at least
 \begin{equation*}
 1- {\bf C} \exp\left(- {\bf C} N^{1/{2\kappa}} (\rho^* [w(B)+x])^{(2\kappa-1)/{\kappa}}\right)
 \end{equation*}
we have in particular
 \begin{align*}
 \norm{\hat f - f^*}_{L_2} &\leq r(2\rho^*) = \left[\frac{A C L [w(B)+x] 2 \rho^*}{ \sqrt{N}}\right]^{1/{2\kappa}} \mbox{ and }
 \cE(\hat f) \leq  \frac{r(2\rho^*)^{2\kappa}}{A} = \frac{C L [w(B)+x]2 \rho^*}{\sqrt{N}}.
 \end{align*}
 \end{rmk}

 \begin{rmk}[Norm and sparsity dependent error rates]
 Theorem~\ref{thm:main_subgaussian_simple} holds for any radius $\rho^*$ satisfying the sparsity equation~\eqref{eq:sparsity_equation}. We have noticed in Section~\ref{sec:sparsity_param} that $\rho^*=20\norm{f^*}$ satisfies the sparsity equation since in that case $0\in\Gamma_{f^*}(\rho^*)$ and so $\Delta(\rho^*)=\rho^*$. Therefore, one can apply  Theorem~\ref{thm:main_subgaussian_simple} to both $\rho^*=20\norm{f^*}$ (this leads to \emph{norm dependent} upper bounds) and to the smallest $\rho^*$ satisfying the sparsity equation~\eqref{eq:sparsity_equation} (this leads to \emph{sparsity dependent} upper bounds) at the same time. Both will lead to meaningful results (a typical example of such a combined result is Theorem~9.2 from \cite{MR2829871} or Theorem~\ref{theo:logistic_lasso} below).
 \end{rmk}

\subsection{Theorem in the bounded setting}
\label{sec:bounded}

We now turn to the \textit{bounded framework}; that is we assume that all the functions in $F$ are uniformly bounded in $L_\infty$. This assumption is very different in nature than the subgaussian assumption which is in fact a norm equivalence assumption (i.e. Definition~\ref{def:subgaussian} is equivalent to $\norm{f}_{L_2}\leq \norm{f}_{\psi_2}\leq L \norm{f}_{L_2}$ for all $f\in \cF$ where $\norm{\cdot}_{\psi_2}$ is the $\psi_2$ Orlicz norm, cf. \cite{MR1113700}). 

 \begin{assumption}[Boundedness assumption]\label{ass:bounded}
 There exist a constant $b>0$ such that for all $f\in F$, $\norm{f}_{L_\infty}\leq b$.
 \end{assumption}

The main motivation to consider the \textit{bounded setup} is for sampling over the canonical basis of a finite dimensional space like $\R^{m\times T}$ or $\R^p$. Note that this type of sampling is {\it stricto sensu} subgaussian, but with a constant $L$ depending on the dimensions $m$ and $T$, which yields sub-optimal rates. This is the reason why the results in the bounded setting are more relevant in this situation. This is especially true for the $1$-bit matrix completion problem as introduced in Section~\ref{sec:intro}. For this example, the $X_i$'s are chosen randomly in the canonical basis $(E_{1,1}, \cdots, E_{m,T})$ of $\R^{m\times T}$. Moreover, in that example, the class $F$ is the class of all linear functionals indexed by $b B_\infty$: $F=\{\inr{\cdot, M}: \max_{p,q}|M_{pq}| \leq b\}$ and therefore the study of this problem falls naturally in the bounded framework studied in this section.

Under the boundedness assumption, the natural way to measure the "statistical complexity" cannot be anymore characterized by Gaussian mean width. We therefore introduce another complexity parameter known as Rademacher complexities. This complexity measure has been extensively studied in the learning theory literature (cf., for instance,  \cite{MR2329442,MR2829871,MR2166554}).

\begin{definition}\label{def:rademacher_var}
Let $H$ be a subset of $L_2$. Let $(\eps_i)_{i=1}^N$ be $N$ i.i.d. Rademacher variables (i.e. $\bP[\eps_i=1] = \bP[\eps_i=-1]=1/2$) independent of the $X_i$'s. The \textbf{Rademacher complexity} of $H$ is 
\begin{equation*}
 {\rm Rad}(H) =  \E \sup_{f\in H}\Big|\frac{1}{\sqrt{N}}\sum_{i=1}^N \eps_i f(X_i)\Big|.
 \end{equation*}
\end{definition}

Note that when $(f(X))_{f\in H}$ is a version of the isonormal process over $L_2$ (cf. Chapter~12 in \cite{MR1932358}) restricted to $H$ then the Gaussian mean-width and the Rademacher complexity coincide: $w(H)={\rm Rad}(H)$. But, in that case, $H$ is not bounded in $L_\infty$ and, in general, the two complexity measures are different. 

There are many examples where Rademacher complexities have been computed  (cf. \cite{MR2075996}). Like in the previous \textit{subgaussian} setting the statistical complexity is given by a function $r(\cdot)$ (we use the same name in the two \textit{bounded} and \textit{subgaussian} setups because this $r(\cdot)$ function plays exactly the same role in both scenarii  even though it uses different notion of complexity).

\begin{definition}\label{def:function_r_boundedness_simple}
The \textbf{complexity parameter} is the non-decreasing function $r(\cdot)$ defined for every $\rho\geq0$  by 
\begin{equation*}
r(\rho) = \left(\frac{ C A {\rm Rad}(B)\rho}{ \sqrt{N}}\right)^\frac{1}{2\kappa}
\text{, where } C = \frac{1920}{7}.
\end{equation*}
\end{definition}

\begin{theo}\label{thm:main_bounded_simple}
Assume that Assumption~\ref{ass:lipschitz}, Assumption~\ref{ass:margin} and Assumption~\ref{ass:bounded} hold. Let the regularization parameter $\lambda$ be chosen as $\lambda = 720{\rm Rad}(B)/7\sqrt{N}$. Then, with probability larger than
 \begin{equation}\label{eq:proba_2_bounded}
 1- {\bf C} \exp\left(- {\bf C} N^{1/{2\kappa}} (\rho^* {\rm Rad}(B))^{(2\kappa-1)/{\kappa}}\right)
 \end{equation}
 we have
 \begin{align*}
 \norm{\hat f - f^*} \leq \rho^*,\quad
 \norm{\hat f - f^*}_{L_2} &\leq r(2\rho^*) = \left[\frac{C A {\rm Rad}(B)2 \rho^* }{ \sqrt{N}}\right]^{1/{2\kappa}} \mbox{ and }
 \cE(\hat f) \leq  \frac{r(2\rho^*)^{2\kappa}}{A} = \frac{C {\rm Rad}(B)2 \rho^* }{ \sqrt{N}},
 \end{align*}
   where ${\bf C}$ denotes positive constants that might change from one instance to the other and depend only on $A$, $b$, $\kappa$ and $r(\cdot)$ is the function introduced in Definition~\ref{def:function_r_boundedness_simple}.
\end{theo}

In the next Sections~\ref{sec:the_logistic_lasso}, \ref{sec:hinge_} and \ref{sec:any_lipschitz_loss_rkhs_norm} we compute $r(\rho)$ either in the subgaussian setup or in the bounded setup and solve the sparsity equation in various examples, showing the versatility of the main strategy.

\section{Application to logistic LASSO and logistic SLOPE} 
\label{sec:the_logistic_lasso}
The first example of application of the main results in Section~\ref{sec:theoretical_results} involves one very popular method  developed during the last two decades in binary classification which is the Logistic LASSO procedure (cf. \cite{MR2699823,MR2412631,MR2427362,MR2721710,MR3073790}).    

We consider the \emph{vector} framework, where $(X_1, Y_1),\ldots, (X_N, Y_N)$ are $N$ i.i.d. pairs with values in $\R^p\times \{-1, 1\}$ distributed like $(X, Y)$. Both bounded and subgaussian framework can be analyzed in this example. For the sake of shortness and since an example in the bounded case is provided in the next section, only the subgaussian case is considered here and we leave the bounded case to the interested reader. We therefore shall apply Theorem~\ref{thm:main_subgaussian_simple} to get estimation and prediction bounds for the well known logistic LASSO and the new logistic SLOPE.

In this section, we consider the class of linear functional indexed by $R B_{l_2}$ for some radius $R \geq 1$ and the logistic loss:
\begin{equation*}
F=\left\{\inr{\cdot, t} : t\in R B_{l_2}\right\}, \ell_f(x, y) = \log(1+\exp(-yf(x))).
\end{equation*}
As usual the oracle is denoted by $f^*= \argmin_{f\in F} \bE \ell_f(X,Y)$, we also introduce $t^*$ such that $f^* = \inr{\cdot,t^*} $.

\subsection{Logistic LASSO}

The logistic loss function is Lipschitz with constant $1$, so Assumption~\ref{ass:lipschitz} is satisfied. It follows from Proposition~\ref{prop:bernstein_logistic_subgaussian} in Section~\ref{sub:logistic_loss} that Assumption~\ref{ass:margin} is satisfied when the design $X$ is the standard Gaussian variable in $\R^p$ and the considered class $F$. In that case, the Bernstein parameter is $\kappa=1$, and we have $A=c_0/R^3$ for some absolute constant $c_0>0$ which can be deduced from the proof of Proposition~\ref{prop:bernstein_logistic_subgaussian}. We consider the $l_1$ norm $\norm{\inr{\cdot, t}} = \norm{t}_{l_1}$ for regularization. We will therefore obtain statistical results for the RERM estimator $\widehat{f_L} = \inr{\widehat{t_L},\cdot}$  that is defined by
\begin{equation*}
\widehat{t_L} \in\argmin_{t\in R B_{l_2}}\left(\frac{1}{N}\sum_{i=1}^N \log\left(1+\exp(-Y_i \inr{X_i, t}\right) + \lambda \norm{t}_{l_1}\right)
\end{equation*}
where $\lambda$ is a regularization parameter to be chosen according to Theorem~\ref{thm:main_subgaussian_simple}. 

 The two final ingredients needed to apply Theorem~\ref{thm:main_subgaussian_simple} are 1) the computation of the Gaussian mean width of the unit ball $B_{l_1}$ of the regularization function $\norm{\cdot}_{l_1}$ 2) find a solution $\rho^*$ to the sparsity equation~\eqref{eq:sparsity_equation}. 

Let us first deal with the complexity parameter of the problem. If one assumes that the design vector $X$ is \textbf{isotropic}, i.e. $\E \inr{X, t}^2 = \norm{t}_{l_2}^2$ for every $t\in\R^p$ then the metric naturally associated with $X$ is the canonical $l_2$-distance in $\R^p$. In that case, it is straightforward to check that $w(B_{l_1}) \leq c_1 \sqrt{\log p}$ for some (known) absolute constant $c_1>0$ and so we define, for all $\rho\geq0$,
\begin{equation}\label{eq:comp_param_logistic}
r(\rho) = {\bf C}  \left(\rho \sqrt{\frac{\log p}{N}}\right)^{1/2}
\end{equation}
for the complexity parameter of the problem (from now and until the end of Section~\ref{sec:the_logistic_lasso}, the constants ${\bf C}$ depends only on $L$, $C$, $c_0$ and $c_1$).

 Now let us turn to a solution $\rho^*$ of the sparsity equation~\eqref{eq:sparsity_equation}. First note that when the design is isotropic the sparsity parameter is the function
 \begin{equation*}
 \Delta(\rho) = \inf \left\lbrace \sup_{g\in\Gamma_{t^*}(\rho)}\inr{h,g} : h \in \rho S_{l_1} \cap r(2\rho)B_{l_2} \right\rbrace 
 \end{equation*}where $\Gamma_{t^*}(\rho) = t^*+(\rho/20)B_{l_1}$. 
 
A first solution to the sparsity equation is $\rho^* = 20 \norm{t^*}_{l_1}$ because it leads to $0 \in \Gamma_{t^*}(\rho^*)$. This solution is called \emph{norm dependent}. 
 
Another radius $\rho^*$ solution to the sparsity equation~\eqref{eq:sparsity_equation} is obtained when $t^*$ is close to a sparse-vector, that is a vector with a small support. We denote by $\norm{v}_0 :=|{\rm supp}(v)|$ the size of the support of $v\in\R^p$. Now, we recall a result from \cite{LM_sparsity}.
\begin{lemma}[Lemma~4.2 in  \cite{LM_sparsity}]
If there exists some  $v\in t^* +(\rho/20)B_{l_1}$ such that $\norm{v}_0\leq c_0(\rho/r(\rho))^2$ then $\Delta(\rho)\geq 4 \rho/5$ where $c_0$ is an absolute constant.
\end{lemma}In particular, we get that $\rho^*\sim s \sqrt{(\log p)/N}$ is a solution to the sparsity equation if there is a $s$-sparse vector which is $(\rho^*/20)$-close to $t^*$ in $l_1$. This radius leads to the so-called \emph{sparsity dependent} bounds.

After the derivation of the Bernstein parameter $\kappa=1$, the complexity $w(B)$ and a solution $\rho^*$ to the sparsity equation,  we are now in a position to apply Theorem~\ref{thm:main_subgaussian_simple} to get statistical bounds for the Logistic LASSO.

\begin{theo}\label{theo:logistic_lasso}
Assume that $X$ is a standard Gaussian vector in $\R^p$. Let $s\in\{1, \ldots, p\}$. Assume that there exists a $s$-sparse vector in $t^*+{\bf C} s \sqrt{(\log p)/N}B_{l_1}$. Then, with probability larger than $1-{\bf C}\exp\left(-{\bf C} s\log p\right)$, for every $1\leq q\leq 2$, the logistic LASSO estimator $\widehat{t_L}$ with regularization parameter
\begin{equation*}
\lambda = \frac{5 c_1 C L}{8} \sqrt{\frac{\log p}{N}}
\end{equation*}
satisfies
\begin{equation*}
\norm{\widehat{t_L}-t^*}_{l_q}\leq {\bf C} \min\left(s^{1/q}\sqrt{\frac{\log p}{N}}, \norm{t^*}_{l_1}^{1/q}\left(\frac{\log p}{N}\right)^{\frac{1}{2} - \frac{1}{2q}}\right)
\end{equation*}and the excess logistic risk of $\widehat{t_L}$ is such that
\begin{equation*}
\cE_{logistic}= R(\widehat{t_L}) -  R(t^*) \leq {\bf C} \min\left(\frac{s \log(p)}{N}, \norm{t^*}_{l_1}\sqrt{\frac{\log(p)}{N}}\right).
\end{equation*}
\end{theo}
Note that an estimation result for any $l_q$-norm for $1\leq q\leq2$ follows from results in $l_1$ and $l_2$ and the interpolation inequality $\norm{v}_{l_q} \leq \norm{v}_{l_1}^{-1+2/q} \norm{v}_{l_2}^{2-2/q}$.

Estimation results for the logistic LASSO estimator in the generalized linear model have been obtained in \cite{van2008high} under the assumption that the basis functions and the oracle are bounded. This assumption does not hold here since the \emph{basis functions} -- defined here by $\psi_k(\cdot)=\inr{e_k, \cdot}$ where $(e_k)_{k=1}^d$ is the canonical basis of $\R^p$ -- are not bounded when the design is $X\sim\cN(0, I_{d\times p})$. Moreover, we do not make the assumption that $f^*$ is bounded in $L_\infty$. Nevertheless, we recover the same estimation result for the $l_2$-loss and $l_1$-loss as in \cite{van2008high}. But we also provide a prediction result since an excess risk bound is also given in Theorem~\ref{theo:logistic_lasso}.  

Note that Theorem~\ref{theo:logistic_lasso} recovers the classic rates of convergence for the logistic LASSO estimator that have been obtained in the literature so far. This rates is the minimax rate as long as $\log(p/s)$ behaves like $\log p$. This is indeed the case when $s\ll p$ which is the classic setup in high-dimensional statistics. But when $s$ is proportional to $p$ this rate is not minimax since there is a logarithmic loss. To overcome this issue we introduce a new estimator: the logistic SLOPE.

\subsection{Logistic Slope}

The construction of the logistic Slope is similar to the one of the logistic LASSO except that the regularization norm used in this case is the SLOPE norm (cf. \cite{MR3485953,MR3418717}): for every $t=(t_j)\in\R^p$,
\begin{equation}\label{eq:slope_norm}
\norm{t}_{SLOPE} = \sum_{j=1}^p \sqrt{\log(ep/j)}t_j^\sharp
\end{equation}
where $t_1^\sharp\geq t_2^\sharp\geq \cdots\geq 0$ is the non-increasing rearrangement of the absolute values of the coordinates of $t$ and $e$ is the base of the natural logarithm. Using this estimator with a regularization parameter $\lambda\sim 1/\sqrt{N}$ we recover the same result as for the Logistic LASSO case except that one can get, in that case, the optimal minimax rate for any $s\in\{1, \ldots, p\}$:
\begin{equation*}
\sqrt{\frac{s}{N}\log\left(\frac{ep}{s}\right)}.
\end{equation*}

Indeed, it follows from Lemma~5.3 in \cite{LM_sparsity} that the Gaussian mean width of the unit ball $B_{SLOPE}$ associated with the SLOPE norm is of the order of a constant. The \emph{sparsity dependent} radius satisfies 
\begin{equation}\label{eq:rho_star_slope}
\rho^*\sim \frac{s}{\sqrt{N}}\log\left(\frac{ep}{s}\right)
\end{equation}as long as there is a $s$-sparse vector in $t^*+(\rho^*/20) B_{SLOPE}$. The \emph{norm dependent} radius is as usual of order $\norm{t^*}_{SLOPE}$. Then, the next result follows from Theorem~\ref{thm:main_subgaussian_simple}. It improves the best known bounds on the logistic LASSO.

\begin{theo}\label{theo:logistic_slope}
Assume that $X$ is a standard Gaussian vector in $\R^p$. Let $s\in\{1, \ldots, p\}$. Assume that there exists a $s$-sparse vector in $t^*+(\rho^*/20)B_{SLOPE}$ for $\rho^*$ as in \eqref{eq:rho_star_slope}. Then, with probability larger than $1-{\bf C} \exp\left(-{\bf C} s\log(p/s)\right)$, the logistic SLOPE estimator 
\begin{equation*}
\widehat{t_S} \in\argmin_{t\in R B_{l_2}}\left(\frac{1}{N}\sum_{i=1}^N \log\left(1+\exp(-Y_i \inr{X_i, t}\right) + \frac{{\bf C}}{\sqrt{N}} \norm{t}_{SLOPE}\right)
\end{equation*}
satisfies
\begin{equation*}
 \norm{\widehat{t_S}-t^*}_{SLOPE}\leq {\bf C} \min\left(\frac{s}{\sqrt{N}}\log\left(\frac{ep}{s}\right), \norm{t^*}_{SLOPE}\right)
 \end{equation*} and
\begin{equation*}
\norm{\widehat{t_S}-t^*}_{l_2}\leq {\bf C} \min\left(\sqrt{\frac{s}{N}\log\left(\frac{ep}{s}\right)}, \sqrt{\frac{\norm{t^*}_{SLOPE}}{\sqrt{N}}}\right)
\end{equation*}and the excess logistic risk of $\widehat{t_S}$ is such that
\begin{equation*}
\cE_{logistic}(\widehat{t_S}) =  R(\widehat{t_S}) - R(t^*) \leq {\bf C} \min\left(\frac{s \log ep/s}{N}, \norm{t^*}_{l_1}\sqrt{\frac{\log ep/s}{N}}\right).
\end{equation*}
\end{theo}

Let us comment on Theorem~\ref{theo:logistic_slope} together with the fact that we do not make any assumption on the output $Y$ all along this work. Theorem~\ref{theo:logistic_slope} proves that there exists an estimator achieving the minimax rate $s\log(ep/s)/N$ for the $\ell_2$-estimation risk (to the square) with absolutely no assumption on the output $Y$. In the case where a statistical model $Y={\rm sign}(\inr{X, t^*} + \xi)$ holds, where $\xi$ is independent of $X$ then Theorem~\ref{theo:logistic_slope} shows that the RERM with logistic loss and SLOPE regularization achieves the minimax rate $s\log(ep/s)/N$ under no assumption on the noise $\xi$. In particular, $\xi$ does not need to have any moment and, for instance, the mimimax rate $s\log(ep/s)/N$ can still be achieved when the noise has a Cauchy distribution. Moreover, this estimation rate holds with exponentially large probability as if the noise had a Gaussian distribution (cf. \cite{LM13}). This is a remarkable feature of Lipschitz loss functions genuinely understood in Huber's seminal paper \cite{MR0161415}.

\begin{table}[h!]
\centering
\begin{tabular}{l|cc}
 & LASSO & SLOPE \\
 \hline
 $w(B)$ & $\sqrt{\log p}$&$1$ \\[.2cm] 
 $\rho^*$ & $ \displaystyle \frac{s}{\sqrt{N}}\sqrt{\log p}$ & $ \displaystyle \frac{s}{\sqrt{N}} \log \frac{ep}{s}$ \\[.4cm]
 $ r(\rho^*)$ & $ \displaystyle \frac{s}{N}\log p$ & $ \displaystyle \frac{s}{N}\log\frac{ep}{s}$
\end{tabular}
\caption{\label{tab:comparison_lasso_slope} Comparison of the key quantities involved in our study for the $\ell_1$ (LASSO) and SLOPE norms}
\end{table} 

\begin{figure}[!h]
\centering
\begin{tikzpicture}[scale=0.3]
\filldraw[fill=blue!20, draw=black] (-25,0) -- (-15,10) -- (-5,0) -- (-15,-10) -- cycle;
\filldraw[fill=blue!20, draw=black] (5,0) -- (9, 6) --  (15, 10) -- (21, 6) -- (25,0) -- (21, -6) -- (15, -10) -- (9, -6) -- cycle;
\draw (-7,0) node {$t^*$};
\draw (23,0) node {$t^*$};
\filldraw (25,0) circle (0.2cm);
\filldraw (-5,0) circle (0.2cm);
\draw (-16,0) node {$0$};
\draw (14,0) node {$0$};
\filldraw (15,0) circle (0.2cm);
\filldraw (-15,0) circle (0.2cm);
\draw[->, thick] (-15, 0) -- (-6, 8);
\draw (-5,8) node {$G$};
\draw[->, thick] (15, 0) -- (24, 8);
\draw (25,8) node {$G$};
\draw[draw=black] (-8,3) -- (-8,-3) -- (-2,-3) -- (-2,3) -- cycle;
\draw[->] (-5,0) -- (-2,-3);
\draw[->] (-5,0) -- (-2,-2);
\draw[->] (-5,0) -- (-2,-1);
\draw[->] (-5,0) -- (-2,0);
\draw[->] (-5,0) -- (-2,1);
\draw[->] (-5,0) -- (-2,2);
\draw[->] (-5,0) -- (-2,3);
\draw (-5, -4) node {$\partial\norm{\cdot}_1(t^*)$};
\draw[draw=black] (22,2) -- (23, 3) --  (27, 3) -- (28, 2) -- (28,-2) -- (27, -3) -- (23, -3) -- (22, -2) -- cycle;
\draw[->] (25,0) -- (28,-2);
\draw[->] (25,0) -- (28,-1);
\draw[->] (25,0) -- (28,0);
\draw[->] (25,0) -- (28,1);
\draw[->] (25,0) -- (28,2);
\draw (27, -4) node {$\partial\norm{\cdot}_{SLOPE}(t^*)$};
\draw[style = dashed] (0,-10) -- (0,10);
\end{tikzpicture}
\caption{\textbf{Gaussian complexity and size of the sub-differential for the $\ell_1$ ans SLOPE norms:} A ``large'' sub-differential at sparse vectors and a small Gaussian mean width of the unit ball of the regularization norm is better for sparse recovery. In this figure, $G$ represents a ``typical'' Gaussian vector used to compute the Gaussian mean width of the unit regularization norm ball.}
\label{fig:trade_off_subgradient_comp}
\end{figure}
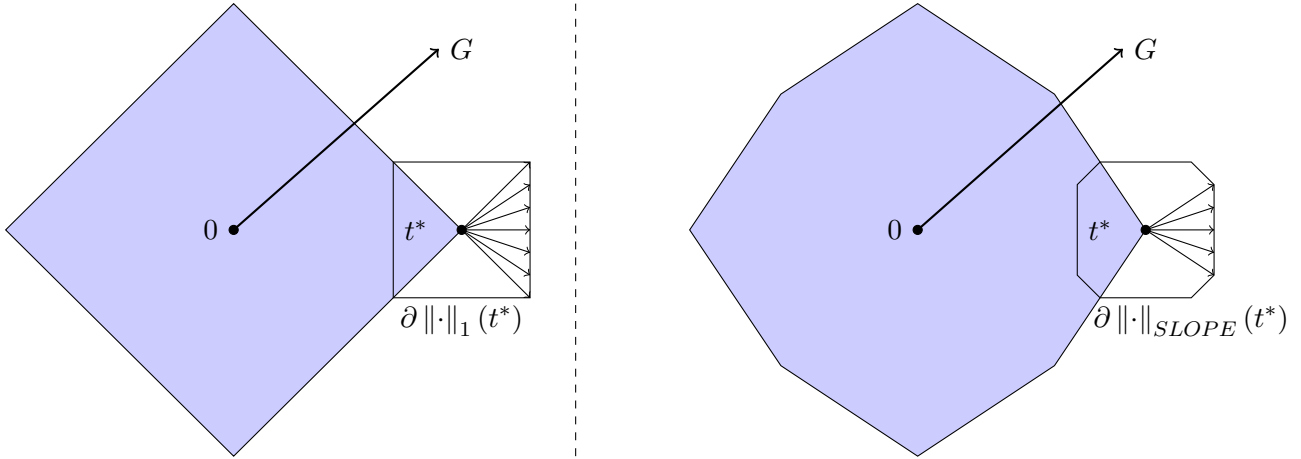

In Table~\ref{tab:comparison_lasso_slope}, the different quantities playing an important role in our analysis have been collected for the $\ell_1$ and SLOPE norms: the Gaussian mean width $w(B)$ of the unit ball $B$ of the regularization norm, a radius $\rho^*$ satisfying the sparsity equation and finally the $L_2$ estimation rate of convergence $r(\rho^*)$ summarizing the two quantities. As mentioned in Figure~\ref{fig:trade_off_subgradient_comp}, having a large sub-differential at sparse vectors and a small Gaussian mean-width $w(B)$ is a good way to construct ``sparsity inducing''regularization norms as it is, for instance the case of ``atomic norms'' (cf. \cite{MR2989474}).   

\section{Application to matrix completion via $S_1$-regularization} 
\label{sec:hinge_}
  
The second example involves matrix completion and uses the bounded setting from Section~\ref{sec:bounded}. The goal is to derive new results on two ways: the 1-bit matrix completion problem where entries are binary, and the quantile completion problem. The main theorems in this section yield upper bounds on completion in $S_p$ norms ($1\leq p\leq 2$) and on various excess risks. We also propose algorithms in order to compute efficiently the RERM in the matrix completion issue but with non differentiable loss and provide a simulation study. We first present a general theorem and then turn to specific loss functions because they induce a discussion about the Bernstein assumption and the $\kappa$ parameter and lead to more particular theorems. 
  
\subsection{General result}
In this section, we consider the matrix completion problem. Contrary to the introduction, we do not immediately focus on the case $Y \in\{-1,+1\}$. So for the moment, $Y$ is a general real random variable and $\ell$ is any Lipschitz loss. The class is $F= \left\{\inr{\cdot, M}: M\in bB_\infty \right\}$, where $bB_\infty = \{M=(M_{pq}):\max_{p,q}|M_{pq}|\leq b\}$ and $b>0$. As the design $X$ takes its values in the canonical basis of $\R^{m\times T}$, the boundedness assumption is  satisfied. Apart from that, the notations and assumptions are  as in the introduction, that is, we assume that $X$ satisfies Assumption~\ref{ass:matrix_design}, with parameters $(\underbar{c},\bar{c})$, and the penalty is the nuclear norm. Thus, the RERM is given by
\begin{equation}
\label{eq:estimator-matrix}
\widehat{M} \in\argmin_{M\in b B_\infty}\left(\frac{1}{N}\sum_{i=1}^N \ell(\inr{X_i,M},Y_i) + \lambda \norm{M}_{S_1}\right).
\end{equation}
Statistical properties of \eqref{eq:estimator-matrix} will follow from Theorem~\ref{thm:main_bounded_simple} since one can recast this problem in the setup of Section~\ref{sec:bounded}. The oracle matrix $M^*$ is defined by $f^*=\left<\cdot,M^*\right>$, that is, $M^* = \argmin_{M \in bB_\infty} \bE \ell(\inr{M,X},Y)$.

Let us also introduce the matrix $\overline{M}=\argmin_{M\in\mathbb{R}^{m\times T}} \bE \ell(\inr{M,X},Y)$. Note that $\left<\overline{M},\cdot\right>=\overline{f}=\arg\min_{f \text{ measurable}} \bE \ell(f(X),Y)$. Our general results usually are on $f^*$ rather than on $\overline{f}$ as it is usually impossible to provide rates on the estimation of $\overline{f}$ without stringent assumptions on $Y$ and $F$. However, as noted in the introduction, in 1-bit matrix completion with the hinge loss, we have $\overline{M}=M^*$ without any extra assumption when $b=1$ (this is a favorable case). On the other hand, to get fast rates in matrix completion with quantile loss requires that $\overline{M}=M^*$ (which is a stringent assumption in this setting).

\paragraph{Complexity function} We first compute the complexity parameter $r(\cdot)$ as introduced in Definition~\ref{def:function_r_boundedness_simple}. To that end one just needs to compute the global Rademacher complexity of the unit ball of the regularization function which is $B_{S_1} = \{A\in\R^{m\times T}:\norm{A}_{S_1}\leq1\}$:  
\begin{equation}
\label{eq:rademacher-matrices}
{\rm Rad}(B_{S_1}) = \E \sup_{\norm{A}_{S_1}\leq1}\Big|\frac{1}{\sqrt{N}}\sum_{i=1}^N \eps_i \inr{X_i, A}\Big| = \E \norm{\frac{1}{\sqrt{N}}\sum_{i=1}^N \eps_i X_i}_{S_\infty}\leq 
c_0(\underbar{c},\bar{c}) \sqrt{\frac{\log(m+T)}{\min(m,T)}}
\end{equation}
where $\norm{\cdot}_{S_\infty}$ is the operator norm (i.e. the largest singular value), the last inequality follows from Lemma~1 in \cite{MR2906869} and $c_0(\underbar{c},\bar{c})>0$ is some constant that depends only on $\underbar{c}$ and $\bar{c}$.

The complexity parameter $r(\cdot)$ is derived from Definition~\ref{def:function_r_boundedness_simple}: for any $\rho\geq0$, 
\begin{equation}\label{eq:choice_r_hinge_S1}
r(\rho) = \left[ \frac{C A \rho {\rm Rad}(B_{S_1})}{\sqrt{N}} \right]^{\frac{1}{2\kappa}}
     = {\bf C} \left[ \rho  \sqrt{\frac{\log(m+T)}{N \min(m,T)}} \right]^{\frac{1}{2\kappa}}
\end{equation}
where from now the constants ${\bf C}$ depend only on $\underbar{c}$, $\bar{c}$, $b$, $A$ and $\kappa$.

\paragraph{Sparsity parameter} The next important quantity is the sparsity parameter. Its expression in this particular case is, for any $\rho>0$,
\begin{equation*}
\Delta(\rho) = \inf\left\{\sup_{G \in\Gamma_{M^*}(\rho)}\inr{H,G}: H\in \rho S_{S_1}\cap ((\sqrt{mT}/\underbar{c}) r(2\rho))B_{S_2}\right\}
\end{equation*}
where $\Gamma_{M^*}(\rho)$ is the union of all the sub-differential of $\norm{\cdot}_{S_1}$ in a $S_1$-ball of radius $\rho/20$ centered in $M^*$. Note that the normalization factor $\sqrt{mT}$ in the localization $(\sqrt{mT} r(2\rho))B_{S_2}$ comes from the ``non normalized isotropic'' property of $X$: $ \underbar{c}\norm{M}_{S_2}^2/(mT) \leq \E \inr{X,M}^2\leq \bar{c} \norm{M}_{S_2}^2/(mT)$ for all $M\in\R^{m\times T}$. Now, we use a result from \cite{LM_sparsity} to find a solution to the sparsity equation.

\begin{lemma}[Lemma~4.4 in \cite{LM_sparsity}]\label{lem:sub_diff_S1}
There exists an absolute constant $c_1>0$ for which the following holds. If there exists $V\in M^*+(\rho/20)B_{S_1}$ such that $
{\rm rank}(V)\leq \left(c_1 \rho/(\sqrt{mT}r(\rho))\right)^2$ then $\Delta(\rho)\geq 4\rho/5$.
\end{lemma}

It follows from Lemma~\ref{lem:sub_diff_S1} that the sparsity equation~\eqref{eq:sparsity_equation} is satisfied by $\rho^*$ when it exists $V\in M^*+(\rho^*/20)B_{S_1}$ such that  $\rank(V) =  c_1 \left(\rho^*/(\sqrt{mT}r(\rho^*))\right)^2$. Note obviously that $V$ can be $M^*$ itself, in this case, $\rho^*$ can be taken such that $\rank(M^*) = c_1 \left(\rho^*/(\sqrt{mT}r(\rho^*))\right)^2$. However, when $M^*$ is not low-rank, it might still be that a low-rank approximation $V$ of $M^*$ is close enough to $M^*$ w.r.t. the $S_1$-norm. As a consequence, if for some $s\in\{1, \ldots, \min(m,T)\}$ there exists a matrix $V$ with rank at most $s$ in $M^*+(\rho_s^*/20)B_{S_1}$ where
  \begin{equation}\label{eq:choice_rho_star}
   \rho_s^* = {\bf C}  \left( s m T \right)^{\frac{\kappa}{2\kappa-1}}
         \left(\frac{\log(m+T)}{N \min(m,T)}\right)^{\frac{1}{2(2\kappa-1)}}.
   \end{equation}
then $\rho_s^*$ satisfies the sparsity equation. 

Following the remark at the end of Subsection~\ref{sec:sparsity_param}, another possible choice is $\rho^* = 20 \|M^*\|_{S_1}$ in order to get \emph{norm dependent} rates. In the end, we choose $\rho^* = {\bf C} \min \left[ \rho_s^* , \|M^*\|_{S_1} \right]$.
We are now in a position to apply Theorem~\ref{thm:main_bounded_simple} to derive statistical properties for the RERM $\widehat{M}$ defined in \eqref{eq:estimator-matrix}.

\begin{theo}\label{theo:S1_isotrope}
Assume that Assumption~\ref{ass:lipschitz}, \ref{ass:matrix_design} and \ref{ass:margin} hold. Consider the estimator in~\eqref{eq:estimator-matrix} with regularization parameter
\begin{equation}
\label{eq:lambda_matrix}
\lambda = \frac{c_0(\underbar{c},\bar{c})720}{7} \sqrt{\frac{\log(m+T)}{N \min(m,T)}}
\end{equation}
where $c_0(\underbar{c},\bar{c})$ are the constants in Assumption~\ref{ass:matrix_design}. Let $s\in\{1, \ldots, \min(m, T)\}$ and assume that there exists a matrix with rank at most $s$ in $M^*+(\rho^*_s/20)B_{S_1}$. Then, with probability at least 
\begin{equation*}
 1- {\bf C} \exp\left(-{\bf C}s(m+T) \log(m+T) \right)
 \end{equation*}
 we have
 \begin{align*}
 \norm{\widehat{M} - M^*}_{S_1}& \leq  {\bf C} \min \left\{ \left( s m T \right)^{\frac{\kappa}{2\kappa-1}}
  \left(\frac{\log(m+T)}{N \min(m,T)}\right)^{\frac{1}{2(2\kappa-1)}},\|M^*\|_{S_1}\right\} ,\\ 
 \frac{1}{\sqrt{mT}}\norm{\widehat{M}-M^*}_{S_2}& \leq {\bf C}  \min\left\{ \left(\frac{s(m+T)\log(m+T)}{N}\right)^{\frac{1}{2({2\kappa}-1)}}, \left(\|M^*\|_{S_1} \sqrt{\frac{\log(m+T)}{N \min(m,T)}}\right)^{\frac{1}{2\kappa}}\right\}
 \\  
 \cE(\widehat{M}) & \leq {\bf C} \min\left\{ \left(\frac{s(m+T)\log(m+T)}{N}\right)^{\frac{{\kappa}}{{2\kappa}-1}},\|M^*\|_{S_1} \sqrt{\frac{\log(m+T)}{N \min(m,T)}}\right\}.
 \end{align*}
\end{theo}
Note that the interpolation inequality also allows to get a bound for the $S_p$ norm, when $1\leq p\leq 2$:
\begin{align*}
 \frac{1}{(mT)^{\frac{1}{p}}} \norm{\widehat{M} - M^*}_{S_p} 
 \leq
  {\bf C}
 \min\Biggl\{&
  \Biggl[
 \left(\frac{  s^{2(p-1)+\kappa(2-p)} (m+T)^{p-1} }{ \min(m,T)^{\frac{2-p}{2}}}\right)^{\frac{1}{p}}
 \sqrt{\frac{\log(m+T)}{N}}\Biggr]^{\frac{1}{2\kappa-1}}
 ,
 \\
 &
  \|M^*\|_{S_1}^{\frac{p-1 + \kappa(2-p)}{p\kappa}}\left( \frac{\log(m+T)}{N\min(m,T)}\right)^{\frac{p-1}{2\kappa p}}
    \left(\frac{1}{mT}\right)^{\frac{2-p}{p}} \Biggr\}.
\end{align*}

Theorem~\ref{theo:S1_isotrope} shows that the sparsity dependent error rate in the excess risk bound is (for $s={\rm rank}(M^*)$)
\begin{equation*}
\left(\frac{\rank(M^*)(m+T)\log(m+T)}{N}\right)^{\frac{\kappa}{2\kappa-1}}
\end{equation*} which is the classic excess risk bound under the margin assumption up to a log factor (cf. \cite{MR2336861}).
As for the $S_2$-estimation error, when $\kappa=1$, we recover the classic $S_2$-estimation rate 
\begin{equation*}
\sqrt{\frac{\rank(M^*)(m+T)\log(m+T)}{N}}
\end{equation*}
which is minimax in general (up to log terms, e.g. take the quadratic loss when $Y$ is bounded and compare to~\cite{rohde2011estimation}).

\subsection{Algorithm and Simulation Outlines}

Since this part provides new methods and results on matrix completion, we propose an algorithm in order to compute efficiently the RERM using the hinge loss and the quantile loss. This section explains the structure of the algorithm that is used with specific loss functions in next sections. Although many algorithms exist for the least squares matrix completion, at our knowledge many of them treat only the exact recovery such as in \cite{cai2010singular} and \cite{mazumder2010spectral}, or at least they all deal with differentiable loss functions, see \cite{hsieh2014nuclear}. On the other hand, the two losses that we mainly consider here are non differentiable because they are piecewise linear (in the case of hinge and $0,5$-quantile loss functions): new algorithms are hence needed. It has been often noted that the RERM with respect to the hinge loss or $0.5$-quantile loss can been solved by a semidefinite programming but the cost is prohibitive for large matrices, say dimensions larger than $100$. It actually works for small matrices as we ran SDP solver in Python in very small examples.

We propose here an \emph{alternating direction method of multiplier} (ADMM) algorithm. For a clear and self-contained introduction to this class of algorithms, the reader is referred to \cite{boyd2011distributed} and we do not explain all the details here and we keep the same vocabulary. When the optimization problem is a sum of two parts, the core idea is to split the problem by introducing an extra variable. In our case, the two following problems are equivalent:
\begin{align*}
 \underset{M}{\textrm{minimize}} \,   \left\lbrace \frac{1}{N}\sum_{i=1}^N \ell(\inr{X_i,M},Y_i) + \lambda \norm{M}_{S_1} \right\rbrace, \quad & \underset{M,L}{\textrm{minimize}} \, \left\lbrace\frac{1}{N}\sum_{i=1}^N \ell(\inr{X_i,M},Y_i) + \lambda \norm{L}_{S_1} \right\rbrace \\
& \textrm{subject to } M=L    
\end{align*}

Below, we use the scaled form and the $m\times T$ matrix $U$ is then called the \emph{scaled dual variable}. Note that the $S_2$ norm is also the Froebenius norm and is thus elementwise. We can now exhibit the \emph{augmented Lagrangian}:
\begin{equation*}
L_\alpha(M,L,U) = \frac{1}{N}\sum_{i=1}^N \ell(\inr{X_i,M},Y_i) + \lambda \norm{L}_{S_1}  + \frac{\alpha}{2}\norm{M-L+ U}_{S_2}^2 - \frac{\alpha}{2}\norm{U}_{S_2}^2,
\end{equation*}
where $\alpha$ is a positive constant, called the \emph{augmented Lagrange parameter}. The ADMM algorithm \cite{boyd2011distributed} is then:
\begin{align}
M^{k+1} &= \underset{M}{\textrm{argmin }} \left( \frac{1}{N}\sum_{i=1}^N \ell(\inr{X_i,M},Y_i) + \frac{\alpha}{2} \norm{M-L^k + U^k}_{S_2}^2 \right)  \label{eq:Admm1_gen}\\
L^{k+1} &= \underset{L}{\textrm{argmin }} \left( \lambda \norm{L}_{S_1} + \frac{\alpha}{2} \norm{M^{k+1}-L + U^k}_{S_2}^2 \right) \label{eq:Admm2_gen} \\
U^{k+1} &= U^k + M^{k+1}-L^{k+1} \nonumber
\end{align}

The starting point $(M^0, L^0, U^0)$ uses one random matrix with independent Gaussian entries for $M^0$ and two zero matrices for $L^0$ and $U^0$. Another choice of starting point is to use a previous estimator with a larger $\lambda$. The stopping criterion is, as explained in \cite{boyd2011distributed}, $\norm{M^{k+1} - M^k}_{S_2}^2 + \norm{U^{k+1} - U^k}_{S_2}^2 \leq \varepsilon$ for a fixed threshold $\varepsilon$. It means that it stops when both $(U_k)$ and $(M_k)$ start converging. 

\paragraph{General considerations} The second step \eqref{eq:Admm2_gen} is independent of the loss function. It is well-known that the solution of this problem is $S_{\lambda/\alpha}(M^{k+1}+ U^k)$ when $S_a(M)$ is the soft-thresholding operator with magnitude $a$ applied to the singular values of the matrix $M$. It is defined for a rank $r$ matrix M with SVD $M = U\Sigma V^\top$ where $\Sigma = \textrm{diag}\left( \left(d_i\right)_{1\leq i \leq r} \right)$ by $S_a(M) = U S_a(\Sigma) V^\top$ where $S_a(\Sigma) = \textrm{diag}\left( \left(\max(0, d_i-a)\right)_{1\leq i \leq r} \right)$.

It requires the SVD of a $m\times T$ matrix at each iteration and is the main bottleneck of this algorithm (the other main step \eqref{eq:Admm1_gen} can be performed elementwise since the $X_i$'s take their values in the canonical basis of $\R^{m\times T}$; so it needs only at most $N$ operations). Two methods may be used in order to speed up the algorithm: efficient algorithms for computing the $n$ largest singular values and the associate subspaces, such as the well-known PROPACK routine in Fortran. It can be plugged in order to solve \eqref{eq:Admm2_gen} by computing the $n$ largest and stop at this stage if the lowest computed singular values is lower than the threshold. It is obviously more relevant when the target is expected to have a very small rank. This method has been implemented in Python and works well in practice even though the parameter $n$ has to be tuned carefully. An alternative method is to use approximate SVD such as in \cite{halko2011finding}.    

Moreover, the first step \eqref{eq:Admm1_gen} (which may be performed elementwise) has a closed form solution for hinge and quantile loss: it is a soft-thresholding applied to a specified quantity. 

Simulated observations as well as real-world data (cf. the MovieLens dataset\footnote{available in http://grouplens.org/datasets/movielens/}) are considered in the examples below. Finally note that parameter $\lambda$ is tuned by cross-validation.

\subsection{$1$-bit matrix completion}

In this subsection we assume that $Y\in\{-1,+1\}$, and we challenge two loss functions: the logistic loss, and the hinge loss. It is worth noting that the minimizer $\overline{M}=\argmin_{M\in \R^{m\times T}} \bE \ell(\inr{M,X},Y)$ is not the same for both losses. For the hinge loss, it is known that it is the matrix formed by the Bayes classifier. This matrix has  entries bounded by $1$ so $M^*=\overline{M}$ as soon as $b=1$. In opposite to this case, the logistic loss leads to a matrix $\overline{M}$ with entries formed by the odds ratio. It may even be infinite when there is no noise. 

\paragraph{Logistic loss.} Let us start by assuming that $\ell$ is the logistic loss. Thanks to Proposition~\ref{prop:margin_logistic_loss} we know that $\kappa=1$ for any $b$ ($A$ is also known, $A= 4\exp(2b)$) and therefore next result follows from Theorem~\ref{theo:S1_isotrope}. Note that we do not assume that $\overline{M}$ is in $F$ and therefore our results provides estimation and prediction bounds for the oracle $M^*$.

\begin{theo}[1-bit Matrix Completion with logistic loss]\label{theo:rates_logistic_S_1}
Assume that Assumption~\ref{ass:matrix_design} holds. Let $s\in\{1, \ldots, \min(m, T)\}$ and assume that there exists a matrix with rank at most $s$ in $M^*+(\rho^*_s/20)B_{S_1}$ where $\rho^*_s$ is defined in \eqref{eq:choice_rho_star}.
With probability at least 
\begin{equation*}
 1- {\bf C} \exp\left(-{\bf C} s \max(m,T) \log(m+T) \right)
 \end{equation*}
 the estimator
\begin{equation}
\widehat{M} \in\argmin_{M\in b B_\infty}\left(\frac{1}{N}\sum_{i=1}^N \log\left( 1+\exp\left(-Y_i\inr{X_i,M}\right)\right) + \lambda \norm{M}_{S_1}\right)
\end{equation}
with $\lambda$ as in Equation~\eqref{eq:lambda_matrix} satisfies
 \begin{align*}
  \frac{1}{mT} \norm{\widehat{M} - M^*}_{S_1}& \leq  {\bf C} \min \left\{ s
  \sqrt{\frac{\log(m+T)}{N \min(m,T)}}, \frac{\norm{M^*}_{S_1}}{mT} \right\} ,\\ 
 \frac{1}{\sqrt{mT}}\norm{\widehat{M}-M^*}_{S_2}& \leq {\bf C} \min\left\{ \sqrt{\frac{s\max(m,T) \log(m+T)}{N} }, \norm{M^*}_{S_1}^{\frac{1}{2}} \left(\frac{\log(m+T)}{N \min(m,T)}\right)^{\frac{1}{4}} \right\} \\  
 \cE_{logistic}(\widehat{M}) & \leq {\bf C} \min \left\{  \frac{s\max(m,T)\log(m+T)}{N},\norm{M^*}_{S_1} \sqrt{\frac{\log(m+T)}{N \min(m,T)}}\right\}.
 \end{align*}
\end{theo}
Using an interpolation inequality, it is easy to derive estimation bound in $S_p$ for all $1\leq p\leq 2$ as in Theorem~\ref{theo:S1_isotrope} so we do not reproduce it here. Also, note that our bound on $\norm{\widehat{M}-M^*}_{S_2}$ is of the same order as the one in~\cite{lafond2014probabilistic}. We actually now prove that this rate is minimax-optimal (up to log terms).
\begin{theo}[Lower bound with logistic loss]
\label{lower_bound_logistic_S1}
For a given matrix $M\in B_\infty$, define $\mathbb{P}_{M}^{\otimes N}$ as the probability distribution of the $N$-uplet $(X_i,Y_i)_{i=1}^N$ of i.i.d. pairs distributed like $(X, Y)$ such that $X$ is uniformly distributed on the canonical basis $(E_{p,q})$ of $\R^{m\times T}$ and $\mathbb{P}_M(Y=1|X=E_{p,q})=\exp(M_{pq})/[1+\exp(M_{pq})]$ for every $(p,q)\in\{1, \ldots, m\}\times \{1, \ldots,T\}$. Fix $s\in\{1,\dots,\min(m,T)\}$ and assume that $N\geq s (m+T) \log(2)/(8b^2)$. Then
 \begin{align*}
   \inf_{\widehat{M}} \sup_{
 \tiny{
 \begin{array}{c}
 M^*\in b B_\infty
 \\
 {\rm rank}(M^*) \leq s
 \end{array}
 }}
  \mathbb{P}_{M^*}^{\otimes N}\left(
  \frac{1}{\sqrt{mT}}\norm{\widehat{M}-M^*}_{S_2}
  \geq c \sqrt{ \frac{(m+T) s }{N}}
  \right)  \geq \beta 
\end{align*}
  for some universal constants $\beta,c>0$.
\end{theo}

Also, as pointed out in the introduction, the quantity of interest is not the logistic excess risk, but the classification excess risk: let us remind that $R_{0/1}(M)= \mathbb{P}[(Y \neq {\rm sign}(\left<M,X\right>) ] $ for all $M\in\R^{m\times T}$. Even if we assume that $M^*=\overline{M}$, all that can be deduced from Theorem 2.1 in~\cite{zhang2004statistical} is that
\begin{equation*}
\cE_{0/1}(\widehat{M}) = R_{0/1}(\widehat{M}) - \inf_{M \in \mathbb{R}^{m\times T} } R_{0/1}(M) \leq {\bf C} \sqrt{\cE_{logistic}(\widehat{M})}  \leq  {\bf C} \sqrt{\frac{\rank(\overline{M})(m+T)\log(m+T)}{N}} .
\end{equation*}
But this rate on the excess $0/1$-risk may be much better under the margin assumption \cite{Mammen1999,Tsy04} (cf. Equation~\eqref{eq:margin_tsybakov} below). This motivates the use of the hinge loss instead of the logistic loss, for which the results in~\cite{zhang2004statistical} do not lead to a loss of a square root in the rate.

\paragraph{Hinge loss.} As explained above, the choice $b=1$ ensures $\overline{M}=M^*$ without additional assumption. Thanks to Proposition~\ref{prop:margin_hinge_loss} we know that as soon as $\inf_{p,q}|\overline{M}_{p,q}-1/2| \geq \tau $ for some $\tau>0$, the Bernstein assumption is satisfied by the hinge loss with $\kappa=1$ and $A=1/(2\tau)$. This assumption seems very mild in many situations and we derive the results with it. 

\begin{theo}[1-bit Matrix Completion with hinge loss]\label{theo:upper_bound_hinge_S1}
Assume that Assumption~\ref{ass:matrix_design} holds. Assume that $\inf_{p,q}|P(Y=1|X=E_{p,q})-1/2| \geq \tau $ for some $\tau>0$. Let $s\in\{1, \ldots, \min(m, T)\}$ and assume that there exists a matrix with rank at most $s$ in $\overline{M}+(\rho^*_s/20)B_{S_1}$ where $\rho^*_s$ is defined in \eqref{eq:choice_rho_star}. With probability at least 
\begin{equation*}
 1- {\bf C}  \exp\left(-{\bf C}  s \max(m,T) \log(m+T) \right)
 \end{equation*}
 the estimator
\begin{equation}
\label{hinge_minimization_matrix}
\widehat{M} \in\argmin_{M\in B_\infty}\left(\frac{1}{N}\sum_{i=1}^N \left(1-Y_i\inr{X_i,M}\right)_{+} + \lambda \norm{M}_{S_1}\right)
\end{equation}
with $\lambda$ as in Equation~\eqref{eq:lambda_matrix} satisfies
 \begin{align*}
  \frac{1}{mT} \norm{\widehat{M} - \overline{M}}_{S_1}& \leq  {\bf C}  \min \left\{s
  \sqrt{\frac{\log(m+T)}{N \min(m,T)}}, \frac{\norm{\overline{M}}_{S_1}}{mT} \right\} ,\\ 
 \frac{1}{\sqrt{mT}}\norm{\widehat{M}-\overline{M}}_{S_2}& \leq {\bf C}  \min\left\{ \sqrt{\frac{s(m+T) \log(m+T)}{N} }, \norm{\overline{M}}_{S_1}^{\frac{1}{2}} \left(\frac{\log(m+T)}{N \min(m,T)}\right)^{\frac{1}{4}} \right\} \\  
 \cE_{hinge}(\widehat{M}) & \leq {\bf C}  \min \left\{  \frac{s(m+T)\log(m+T)}{N},\norm{\overline{M}}_{S_1} \sqrt{\frac{\log(m+T)}{N \min(m,T)}}\right\}.
 \end{align*}
\end{theo}

In this case, \cite{zhang2004statistical} implies that the excess risk bound for the classification error (using the $0/1$-loss) is the same as the one for the hinge loss: it is therefore of the order of $\rank(\overline{M})\max(m,T)/N$.

Note that the rate $\rank(\overline{M})\max(m,T)/N$ for the classification excess error was only reached in~\cite{Cottet2016} up to our knowledge (using the PAC-Bayesian technique from~\cite{catoni2004statistical,catoni2007pac,mai2015Bayesian,alquier2015properties}), in the very restrictive noiseless setting - that is, $P(Y=1|X=E_{p,q}) \in\{0,1\}$ which is equivalent to $P(Y={\rm sign}(\inr{\overline{M},X})=1$. Here this rate is proved to hold in the general case. Other works, including~\cite{srebro2004maximum}, obtained only rates in $1/\sqrt{N}$. Finally, we prove that this rate is the minimax rate in the next result.

\begin{theo}[Lower bound with hinge loss]
\label{theo:minimax_lower_bound_hinge_S1}
For a given matrix $M\in B_\infty$, let $\mathbb{E}_{M}^{\otimes N}$ be the expectation w.r.t. the $N$-uplet $(X_i,Y_i)_{i=1}^N$ of i.i.d. pairs distributed like $(X, Y)$ such that $X$ is uniformly distributed on the canonical basis $(E_{p,q})$ of $\R^{m\times T}$ and $\mathbb{P}_M(Y=1|X=E_{p,q})= M_{pq}$ for every $(p,q)\in\{1, \ldots, m\}\times \{1, \ldots, T\}$. Fix $s\in\{1,\dots,\min(m,T)\}$ and assume that $N\geq s \max(m,T) \log(2)/8$. Then
 \begin{equation*}
\inf_{\widehat{M}} \sup_{
 \tiny{
 \begin{array}{c}
 M^*\in B_\infty
 \\
 {\rm rank}(M^*) \leq r
 \end{array}
 }}
  \mathbb{E}_{M^*}^{\otimes N}\left( \cE_{hinge}(\widehat{M})\right)
  \geq c   \frac{s \max(m,T)}{N}
 \end{equation*}  
  for some universal constants $c>0$.
\end{theo}

Theorem~\ref{theo:minimax_lower_bound_hinge_S1} provides a minimax lower bound in expectation whereas Theorem~\ref{theo:upper_bound_hinge_S1} provides an excess risk bound with large deviation. The two residual terms of the excess hinge risk from Theorem~\ref{theo:minimax_lower_bound_hinge_S1} and Theorem~\ref{theo:upper_bound_hinge_S1} match up to the $\log(m+T)$ factor.   

\paragraph{Simulation Study.} As the hinge loss has not been often studied in the matrix context, we provide many simulations in order to show the robustness of our method and the opportunity of using the hinge loss rather than the logistic loss. We follow the simulations ran in \cite{Cottet2016} and compare several methods. An estimator based on the logistic model, studied in \cite{Davenport14}, is also challenged\footnote{In the followings, the four estimators will be referred to \emph{Hinge} for  estimator given in~\eqref{hinge_minimization_matrix}, \emph{Hinge Bayes} and \emph{Logit Bayes} for the two Bayesian estimators from \cite{Cottet2016} with respectively hinge and logistic loss functions, and \emph{Logit} for the estimator from \cite{Davenport14}. The Bayesian estimators use the Gammma prior distribution.}. 

\paragraph{A first set of simulations.}

The simulations are all based on a low-rank $200 \times 200$ matrix $M^\star$ from which the data are generated and which is the target for the predictions. $M^\star$ is also a minimizer of $R_{0/1}$ so the error criterion that we will report for a matrix $M$ is the difference of the predictions between $M^\star$ and $M$, which is $\bP[\textrm{sign}(\inr{M^\star,X}) \neq \textrm{sign}(\inr{M,X})]$. The $X_i$'s correspond to $20 \%$ of the entries randomly picked so the misclassification rate is also $1/mT \sum_{p,q} I\{\textrm{sign}(M_{p,q}) \neq \textrm{sign}(M^\star_{p,q})\}$. 

Two different scenarios are tested: the first one (called A), involves a matrix $M^\star$ with only entries in $\{-1,+1\}$ so the Bayes classifier is low rank and favors the hinge loss. The second test (called B) involves a matrix $M^\star=LR^\top$ where $L,R$ have i.i.d. Gaussian entries and the rank is the number of columns. In this case, the Bayes matrix contains the signs of a low-rank matrix, but it is not itself low rank in general. We also test the impact of the noise structure on the results:
\begin{enumerate}
\item (noiseless) $Y_i = {\rm sign}(\inr{M^\star,X_i})$ 
\item (logistic) $Y_i = {\rm sign}(\inr{M^\star,X_i}+Z_i)$, where $Z_i$ follows a logistic distribution 
\item (switch) $Y_i = \epsilon_i {\rm sign}(\inr{M^\star,X_i})$ where $\epsilon_i = (1-p)\delta_1 + p\delta_{-1}$
\end{enumerate}
Finally, we run all the simulations on rank $3$ and rank $5$ matrices. $\lambda$ is tuned by cross validation. All the simulations are run one time.

\begin{table}[h!]
\centering
\begin{tabular}{ll|ccc|ccc}
& Model & A1 & A2 ($p=.1$) & A3 & B1 & B2 ($p=.1$) & B3 \\
\hline \hline
\multirow{4}{*}{Rank $3$} 	& Hinge 			& 0	& 0		& 14.5	& 6.7	& 10.9 	& 21.0	\\
							& Logit			& 0	& 0.5	& 17.3	& 5.1	& 10.7 	& 19.8  	\\
							& Hinge Bayes 	& 0	& 0.1	& 8.5	& 5.3	& 10.8 	& 22.1	\\  
							& Logit Bayes	& 0	& 0.5	& 16.0	& 4.1	& 10.1 	& 16.0 	\\
\hline
\multirow{4}{*}{Rank $5$} 	& Hinge 			& 0	& 0.8	& 29.0	& 11.7	& 19.3	& 23.3	\\
							& Logit	 		& 0	& 3.1	& 30.1	& 9.0	& 18.3 	& 22.1  	\\
							& Hinge Bayes 	& 0	& 0.5	& 27		& 9.4	& 17.9 	& 24.4	\\  
							& Logit Bayes	& 0	& 4.4	& 32.5	& 7.8	& 17.3 	& 21.5 	
\end{tabular}
\caption{Misclassification error rates on simulated matrices in various cases. Model $\in\{A, B\}\{1, 2, 3\}$ refers to scenario $\in\{A, B\}$ and noise structure $\in\{1, 2, 3\}$. For the noise-free Model $=A0$, the $0$ column shows the exact reconstruction property of all procedures.}
\label{tab:hingesim}
\end{table} 

The results are very similar among the methods, see Table~\ref{tab:hingesim}. The logistic loss performs better for matrices of type B and especially for high level of noise in the logistic data generation as expected. For type A matrices, the hinge loss performs slightly better. The Bayesian models performs as good as the frequentist estimators even though the program solved is not convex.

\begin{figure}[h!]
\centering
\includegraphics[scale=.5]{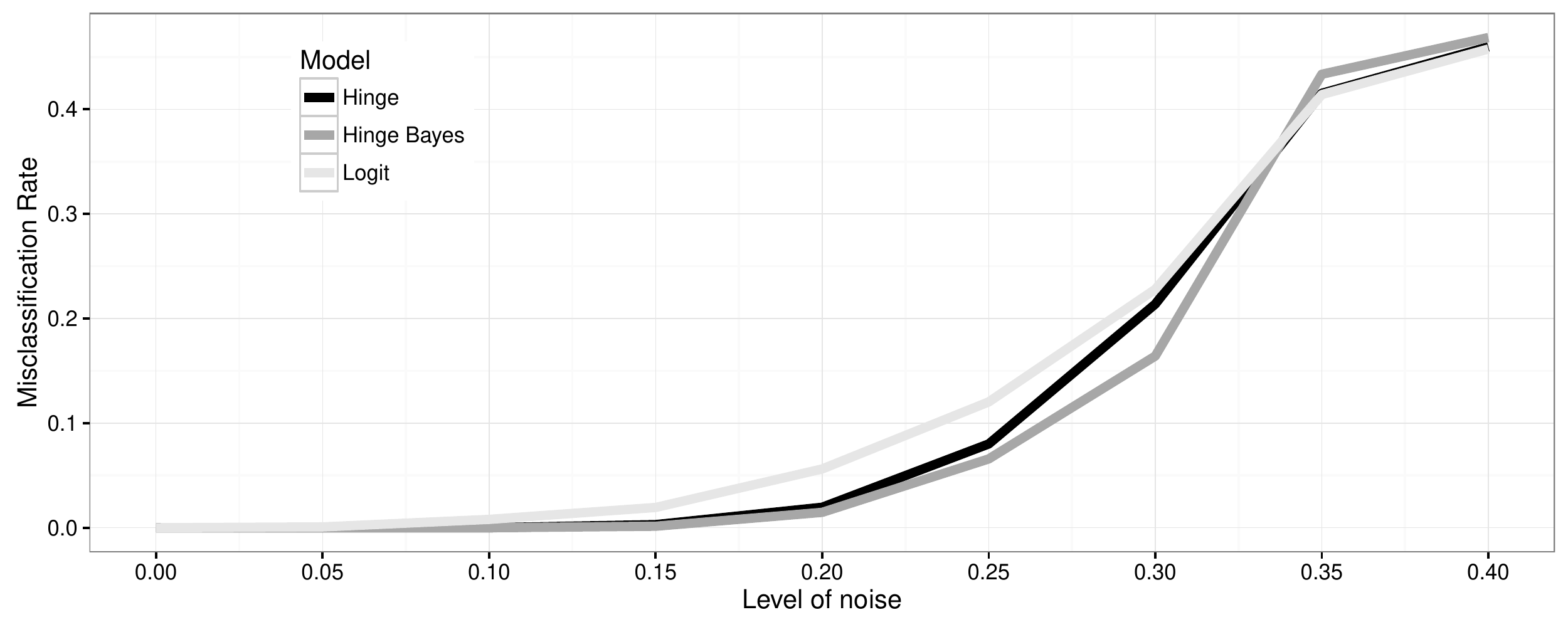}
\label{fig:hinge}
\caption{Misclassification error rates for a large range of switch noise (noise structure number $3$).}
\end{figure}

\paragraph{Impact of the noise level.} The second experiment is a focus on the switch noise and matrices that are well separated (as A2 in the previous example). The noise lies between $p=0$ and almost full noise ($p=.4$). The performance of the RERM with the hinge loss is slightly worse than the Bayesian estimator with hinge loss but always better than the RERM with the logistic loss, see Figure \ref{fig:hinge}.  

\paragraph{Real dataset.} We finally run the hinge loss estimator on the MovieLens dataset. The ratings, that lie in $\{1,2,3,4,5\}$, are split between good ratings ($4,5$) and bad ratings (others). The goal is therefore to predict whether the user will like a movie or not. On a test set that contains $20\%$ of the data, the misclassification rate in prediction are almost the same for all the methods (Table \ref{tab:hingeML}).

\begin{table}[h!]
\centering
\begin{tabular}{l|ccc}
Model & Hinge Bayes & Logit & Hinge \\
\hline
misclassification rate & .28 & .27 & .28 
\end{tabular}
\caption{Misclassification Rate on MovieLens 100K dataset}
\label{tab:hingeML}
\end{table} 
 
\subsection{Quantile loss and median matrix completion}

The matrix completion problem with continuous entries has almost always been tackled with a penalized least squares estimator~\cite{CandesP10,MR2906869,klopp2014noisy,LM_sparsity,mai2015Bayesian}, but the use of other loss functions may be very interesting in this case too. Our last result on matrix completion is a result for the quantile loss $\rho_{\tau}$ for $\tau\in(0,1)$. Let us recall that $\rho_\tau(u) = u(\tau - I(u\leq0))$ for all $u\in\R$ and $\ell_M(x,y) = \rho_\tau(y-\inr{M, x})$. While the aforementionned references provided ways to estimate the conditional mean of $Y|X=E_{p,q}$, here, we thus provide a way to estimate conditional quantiles of order $\tau$. When $\tau=0.5$, it actually estimates the conditional median, which is known to be an indicator of central tendency that is more robust than the mean in the presence of outliers. On the other hand, for large and small $\tau$'s (for example the $0.05$ and $0.95$ quantiles), this allows to build confidence intervals for $Y|X=E_{p,q}$. Confidence bounds for the entries of matrices in matrix completion problems are something new up to our knowledge. 

The following result studies a particular case in which the Bernstein Assumption is proved in Proposition~\ref{prop:quantile_loss}. Following \cite{MR3526202}, it assumes that the conditional distribution of $Y$ given $X$ is continuous and that the density is not too small on the domain of interest -- this ensures that Bernstein's condition is satisfied with $\kappa=1$ and $A$ depending on the lower bound on the density, see Section~\ref{sec:bernstein_margin_condition} for more details. It can easily be derived for a specific distribution such as Gaussian, Student and even Cauchy. But we also have to assume that $\overline{M}\in b B_\infty$, or in other words $\overline{M}=M^*$, which is a more stringent assumption: in practice, it meands that we should know {\it a priori} an upper bound $b$ on the quantiles to be estimated.

\begin{theo}[Quantile matrix completion]\label{theo:matrix_quantile}
Assume that Assumption~\ref{ass:matrix_design} holds. Let $b>0$ and assume that $\overline{M} \in b B_\infty$. Assume that for any $(p,q)$, $Y|X=E_{p,q}$ has a density with respect to the Lebesgue measure, $g$, and that $g(u)>1/c$ for some constant $c>0$ for any $u$ such that $|u-\overline{M}_{i,j}|\leq 2b $. Let $s\in\{1, \ldots, \min(m,T)\}$ and assume that there exists a matrix with rank at most $s$ in $\overline{M}+(\rho^*_s/20)B_{S_1}$ where $\rho^*_s$ is defined in \eqref{eq:choice_rho_star}. Then, with probability at least 
\begin{equation*}
 1- {\bf C} \exp\left(-{\bf C} s \max(m,T) \log(m+T) \right)
 \end{equation*}
 the estimator
\begin{equation}
\widehat{M} \in\argmin_{M\in b B_\infty}\left(\frac{1}{N}\sum_{i=1}^N \rho_{\tau}(Y_i-\inr{X_i,M}) + \lambda \norm{M}_{S_1}\right)
\end{equation}
with $\lambda = c_0(\underbar{c},\bar{c}) \sqrt{\log(m+T)/(N \min(m,T))}$ satisfies
 \begin{align*}
  \frac{1}{mT} \norm{\widehat{M} - \overline{M}}_{S_1}& \leq  {\bf C} \min \left\{s
\sqrt{\frac{\log(m+T)}{N \min(m,T)}},\frac{\norm{\overline{M}}_{S_1}}{mT} \right\} ,\\ 
 \frac{1}{\sqrt{mT}}\norm{\widehat{M}-\overline{M}}_{S_2}& \leq {\bf C} \min\left\{ \sqrt{\frac{s(m+T) \log(m+T)}{N} }, \norm{\overline{M}}_{S_1}^{\frac{1}{2}} \left(\frac{\log(m+T)}{N \min(m,T)}\right)^{\frac{1}{4}} \right\} \\  
 \cE_{quantile}(\widehat{M}) & \leq {\bf C} \min \left\{  \frac{s(m+T)\log(m+T)}{N},\norm{\overline{M}}_{S_1} \sqrt{\frac{\log(m+T)}{N \min(m,T)}}\right\}.
 \end{align*}
\end{theo}

We obtain the same rate as for the penalized least squares estimator that is $\sqrt{s(m+T) \log(m+T)/N}$ (cf. \cite{rohde2011estimation,MR2906869}). 

\subsubsection*{Simulation study.} The goal of this part is to challenge the regularized least squares estimator by the RERM with quantile loss. The quantile used here is therefore the median. The main conclusion of our study is that median based estimators are more robust to outliers and noise than mean based estimators. We first test them on simulated datasets and then turn to use a real dataset.

\paragraph{Simulated matrices.}
The observations come from a base matrix $M^\star$ which is a $200\times 200$ low rank matrix. It is built by $M^\star=L R^\top$ where the entries of $L, R$ are i.i.d. gaussian and $L,R$ have $3$ columns (and therefore, the rank of $M^\star$ is $3$). The $X_i$'s correspond to $20\%$ randomly picked entries. The criterion that we retain is the $l_1$ reconstruction of $M^\star$ that is: $1/mT \sum_{p,q} |M^\star_{p,q}-M_{p,q}|$.

The observations are made according to this flexible model:
\begin{equation*}
Y_i = \inr{M^\star,X} + z_i + o \zeta_i.
\end{equation*}
$z_i$ is the noise, $o$ is the magnitude of outliers and $\zeta_i$ is the outlier indicator parametrized by the share $p$ such that $\zeta_i = p/2 \delta_{-1} + (1-p)\delta_0 + p/2\delta_1$.
The different parameters for the different scenarios are summarized in Table~\ref{tab:scheme_simulations_quantile}.

\begin{figure}[h!]
\centering
\includegraphics[scale=.5]{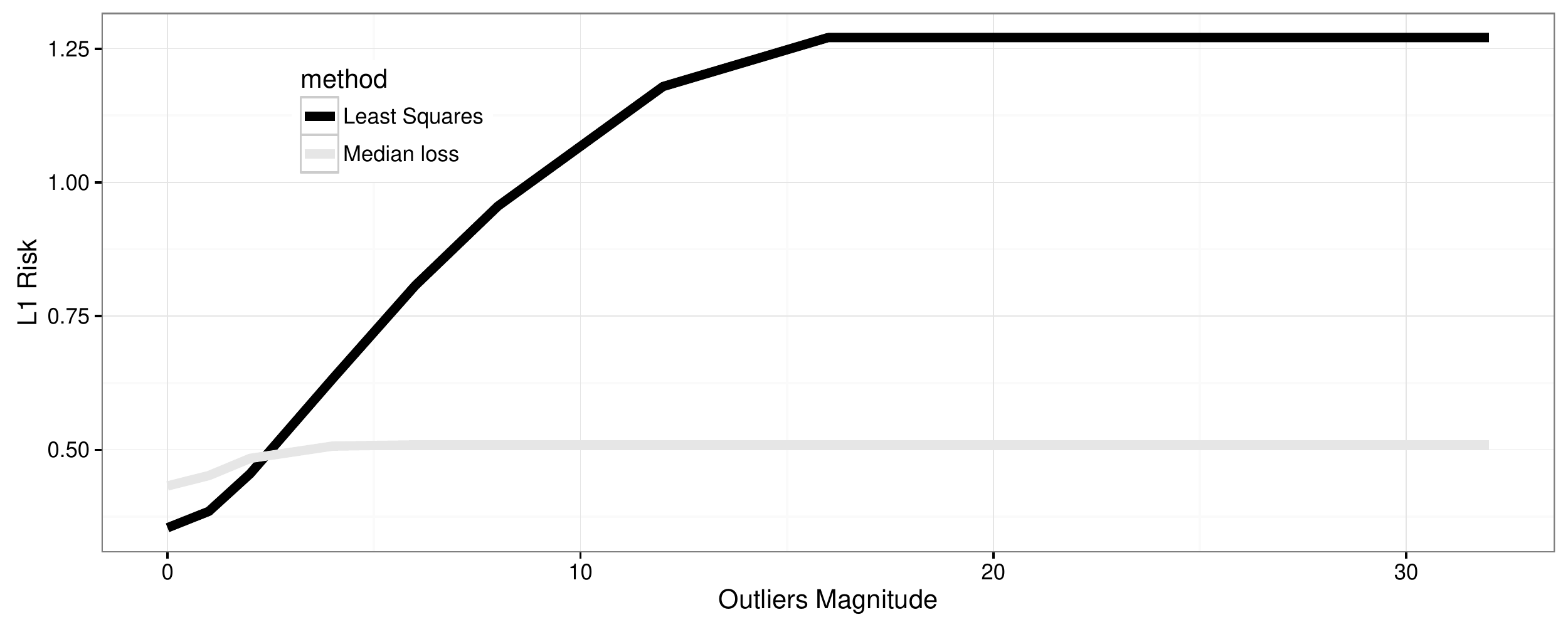}
\caption{\label{fig:MedianSimu}$l_1$ reconstruction for different magnitude of outliers}
\includegraphics[scale=.5]{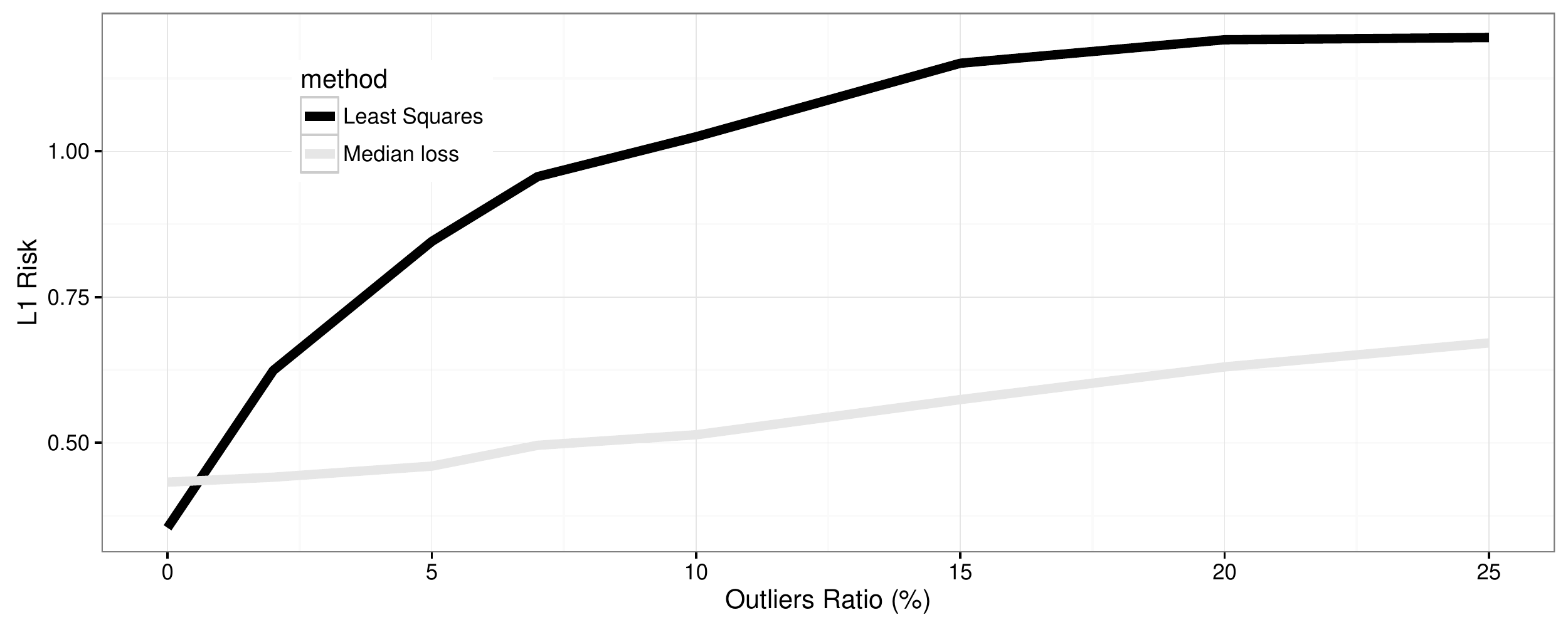}
\caption{\label{fig:MedianSimuPCT}$l_1$ reconstruction for different percentage of outliers}
\includegraphics[scale=.5]{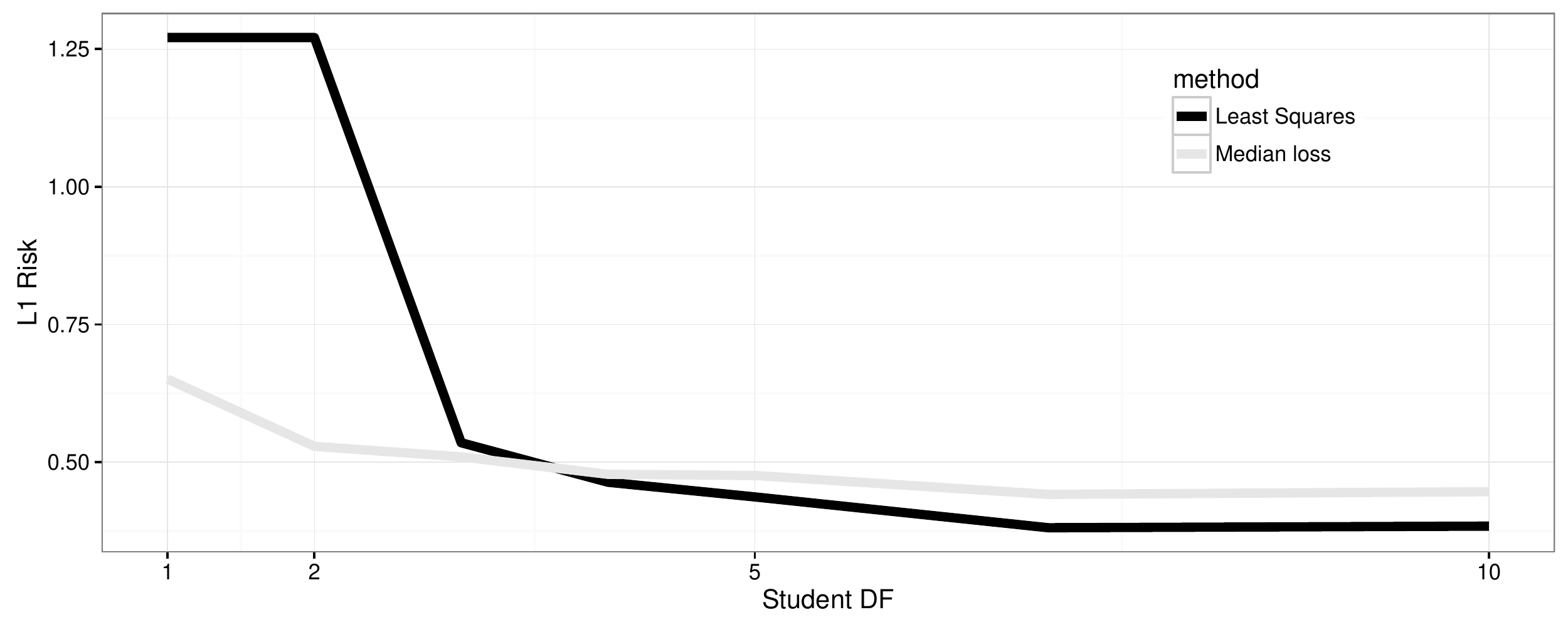}
\caption{\label{fig:MedianSimu_Student}$l_1$ reconstruction for student noise with various magnitude degrees of freedom}
\end{figure}

On the first experiment, $p$ is fixed to $10\%$ and the magnitude $o$ increases. As expected for least squares, the results are better for low magnitude of outliers (it corresponds to the penalized maximum likelihood estimator), see Figure \ref{fig:MedianSimu}. Quickly, the performance of the least squares estimator is getting worse and when the outliers are large enough, the best least squares predictor is a matrix with null entries. In opposite to this estimator, the median of the distribution is almost not affected by outliers and it is completely in line with the results: the performances are strictly the same for mid-range to high-range magnitude of outliers. The robustness of the quantile reconstruction is totally independent to the magnitude of the outliers. 

 \begin{table}[h!]
 \centering
 \begin{tabular}{lcccc}
  & $z_i$ & $o$ & $\zeta_i$ \\
 Figure~\ref{fig:MedianSimu} & $\cN(0,1/4)$ & $o=0..30$ & $p=0.1$ \\
 Figure~\ref{fig:MedianSimuPCT} &  $\cN(0,1/4)$ & $10$ & $p=0..0.25$ \\
 Figure~\ref{fig:MedianSimu_Student} & $t_\alpha, \alpha = 1..10$ & $0$ & $p=0$  
 \end{tabular}\\
 
 $t_\alpha$: t-distribution with $\alpha$ degrees of freedom.
 \caption{\label{tab:scheme_simulations_quantile}Parameters and distributions of the simulations}
 \end{table}

A second experiment involves fixed magnitude of outliers but the share of them increases, see Figure \ref{fig:MedianSimuPCT}. The median completion is, as expected, more robust and the results deteriorate less than the ones from least squares. When the outliers ratio is greater than $20\%$, the least squares estimator completely fails while the median completion still works. 

The third simulation involves non gaussian noise without outliers: we use the t-distribution, that has heavy tails. In this challenge, a lower degree of freedom involves heavier tails and the worst case is for Student distribution with degree $1$. We can see that the least squares is inadequate for small degrees of freedom ($1$ to $2$) and behaves better than the median completion for larger degrees of freedom, see Figure \ref{fig:MedianSimu_Student}.

\paragraph{Real dataset.} The last experiment involves the MovieLens dataset. We keep one fifth of the sample for test set to check the prediction accuracy. Even though the least squares estimator remains very efficient in the standard case, see Table \ref{tab:MedianMovieLens}, the results are quite similar for the MAE criterion. In a second step, we add artificial outliers. In order to do that, we change $20\%$ of $5$ ratings to $1$ ratings. It can be seen as malicious users that change ratings in order to distort the perception of some movies. As expected, it depreciates the least squares estimator performance but the median estimator returns almost as good performances as in the standard case.

\begin{table}
\centering
\begin{tabular}{lcc}
 & MSE & MAE \\
 Raw Data, LS 		& 0.89 & 0.75 \\
 Raw Data,  Median 	& 0.93	& 0.75 \\
 Outliers, LS 		& 1.04 	& 0.84 \\
 Outliers, Median 	& 0.96	& 0.78
\end{tabular}
\caption{\label{tab:MedianMovieLens}Prediction power of Least Squares and Median Loss on MovieLens 100K dataset}
\end{table}

\section{Kernel methods via the hinge loss and a RKHS-norm regularization} 
\label{sec:any_lipschitz_loss_rkhs_norm}
In this section, we consider regularization methods in some general
Reproducing Kernel Hilbert Space (RKHS) (cf. \cite{MR1864085}, Chapter~4 in \cite{MR2450103}  or Chapter~3 of \cite{MR1641250} for  general references on RKHS).

Unlike the previous examples, the regularization norm here, which is the norm $\norm{\cdot}_{\cH_K}$ of a RKHS $\cH_K$,  is not associated with some "hidden" concept of sparsity. In particular, RKHS norms have no singularity since they are differentiable at any point except in $0$. As a consequence the sparsity parameter $\Delta(\rho)$ cannot be larger than $4\rho/5$, i.e. $\rho$ does not satisfy the sparsity equation, unless the set $\Gamma_{f^*}(\rho)$ contains $0$ that is for $\rho\geq20 \norm{f^*}_{\cH_K}$. Indeed, one key observation is that any norm is non differentiable at $0$ and that its subdifferential at $0$ is somehow extremal:
\begin{equation}\label{eq:non-diff-norm}
 \partial\norm{\cdot}(0) = B_*:=\{f:\norm{f}_*\leq1\},
 \end{equation} where $\norm{\cdot}_*$ is the dual norm.
 
As a consequence, the rates obtained in this section do not depend on some \textit{hidden sparsity parameter} associated with the oracle $f^*$ but on the RKHS norm at $f^*$, that is $\norm{f^*}_{\cH_K}$. The aim of this section is therefore to show that our main results apply beyond ``sparsity inducing regularization methods'' by showing that ``classic'' regularization method, inducing smoothness for instance, may also be analyzed the same way and fall into the scope of Theorem~\ref{thm:main_subgaussian_simple} and Theorem~\ref{thm:main_bounded_simple}. This section also shows an explicit expression for the Gaussian mean-width with localization as used in Definition~\ref{def:function_r} (a sharper way to measure statistical complexity via a local $r(\cdot)$ function provided below).

\paragraph{Mathematical background} In this setup, the data are still $N$ i.i.d. pairs $(X_i,Y_i)_{i=1}^N$ where the $X_i$'s take their values in some set $\cX$ and $Y_i\in\{-1,+1\}$. A "similarity measure" is provided over the set $\cX$ by means of a kernel $K:\cX\times \cX\to \R$ so that $x_1,x_2\in\cX$ are "similar" when $K(x_1,x_2)$ is small. One can think for instance of $\cX$  the set of all DNA sequences (that is finite words over the alphabet $\{A, T, C, G \}$) and $K(w_1, w_2)$ is the minimal number of changes like insertion, deletion and mutation needed to transform word $w_1\in\cX$ into word $w_2\in\cX$.

The core idea behind kernel methods is to transport the design data $X_i$'s from $\cX$ to a Hilbert space via the application $x\to K(x,\cdot)$ and then construct statistical procedures based on the "transported" dataset $(K(X_i,\cdot), Y_i)_{i=1}^N$. The advantage of doing so is that the space where the $K(X_i,\cdot)$'s belong have much structure than the initial set $\cX$ which may have no algebraic structure at all. The first thing to set is to define somehow the "smallest" Hilbert space containing all the functions $x\to K(x,\cdot)$. We recall now one classic way of doing so that will be used later to define the objects that need to be considered in order to construct RERM in this setup and to obtain estimation rates for them via Theorem~\ref{thm:main_subgaussian_simple} and Theorem~\ref{thm:main_bounded_simple}.

Recall that if $K:\cX\times
\cX\rightarrow \R$ is a positive definite kernel such that $\norm{K}_{L_2}<\infty$, then by Mercer's
theorem, there is an orthogonal basis $\left(\phi_i\right)_{i\in\N}$ of
$L_2$ such that $\mu\otimes\mu$-almost surely,
$K(x,x')=\sum_{i=1}^\infty \lambda_i \phi_i(x)\phi_i(x')$ where
$(\lambda_i)_{i\in\N}$ is the sequence of eigenvalues of the positive self-adjoint integral
operator $T_K$ (arranged in a non-increasing order) defined for every
$f\in L_2$ and $\mu$-almost every $x\in\cX$ by
\begin{equation*}
 (T_K f)(x)=\int K(x,x')f(x')d\mu(x').
\end{equation*}In particular, for all $i\in\N$, $\phi_i$ is an eigenvector of $T_K$
corresponding to the eigenvalue $\lambda_i$; and $(\phi_i)_i$ is an orthonormal system in $L_2$.

The reproducing kernel Hilbert space $\cH_K$ is the set of all
function series $\sum_{i=1}^\infty a_i K(x_i,\cdot)$ converging in
$L_2$ endowed with the inner product
\begin{equation*}
 \inr{\sum a_i K(x_i,\cdot),\sum b_j K(x'_j,\cdot)}=\sum_{i,j}a_i b_j K(x_i,x'_j)
\end{equation*}where $a_i,b_j$'s are any real numbers and the $x_i$'s and $x'_j$'s are any points in $\cX$.

\paragraph{Estimator.} The RKHS $\cH_K$ is therefore a class of functions from $\cX$ to $\R$ that can be used as a learning model and the norm naturally associated to its Hilbert structure can be used as a regularization function. Given a Lipschitz loss function $\ell$,  the oracle is defined as
\begin{equation*}
 f^*\in\argmin_{f\in\cH_K}\E \ell_f(X, Y)
\end{equation*}and it is believed that $\norm{f^*}_{\cH_K}$ is small
which justified the use of the RERM with regularization function given
by the RKHS norm $\norm{\cdot}_{\cH_K}$:
 \begin{equation*}
     \hat f \in\argmin_{f\in
   \cH_K}\left(\frac{1}{N}\sum_{i=1}^N\ell_f(X_i,Y_i)+\lambda \norm{f}_{\cH_K}\right)
 \end{equation*}

Statistical properties of this RERM may be obtained from Theorem~\ref{thm:main_subgaussian_simple} in the subgaussian case and from Theorem~\ref{thm:main_bounded_simple} in the bounded case. To that end, we only have to compute the Gaussian mean width and/or the Rademacher complexities  of $B_{\cH_K}$. In this example, we rather compute the localized version of those quantities because it is possible to derive explicit formula. They are obtained by intersecting the ball with $r \cE$. In order not to induce any confusion, we still use the global ones in estimation bounds. 

\paragraph{Localized complexity parameter.} The goal is to compute $w(\rho B_{\cH_K}\cap r \cE)$ and ${\rm Rad}(\rho B_{\cH_K}\cap r \cE)$ for all $\rho,r>0$ where $B_{\cH_K}=\{f\in\cH_K:\norm{f}_{\cH_K}\leq1\}$ is the unit ball of the RKHS and $\cE=\{f\in\cH_K:\E f(X)^2\leq1\}$ is the ellipsoid associated with $X$. In the following, we embed the two sets $B_{\cH_K}$ and $\cE$ in $l_2=l_2(\N)$ so that we simply have to compute the Gaussian mean width and the Rademacher complexities of the intersection of two ellipsoids sharing the same coordinates structure.

The unit ball of $\cH_K$ can be constructed from the
eigenvalue decomposition of $T_K$ by considering the feature map
$\Phi:\cX\rightarrow l_2$ defined by
$\Phi(x)=\left(\sqrt{\lambda_i}\phi_i(x)\right)_{i\in\N}$ and then the
unit ball of $\cH_K$ is just
\begin{equation*}
 B_{\cH_K}=\big\{f_\beta(\cdot)=\inr{\beta,\Phi(\cdot)}:\norm{\beta}_{l_2}\leq1\big\}.
\end{equation*}One can use the feature map $\Phi$ to show that there is an isometry between the two Hilbert spaces
$\cH_K$ and $l_2$ endowed with the norm $\norm{\beta}_K=\left(\sum
\beta_i^2/\lambda_i\right)^{1/2}$. The unit ball of $l_2$ endowed with the norm $\norm{\cdot}_K$ is an ellipsoid denoted by $\cE_K$.

Let us now determine the ellipsoid in $l_2$  associated with the design $X$ obtained via this natural isomorphism $\beta\in l_2\to f_\beta(\cdot) = \inr{\beta, \Phi(\cdot)}\in\cH_K$ between  $l_2$ and $\cH_K$. Since $(\phi_i)_i$ is an orthonormal system
in $L_2$, the covariance operator of $\Phi(X)$ in $l_2$ is simply the
diagonal operator with diagonal elements $(\lambda_i)_i$. As a consequence the
ellipsoid associated with $X$ is isomorphic to
$\tilde\cE=\{\beta \in l_2:\E\inr{\beta,\Phi(X)}^2\leq1\}$; it has the same coordinate structure as the canonical one in
$l_2$ endowed with $\norm{\cdot}_K$: $\tilde\cE=\{\beta\in l_2:\sum \lambda_i\beta_i^2\leq1\}$.   So that, we obtain
\begin{equation}\label{eq:gauss_mean_width_rkhs}
w(K_\rho(f^*)\cap r \cE_{f^*})= w(\rho \cE_{K}\cap
r\tilde\cE)\sim \left(\sum_{j}(\rho^2\lambda_j)\wedge r^2\right)^{1/2}
\end{equation}
where the last inequality follows from
Proposition~2.2.1 in \cite{Talagrand:05} (note that we defined the
Gaussian mean widths in Definition~(\ref{def:gauss_mean_width})
depending on the covariance of $X$). We also get from Theorem~2.1 in \cite{MR2075996} that
\begin{equation}\label{eq:rademacher_rkhs}
{\rm Rad}(K_\rho(f^*)\cap r \cE_{f^*})\sim \left(\sum_{j}(\rho^2\lambda_j)\wedge r^2\right)^{1/2}.
\end{equation} 
Note that unlike the previous examples, we do not have to assume isotropicity of the design. Indeed, in the RKHS case, the unit ball of the regularization function is isomorphic to the  ellipsoid $\cE_K$. Since $\cE$ is also an ellipsoid having the same coordinates structure as $\cE_K$ (cf. paragraph above), for all $\rho, r>0$, the intersection $\rho B_{\cH_K}\cap r \cE$ is equivalent to an ellipsoid, meaning that, it contains an ellipsoid  and is contained in a multiple of this ellipsoid. Therefore, the Gaussian mean width and the Rademacher complexity of $\rho B_{\cH_K}\cap r \cE$ has  been computed without assuming isotropicity (thanks to general results on the complexity of Ellipsoids from Proposition~2.2.1 in \cite{Talagrand:05} and Theorem~2.1 in \cite{MR2075996}). 

It follows from \eqref{eq:gauss_mean_width_rkhs} and \eqref{eq:rademacher_rkhs} that the Gaussian mean width and the Rademacher complexities are equal. Therefore, up to constant ($L$ in the subgaussian case and $b$ in the bounded case), the two subgaussian and bounded setups may be analyzed at the same time. Nevertheless, since we will only consider in this setting  the hinge loss and that the Bernstein condition (cf.  Assumption~\ref{ass:margin}) with respect to the hinge loss has been studied in Proposition~\ref{prop:margin_hinge_loss} only in the bounded case. We therefore continue the analysis only for the bounded framework. 

We are now able to identify the complexity parameter of the problem. We actually do not use the localization in this and rather use only the global complexity parameter as defined in Definition~\ref{def:function_r_boundedness_simple}: for all $\rho>0$:
\begin{equation}
 \label{eq:RKHS_function_r}
r(\rho)=\left[\frac{ {\bf C} \rho\left(\sum_j\lambda_j\right)^{1/2}}{\sqrt{N}}\right]^{\frac{1}{2\kappa}}
\end{equation}
where $\kappa\geq1$ is the Bernstein parameter.

\paragraph{Results in the bounded setting} Finally, let us discuss about the boundedness assumption. It is known (cf., for instance, Lemma~4.23 in \cite{MR2450103}) that if the kernel $K$ is bounded then the functions in the RKHS $\cH_K$ are bounded: for any $f\in\cH_K$, $\norm{f}_{L_\infty}\leq \norm{K}_{\infty} \norm{f}_{\cH_K}$ where $\norm{K}_{\infty} := \sup_{x\in\cX}\sqrt{K(x, x)}$. As a consequence, if one restricts the search space of the RERM to a RKHS ball of radius $R$, one has $F := R B_{\cH_K}\subset \norm{K}_{\infty} B_{L_\infty}$ and therefore the boundedness assumption is satisfied by $F$. However, note that a refinement of the proof of Theorem~\ref{thm:main_bounded} using a boundedness parameter $b$ depending on the radius of the RKHS balls used while performing  the peeling device yields statistical properties for the RERM with no search space constraint. For the sake of shortness, we do not provide this analysis here.
 
We are now in a position to provide estimation and prediction results for the RERM
\begin{equation}\label{eq:rerm_rkhs}
\hat f \in\argmin_{f\in
   R B_{\cH_K}}\left(\frac{1}{N}\sum_{i=1}^N (1-Y_if(X_i))_+ +\frac{{\bf C}\left(\sum_j\lambda_j\right)^{1/2}}{ \sqrt{N}} \norm{f}_{\cH_K}\right)
\end{equation}
where the choice of the regularization parameter $\lambda$ follows from Theorem~\eqref{thm:main_bounded_simple} and \eqref{eq:rademacher_rkhs} (for $r=+\infty$). Note that unlike the examples in the previous sections, we do not have to find some radius $\rho^*$ satisfying the sparsity equation~\eqref{eq:sparsity_equation} to apply Theorem~\ref{thm:main_bounded_simple} since we simply take $\rho^* = 20 \norm{f^*}_{\cH_K}$ to insure that $0\in\Gamma_{f^*}(\rho^*)$.

\begin{theo}\label{theo:rkhs}
 Let $\cX$ be some space, $K:\cX\times\cX\to\R$ be a bounded kernel and denote by $\cH_K$ the associated RKHS. Denote by $(\lambda_i)_i$ the sequence of eigenvalues associated to $\cH_K$ in $L_2$. Assume that the Bayes rule $\overline{f}$ from \eqref{eq:bayes_rule} belongs to $R B_{\cH_K}$ and that the margin assumption \eqref{eq:margin_tsybakov} is satisfied for some $\kappa\geq1$. 

  Then the RERM defined in \eqref{eq:rerm_rkhs} satisfies with probability larger than 
\begin{equation*}
1-{\bf C} \exp\left(-{\bf C} N^{1/2\kappa}\left(\norm{\overline{f}}_{\cH_K}\left(\sum_j\lambda_j\right)^{1/2}\right)^{(2\kappa-1)/\kappa}\right),
\end{equation*}
  that
 \begin{equation*}
   \norm{\hM-\overline{f}}_{L_2}\leq {\bf C} \left[\frac{\norm{\overline{f}}_{\cH_K} \left(\sum_j\lambda_j\right)^{1/2}}{\sqrt{N}}\right]^{1/2\kappa} \mbox{ and } \cE_{hinge}(\hM)\leq {\bf C} \frac{\norm{\overline{f}}_{\cH_K} \left(\sum_j\lambda_j\right)^{1/2}}{\sqrt{N}}
 \end{equation*}where $\cE(\hM)$ is the excess hinge risk of $\hM$.
\end{theo}

Note that classic procedures in the literature on RKHS are mostly developed in the classification framework. They are usually based on the hinge loss
and the regularization function is the square of the RKHS norm. For such procedures,
oracle inequalities have been obtained in Chapter~7 from
\cite{MR2450103} under the margin assumption (cf. \cite{Tsy04}). A result that is close to the one obtained in Theorem~\ref{theo:rkhs} is Corollary~4.12 in \cite{MR2590050}. Assuming that $\norm{Y}_\infty\leq {\bf C}$,  $\cX\subset\R^d$,  $\norm{K}_\infty\leq1$, that the eigenvalues of the integral operator satisfies
\begin{equation}\label{eq:weak_lp_eigenvalues}
\lambda_i\leq c i^{-1/p}
\end{equation}for some $0<p<1$ and that the eigenvectors $(\phi_i)$ are such that $\norm{\phi_i}_\infty\leq A$ for any $i$ and some constant $A$ then the RERM $\tilde f$ over the entire RKHS space, w.r.t. the quadratic loss and for a regularization function of the order of  (up to logarithmic terms)
\begin{equation}\label{eq:reg_func_shahar}
f\mapsto \rho(\norm{f}_\cH):= \max\left(\frac{\norm{f}_{\cH}^{2p/(1+p)}}{N^{1/(1+p)}}, \frac{\norm{f}_{\cH}^2}{N}\right)
\end{equation}satisfies with large probability an oracle inequality like
\begin{equation*}
\E(Y-\tilde f(X))^2\leq \inf_{r\geq 1}\left(\inf_{\norm{f}_{\cH}\leq r}\E(Y-f(X))^2+{\bf C} \rho(r)\right).
\end{equation*} In particular, an error bound (up to log factors) follows from this result: with high probability,
\begin{equation}\label{eq:shahar_kernel_meth}
\norm{\tilde f - f^*}_{L_2}^2\leq {\bf C} \rho(\norm{f^*}_{\cH}) = {\bf C}\max\left(\frac{\norm{f^*}_{\cH}^{2p/(1+p)}}{N^{1/(1+p)}}, \frac{\norm{f^*}_{\cH}^2}{N}\right).
\end{equation}One may compare this result to the one from Theorem~\ref{theo:rkhs} under assumption~\eqref{eq:weak_lp_eigenvalues} even though the two procedures $\tilde f$ and $\hat f$ use different loss functions,  regularization function and different search space. If assumption~\eqref{eq:weak_lp_eigenvalues} holds then $\left(\sum_j \lambda_j\right)^{1/2} \leq c $ and so, one can take $r(\rho)= \left( {\bf C} c \rho /(\theta\sqrt{N}) \right)^{1/(2\kappa)}$ and $\lambda ={\bf C} \sqrt{C/N}$. For such a choice of regularization parameter, Theorem~\ref{theo:rkhs} provides an error bound of the order of
\begin{equation}
\norm{\tilde f - \overline{f}}_{L_2(\mu)}^2\leq {\bf C} \left[\frac{\norm{\overline{f}}_{\cH_K} C}{\sqrt{N}}\right]^{1/\kappa}
\end{equation}which is almost the same as the one obtained in \eqref{eq:shahar_kernel_meth} when $\kappa=1$ and $p$ is close to $1$. But our result is worse when $\kappa>1$ and $p$ is far from $1$. This is the price that we pay by using the hinge loss -- note that the quadratic loss satisfies the Bernstein condition with $\kappa=1$ --  and by fixing a regularization function which is the norm $\norm{\cdot}_{\cH_K}$ instead of fitting the regularization function in a ``complexity dependent way'' as in \eqref{eq:reg_func_shahar}. In the last case, our procedure $\hat f$ does not benefit from the ``real complexity'' of the problem which is localized Rademacher complexities -- note that we used global Rademacher complexities to fit $\lambda$ and construct the complexity function $r(\cdot)$.

\section{A review of the Bernstein  and margin conditions} 
\label{sec:bernstein_margin_condition}

In order to apply the main results from Theorem~\ref{thm:main_subgaussian_simple} and Theorem~\ref{thm:main_bounded_simple}, one has to check the Bernstein condition. This section is devoted to the study of this condition for three loss functions: the hinge loss, the quantile loss and the logistic loss. This condition has been extensively studied in Learning theory (cf. \cite{bartlett2003large,zhang2004statistical,MR2426759,MR2240689,van2008high,elsener2016robust}). We can identify mainly two approaches to study this condition: when the class $F$ is convex and the loss function $\ell$ is ``strongly convex'', then the risk function inherits this property and automatically satisfies the Bernstein condition (cf. \cite{bartlett2003large}). On the other hand, for loss functions like the hinge or quantile loss, that are affine by parts, one has to use a different path. In such cases, one may go back to a statistical framework and try to check the margin assumption. As a consequence, in the latter case, the Bernstein condition is usually more restrictive and requires strong assumptions on the distribution of the observations.

\subsection{Logistic loss} 
\label{sub:logistic_loss}
In this section, we study the Bernstein condition of the logistic loss function which is defined for every $f:\cX\to\R$, $x\in\cX$, $y\in\{-1, 1\}$ and $u\in\R$ by
\begin{equation*}
 \ell_f(x, y)  = \tilde{\ell}(yf(x)) \mbox{ where } \tilde{\ell}(u) = \log(1+\exp(-u)).
 \end{equation*} 
Function  $\tilde{\ell}$ is strongly convex on every compact interval in $\R$. As it was first observed in \cite{bartlett2003large,bartlett2006convexity}, one may use this property to check the Bernstein condition for the loss function $\ell$.  This approach was  extended to the bounded regression problem with respect to $L_p$ loss  functions ($1<p<\infty$) in~\cite{mendelson2002improving} and to non convex classes in~\cite{MR2426759}.

In the bounded scenario, \cite{bartlett2006convexity} proved that the logistic loss function satisfies the Bernstein condition for $\kappa=1$.  One may therefore use that result to apply Theorem~\ref{thm:main_bounded_simple}. The analysis is pretty straightforward in the bounded case. It becomes more delicate in the subgaussian scenario as considered in Theorem~\ref{thm:main_subgaussian_simple}.

\begin{proposition}[\cite{bartlett2003large}]
\label{prop:margin_logistic_loss}
Let $F$ be a convex class of functions from $\cX$ to $\R$. Assume that for every $f\in F$, $\norm{f}_{L_\infty}\leq b$. Then the class $F$ satisfies the Bernstein condition with Bernstein parameter $\kappa=1$ and constant $A= 4\exp(2b)$.
\end{proposition}

This result solves the problem of the Bernstein condition with respect to the logistic loss function over a convex class $F$ of functions as long as all functions in $F$ are uniformly bounded by some constant $b$. We will therefore use this result only in the bounded framework, for instance, when $F$ is a class of linear functional indexed by a bounded set of vectors and when the design takes its values in the canonical basis. 

In the subgaussian framework, one may proceed as in \cite{van2008high} and assume that a statistical model holds. In that case, the Bernstein condition is reduced to the study of the Margin assumption since, in that case, the ``Bayes rule'' $\overline{f}$ (which is called the log-odds ratio in the case of the logistic loss function) is assumed to belong to the class $F$ and so $f^*=\overline{f}$. The margin assumption with respect to the logistic loss function has been studied in Example~1 from \cite{van2008high} but for a slightly different definition of the Margin assumption. Indeed, in \cite{van2008high} only functions $f$ in a $L_\infty$ neighborhood of $\overline{f}$ needs to satisfy the Margin assumption whereas in Assumption~\ref{ass:margin} it has to be satisfied in the non-bounded set $\cC$. 

From our perspective, we do not want to make no ``statistical modeling assumption''. In particular, we do not want to assume that $\overline{f}$ belongs to $F$. We therefore have to prove the Bernstein condition when $\bar f$ may not belong to $F$. We used this result in Section~\ref{sec:the_logistic_lasso} in order to obtain statistical bounds for the Logistic LASSO and Logistic Slope procedures. In those cases, $F$ is a class of linear functionals. We now state that the Bernstein condition is satisfied for a class of linear functional when $X$ is a standard Gaussian vector. 


\begin{proposition}\label{prop:bernstein_logistic_subgaussian}
Let $F=\{\inr{\cdot, t}:t\in R B_{l_2}\}$ be a class of linear functionals indexed by $RB_{l_2}$ for some radius $R\geq1$. Let $X$ be a standard Gaussian vector in $\R^d$ and let $Y$ be a $\{-1, 1\}$ random variable. For every $f\in F$,   the excess logistic risk of $f$, denoted by $P\cL_f$, satisfies
\begin{equation*}
 \cE_{logistic}(f) = P \cL_f\geq \frac{c_0}{R^3} \norm{f-f^*}_{L_2}^2
 \end{equation*} where $c_0$ is some absolute constant.
\end{proposition}

\subsection{Hinge loss} 
\label{sub:hinge_loss}
Unlike the logistic loss function, both the hinge loss and the quantile losses does not enjoy a strong convexity property. Therefore, one has to turn to a different approach as the one used in the previous section to check the Bernstein condition for those two loss functions. 

For the hinge loss function, Bernstein condition is more stringent and is connected to the margin condition in classification. So, let us first introduce some notations specific to classification.  In this setup, one is given $N$ labeled pairs $(X_i, Y_i), i=1, \ldots, N$ where $X_i$ takes its values in $\cX$ and $Y_i$ is a label taking values in $\{-1,+1\}$.  The aim is to predict the label $Y$ associated with $X$ from the data when $(X,Y)$ is distributed like the $(X_i,Y_i)$'s. The classic loss function considered in this setup is the $0-1$ loss function $\ell_f(x,y)= I(y\neq f(x))$ defined for any $f:\cX\to\{-1,+1\}$. The $0-1$ loss function is not convex,  this may result in some computational issues when dealing with it. A classic approach is to use a `` convex relaxation function'' as a surrogate to the $0-1$ loss function: note that this is a way to motivate the introduction of the hinge loss $\ell_f(x, y)  = \max(1-yf(x), 0)$. It is well known that the Bayes rules minimizes both  the standard $0-1$ risk  as well as the hinge risk: put $\eta(x) := \bE[Y|X=x]$ for all $x\in\cX$ and define the Bayes rule as
\begin{equation}\label{eq:bayes_rule}
\overline{f}(x) = {\rm sgn}(\eta(x)),
\end{equation}
then $\overline{f}$ minimizes $f\to P \ell_f$ over all measurable functions from $\cX$ to $\R$ when $\ell_f$ is the hinge loss of $f$.

Let $F$ be a class of functions from $\cX$ to $[-1, 1]$. Assume that $\overline{f}\in F$ so that $\overline{f}$ is an oracle in $F$ and thus (using the notations from Section~\ref{sec:theoretical_results}) $f^*=\overline{f}$. In this situation, Margin assumption with respect to the hinge loss (cf. \cite{Tsy04,MR2364224}) restricted to the class $F$ and Bernstein condition (cf. Assumption~\ref{ass:margin}) coincide. Therefore, Assumption~\ref{ass:margin} holds when the Margin assumption w.r.t. the hinge loss holds. According to Proposition 1 in~\cite{MR2364224}, the Margin assumption with respect to the hinge loss is equivalent the Margin assumption with respect to the $0-1$ loss for a class $F$ of functions with values in $[-1,1]$. Then, according to Proposition 1 in~\cite{Tsy04} and \cite{MR2182250} the margin assumption with respect to the $0-1$ loss with parameter $\kappa$ is equivalent to
\begin{equation}\label{eq:margin_tsybakov}
\left\{
  \begin{array}{cc}
  \bP(|\eta(X)|\leq t) \leq c t^{\frac{1}{\kappa-1}}, \forall 0\leq t\leq1 & \mbox{ when } \kappa>1\\
   |\eta(X)|\geq \tau \mbox{ a.s. for some } \tau>0 & \mbox{ when } \kappa=1.
  \end{array}\right.
\end{equation} 
As a consequence, one can state the following result on the Bernstein condition for the hinge loss in the bounded case scenario.

\begin{proposition}[Proposition 1, \cite{MR2364224}]\label{prop:margin_hinge_loss}
Let $F$ be a class of functions from $\cX$ to $[-1,1]$. Define  $\eta(x)= \bE[Y|X=x]$ for all $x\in\cX$ and assume that the Bayes rule \eqref{eq:bayes_rule} belongs to $F$. If \eqref{eq:margin_tsybakov} is satisfied for some $\kappa\geq1$ then Assumption~\ref{ass:margin} holds with parameter $\kappa$ for the hinge loss, and $A$ depending on $c$, $\kappa$ and $\tau$ (which is explicitly given in the mentioned references). In the special case when $\kappa=1$ then $A=1/(2\tau)$.
\end{proposition}

Note that up to a modification of the constant $A$, the same result holds for functions with values in $[-b,b]$ for $b>0$, a fact we used in Section~\ref{sec:any_lipschitz_loss_rkhs_norm}.

\subsection{Quantile loss} 
\label{sub:quantile_loss}
 In this section, we study the Bernstein parameter of the \textbf{quantile loss} in the bounded regression model, that is when for all $f\in F,\|f\|_{L_\infty} \leq b$ a.s.. Let $\tau\in(0,1)$ and, for all $x\in\cX$, define $\overline{f}(x)$ as the quantile of order $\tau$ of $Y|X=x$ and assume that $\bar f$ belongs to $F$, in that case, $\overline{f} = f^*$ and Bernstein condition and margin assumption are the same. Therefore one may follow the study of the margin assumption for the quantile loss in \cite{elsener2016robust} to obtain the following result.

 \begin{proposition}[\cite{elsener2016robust}]\label{prop:quantile_loss}
 Assume that for any $x\in\cX$, it is possible to define
 a density $f_x$ w.r.t the Lebesgue measure for $Y|X=x$ such that  $f_x(u) \geq 1/C $ for some $C>0$ for all $u\in\R$ with $|u-f^*(x)| \leq 2 b$.
 Then the quantile loss satisfies the Bernstein's assumption with $\kappa=1$ and $A=2C$ over $F$.
 \end{proposition}

 \section{Discussion} 
 \label{sec:discussion}
This paper covers many aspects of the regularized empirical risk estimator (RERM) with Lipschitz loss. This property is commonly shared by many loss functions used in practice such as the hinge loss, the logistic loss or the quantile regression loss. This work offers a general method to derive estimation bounds as well as excess risk upper bounds. Two main settings are covered: the subgaussian framework and the bounded framework. The first one is illustrated by the classification problem with logistic loss. In particular, minimax rates are achieved when using  the SLOPE regularization norm. The second framework is used to derive new results on matrix completion and in kernel methods. 

A possible extension of this work is to study other regularization norms. In order to do that, one has to compute the complexity parameter in one of the settings and a solution of the sparsity equation. The latter usually involves to understand the sub-differential of the regularization norm and in particular its singularity points which are related to the sparsity equation.     
 

\section{Proof of Theorem~\ref{thm:main_subgaussian_simple} and Theorem~\ref{thm:main_bounded_simple}} 
\label{sec:proof_of_theorem_ref_thm_main}

\subsection{More general statements: Theorems~\ref{thm:main_bounded} and~\ref{thm:main_subgaussian}}
First, we state two theorems: Theorem~\ref{thm:main_subgaussian} in the subgaussian setting, and Theorem~\ref{thm:main_bounded} in the bounded setting. These two theorems rely on localized versions of the complexity function $r(\cdot)$ that will be defined first. Note that the localized version of $r(\cdot)$ can always be upper bounded by the simpler version used in the core of the paper. Thus, Theorem~\ref{thm:main_subgaussian_simple} is a direct corollary of Theorem~\ref{thm:main_subgaussian}, and Theorem~\ref{thm:main_bounded_simple} is a direct corollary of Theorem~\ref{thm:main_bounded}.

So let us start with a localized complexity parameters. The "statistical size" of the family of "sub-models" $(\rho B)_{\rho>0}$ is now measured by local Gaussian mean-widths in the subgaussian framework.
\begin{definition}\label{def:function_r}
Let $\theta>0$. The \textbf{complexity parameter} is a non-decreasing function $r(\cdot)$ such that for every $\rho\geq0$, 
\begin{equation*}
C L w\left(\rho B \cap r(\rho) B_{L_2}\right)\leq \theta r(\rho)^{2\kappa} \sqrt{N}
\end{equation*}
\end{definition}

In the boundedness case, it is written as follows.
\begin{definition}\label{def:function_r_boundedness}
Let $\theta>0$.
The \textbf{complexity parameter} is a non-decreasing function $r(\cdot)$ such that for every $\rho\geq0$, 
\begin{equation*}
48 {\rm Rad}(\rho B \cap r(\rho) B_{L_2}) \leq \theta r(\rho)^{2\kappa} \sqrt{N}
\end{equation*}where $\kappa$ is the Bernstein parameter from Assumption~\ref{ass:margin}.
\end{definition}

To obtain the  complexity functions from Definition~\ref{def:function_r_simple} and \ref{def:function_r_boundedness_simple}, we use the fact that $w\left(\rho B \cap r(\rho) B_{L_2}\right) \leq w(\rho B)$ and ${\rm Rad}(\rho B \cap r(\rho) B_{L_2}) \leq {\rm Rad}(\rho B )$: it indeed does not use the localization. We also set $\theta = 7/40A$ in those definitions because it is the largest value allowed in the following theorems.

\begin{theo}\label{thm:main_subgaussian}
Assume that Assumption~\ref{ass:lipschitz}, Assumption~\ref{ass:margin} and Assumption~\ref{ass:sub-gauss} hold where $r(\cdot)$ is a function as in Definition~\ref{def:function_r} for some $\theta$ such that $40A\theta\leq 7$ and assume that $\rho\to r(2\rho)/\rho$ is non-increasing. Let the regularization parameter $\lambda$ be chosen such that
\begin{equation}\label{eq:reg_param_choice}
 \frac{10\theta r(2\rho)^{2\kappa}}{7 \rho}< \lambda < \frac{r(2\rho)^{2\kappa}}{2A \rho}, \quad \forall \rho\geq \rho^*
 \end{equation} 
 where $\rho^*$ satisfies \eqref{eq:sparsity_equation}.  Then, with probability larger than
 \begin{equation}\label{eq:proba}
 1- \sum_{j=0}^\infty \sum_{i\in I_j}\exp\left(- \frac{\theta^2 N (2^{(i-1)\vee0}r(2^{j}\rho^*))^{4\kappa-2}}{4C^2L^2}\right)
 \end{equation}
 where for all $j\in\N$, $I_j=\{1\}\cup\{i\in\N^*: 2^{i-1} r(2^{j} \rho^*)\leq 2^{j} \rho^* d_{L_2}(B)\}$, we have
 \begin{align*}
 \norm{\hat f - f^*} \leq \rho^*,\quad
 \norm{\hat f - f^*}_{L_2} &\leq r(2\rho^*) \mbox{ and }
 \cE(\hat f) \leq  r(2\rho^*)^{2\kappa}/A. 
 \end{align*}
\end{theo}

 {\bf Proof of Theorem~\ref{thm:main_subgaussian_simple}:} Let $r(\cdot)$ be chosen as in \eqref{def:function_r_simple}. For this choice, one can check that the regularization parameter used for the construction of the RERM satisfies \eqref{eq:reg_param_choice} with an adequate constant choice. Moreover, for this choice of function $r(\cdot)$ it is straightforward to lower bound  the sum in the probability estimate in \eqref{eq:proba}. The parameter $\lambda$ is chosen in the middle of the range.     
 \endproof

The bounded case is in the same spirit.

\begin{theo}\label{thm:main_bounded}
Assume that Assumption~\ref{ass:lipschitz}, Assumption~\ref{ass:margin} and Assumption~\ref{ass:bounded} hold where $r(\cdot)$ is a function as in Definition~\ref{def:function_r_boundedness} for some $\theta$ such that $40A\theta\leq7$ and assume that $\rho\to r(2\rho)/\rho$ is non-increasing.
Let the regularization parameter $\lambda$ be chosen such that
\begin{equation}\label{eq:reg_param_choice_bounded}
 \frac{10\theta r(2\rho)^{2\kappa}}{7 \rho}< \lambda < \frac{r(2\rho)^{2\kappa}}{2A \rho}, \quad \forall \rho\geq \rho^*
 \end{equation}
 where $\rho^*$ satisfies \eqref{eq:sparsity_equation}.  Then, with probability larger than
 \begin{equation}\label{eq:proba_bounded}
 1- 2 \sum_{j=0}^\infty \sum_{i\in I_j}\exp\left(- c_0\theta^2 N (2^i r(2^{j+1}\rho^*))^{4\kappa-2}\right)
 \end{equation}
 where $ c_0 = 1/\max\left(48,207 \theta b^{2\kappa-1} \right) $ and for all $j\in\N$, $I_j:=\{1\}\cup\{i\in\N^*: 2^{i-1} r(2^j \rho^*) \leq \min(2^{j} \rho^* d_{L_2}(B), b)\}$, we have
 \begin{align*}
 \norm{\hat f - f^*} \leq \rho^*,\quad
 \norm{\hat f - f^*}_{L_2} &\leq r(2\rho^*) \mbox{ and }
 \cE(\hat f) \leq  r(2\rho^*)^{2\kappa}/A.
 \end{align*}
\end{theo}

 The proof of Theorem~\ref{thm:main_bounded_simple} is identical to the one of Theorem~\ref{thm:main_subgaussian_simple} and we do not reproduce it here. 


\subsection{Proofs of Theorems~\ref{thm:main_bounded} and~\ref{thm:main_subgaussian}}

Proof of Theorem~\ref{thm:main_subgaussian} and and Theorem~\ref{thm:main_bounded} follow the same strategy. They are split into two parts. First, we identify an event onto which the statistical behavior of the regularized estimator $\hM$ can be controlled using only deterministic arguments. Then, we prove that this event holds with a probability at least as large as the one in \eqref{eq:proba} in the case of Theorem~\ref{thm:main_subgaussian}  and as in \eqref{eq:proba_bounded} in the case of Theorem~\ref{thm:main_bounded}. We first introduce this event which is common to the subgaussian and the bounded setups:
\begin{equation*}
\Omega_0:=\left\{ \big|(P-P_N)\cL_f\big|\leq \theta \max\left(r(2 \max(\norm{f-f^*},\rho^*))^{2\kappa}, \norm{f-f^*}_{L_2}^{2\kappa}\right): \mbox{ for all } f\in F\right\}
\end{equation*} where $\theta$ is a parameter appearing in the definition of $r(\cdot)$ in Definition~\ref{def:function_r} and Definition~\ref{def:function_r_boundedness}, $\kappa\geq1$ is the Bernstein parameter from Definition~\ref{ass:margin} and $\rho^*$ is a radius satisfying the sparsity Equation~\eqref{eq:sparsity_equation}.

\begin{proposition}\label{prop:algebra}Let $\lambda$ be as in \eqref{eq:reg_param_choice} (or equivalently as in \eqref{eq:reg_param_choice_bounded})  and let $\rho^*$ satisfy \eqref{eq:sparsity_equation}, on the event $\Omega_0$, one has 
 \begin{align*}
 \norm{\hat f - f^*} \leq \rho^*,\quad
 \norm{\hat f - f^*}_{L_2} &\leq r(2\rho^*) \mbox{ and }
 \cE(\hat f) \leq \theta r(2\rho^*)^{2\kappa}.
 \end{align*}
 \end{proposition} 

\beginproof
Denote $\hat \rho = \norm{\hat f-f^*}$. We first prove that $\hat \rho < \rho^*$. To that end, we assume that the reverse inequality holds and show some contradiction. Assume that $\hat \rho\geq \rho^*$. Since $\rho\to r(2\rho)/\rho$ is non-increasing then by Lemma~\ref{lem:Delta}, $\rho\to \Delta(\rho)/\rho$ is non-decreasing and so we have
\begin{equation*}
\frac{\Delta(\hrho)}{\hrho}\geq \frac{\Delta(\rho^*)}{\rho^*}\geq \frac{4}{5}.
\end{equation*}
Now, we consider two cases: either $\norm{\hat f-f^*}_{L_2}\leq r(2\hrho)$ or $\norm{\hat f-f^*}_{L_2}> r(2\hrho)$. 

First assume that $\norm{\hat f-f^*}_{L_2}\leq r(2\hat \rho)$. Since $\Delta(\hrho)\geq 4\hrho/5$ and $h=\hM-f^*\in\hrho S \cap r(2\hrho)B_{L_2}$, it follows from the definition of the sparsity parameter $\Delta(\hrho)$ that there exists some $f\in F$ such that $\norm{f-f^*}\leq \hrho/20$ and for which
\begin{equation*}
\norm{f+h}-\norm{f}\geq \frac{4\hrho}{5}.
\end{equation*}It follows that
\begin{equation*}
\norm{\hM}-\norm{f^*} = \norm{f^*+h} - \norm{f^*}\geq \norm{f+h} - \norm{f} -2\norm{f-f^*}\geq \frac{4\hrho}{5}-\frac{\hrho}{10} = \frac{7\hrho}{10}.
\end{equation*}

Let us now introduce the excess regularized loss: for all $f\in F$, 
 \begin{equation*}
 \cL_f^\lambda = \cL_f + \lambda(\norm{f}-\norm{f^*}) = \left(\ell_f+\lambda \norm{f}\right) - \left(\ell_{f^*}+\lambda \norm{f^*}\right).
 \end{equation*}
On the event $\Omega_0$, we have
\begin{align*}
P_N \cL_{\hM}^\lambda &= P_N \cL_{\hM} + \lambda\left(\norm{\hM}-\norm{f^*}\right) \geq  (P_N-P)\cL_{\hM} + \lambda\left(\norm{\hM}-\norm{f^*}\right)\\
&\geq -\theta \max\left(r(2\hrho)^{2\kappa}, \norm{\hM-f^*}_{L_2}^{2\kappa}\right)  + \frac{7\lambda\hrho}{10}= -\theta r(2\hrho)^{2\kappa} + \frac{7\lambda\hrho}{10}>0
\end{align*}because by definition of $\lambda$,  $7\lambda\hrho > 10 \theta r(2\hrho)^{2\kappa}$. Therefore, $P_N \cL_{\hM}^\lambda>0$. But, by construction, one has $P_N \cL_{\hM}^\lambda\leq0$.

Then, assume that  $\norm{\hM-f^*}_{L_2}> r(2\hat \rho)$. In particular, $f\in\cC$ where $\cC$ is the set introduced in \ref{eq:set_cC}  below Assumption~\ref{ass:margin}. By definition of $\hM$ we have $P_N \cL^\lambda_{\hM}\leq0$ so it follows from Assumption~\ref{ass:margin} that
\begin{align} \nonumber
&\norm{\hM-f^*}_{L_2}^{2\kappa}\leq A P\cL_{\hM} = A \left[(P-P_N)\cL_{\hM} + P_N \cL_{\hM}^\lambda + \lambda\left(\norm{f^*} - \norm{\hM}\right)\right]\\
&\leq A\theta \max\left(r(2\hrho)^{2\kappa}, \norm{\hM-f^*}_{L_2}^{2\kappa}\right) +  A \lambda \norm{\hM-f^*} = A\theta \norm{\hM-f^*}_{L_2}^{2\kappa} + A \lambda \hrho. \label{proof_eq_1}
\end{align}
Hence, if $A\theta\leq 1/2$ then 
\begin{equation*}
r(2\hrho)^{2\kappa}\leq \norm{\hM-f^*}_{L_2}^{2\kappa}\leq 2 A \lambda \hrho.
\end{equation*}But, by definition of $\lambda$ one has $r(2\hrho)^{2\kappa} > 2 A \lambda \hrho$.

Therefore, none of the two cases is possible when one assumes that $\hrho\geq \rho^*$ and so we necessarily have $\hat \rho<\rho^*$.

Now, assuming that $\norm{\hM - f^*}_{L_2}> r(2\rho^*)$ and following \eqref{proof_eq_1} step by step also leads to a contradiction, so $\norm{\hM - f^*}_{L_2}\leq r(2\rho^*)$.

Next, we prove the result for the excess risk. One has
\begin{align*}
&P_N \cL_{\hM}^\lambda = P_N \cL_{\hM} + \lambda\left(\norm{\hM}-\norm{f^*}\right) = (P_N - P)\cL_{\hM} + P\cL_{\hM} + \lambda\left(\norm{\hM}-\norm{f^*}\right) \\
&\geq -\theta \max\left(r(2\rho^*)^{2\kappa}, \norm{\hM-f^*}_{L_2}^{2\kappa}\right)  + P\cL_{\hM} - \lambda \hrho \geq -\theta r(2 \rho^*)^{2\kappa} - \lambda \rho^* + P\cL_{\hM}\\ 
&\geq -\left(\theta + \frac{1}{2A}\right) r(2\rho^*)^{2\kappa} + P\cL_{\hM}\geq \frac{-r(2 \rho^*)^{2\kappa}}{A} + P \cL_{\hM}.
\end{align*}In particular, if $P\cL_{\hM}> r(2\rho^*)^{2\kappa}/A$ then $P_N \cL_{\hM}^\lambda>0$ which is not possible by construction of $\hM$ so we necessarily have $P\cL_{\hM}\leq r(2\rho^*)^{2\kappa}/A$.
\endproof

Proposition~\ref{prop:algebra} shows that $\hM$ satisfies some estimation and prediction properties on the event $\Omega_0$. Next, we prove that $\Omega_0$ holds with large probability in both subgaussian and bounded frameworks. We start with the subgaussian framework. To that end, we introduce several tools. 

Recall that the $\psi_2$-norm of a real valued random variable $Z$ is defined by 
\begin{equation*}
\norm{Z}_{\psi_2}=\inf\left\{c>0: \E \psi_2(|Z|/c)\leq\psi_2(1)\right\}
\end{equation*}where $\psi_2(u)=\exp(u^2)-1$ for all $u\geq0$. The space $L_{\psi_2}$ of all real valued random variables with finite $\psi_2$-norm is called  the Orlicz space of subgaussian variables. We refer the reader to \cite{MR1113700,MR1890178} for more details on Orlicz spaces. 

We recall several facts on the $\psi_2$-norm and subgaussian processes. First, it follows from Theorem~1.1.5 from \cite{MR3113826} that $\norm{Z}_{\psi_2}\leq \max(K_0, K_1)$ if 
\begin{equation}\label{eq:laplace_psi_2}
\E \exp(\lambda |Z|)\leq \exp\left(\lambda^2 K_1^2\right), \quad \forall \lambda \geq 1/K_0. 
\end{equation}It follows from Lemma~1.2.2 from \cite{MR3113826} that, if $Z$ is a centered $\psi_2$ random variable then, for all $\lambda>0$,
\begin{equation}\label{eq:lemm_psi_2}
\E\exp\left(\lambda Z\right)\leq \exp\left( e \lambda^2 \norm{Z}_{\psi_2}^2\right).
\end{equation}

Then, it follows from Theorem~1.2.1 from \cite{MR3113826} that if $Z_1, \ldots, Z_N$ are independent centered real valued random variables then
\begin{equation}\label{eq:sum_psi_2}
\norm{\sum_{i=1}^N Z_i}_{\psi_2}\leq 16 \left(\sum_{i=1}^N \norm{Z_i}_{\psi_2}^2\right)^{1/2}.
\end{equation}

Finally, let us turn to some properties of subgaussian processes. Let $(T,d)$ be a pseudo-metric space. Let $(X_t)_{t\in T}$ be a random process in $L_{\psi_2}$ such that for all $s,t\in T$, $\norm{X_t-X_s}_{\psi_2}\leq d(s, t)$. It follows from the comment below Theorem~11.2 p.300 in \cite{LT:91} that for all measurable set $A$ and all $s,t\in T$,
\begin{equation*}
\int_A |X_s-X_t|d\bP\leq d(s, t) \bP(A)\psi_2^{-1}\left(\frac{1}{\bP(A)}\right).
\end{equation*}Therefore, it follows from equation~(11.14) in \cite{LT:91} that for every $u>0$,
\begin{equation}\label{eq:dev_subgauss_process}
\bP\left(\sup_{s,t\in T}|X_s-X_t|> c_0 (\gamma_2 + Du)\right)\leq \psi_2(u)^{-1}
\end{equation} where $D$ is the diameter of $(T,d)$, $c_0$ is an absolute constant and $\gamma_2$ is the majorizing measure integral $\gamma(T, d; \psi_2)$ (cf. Chapter~11 in \cite{LT:91}). When $T$ is a subset of $L_2$ and $d$ is the natural metric of $L_2$ it follows from the majorizing measure theorem that $\gamma_2\leq c_1 w(T)$ (cf. Chapter~1 in \cite{Talagrand:05}). 

\begin{lemma}\label{lem:subgauss}
Assume that Assumption~\ref{ass:lipschitz} and Assumption~\ref{ass:sub-gauss} hold. Let $F^\prime\subset F$ then for every $u>0$, with probability at least $1-2\exp(-u^2)$
\begin{equation*}
 \sup_{f,g\in F^\prime}\left|(P-P_N)(\cL_f-\cL_g)\right|\leq \frac{c_0 L}{\sqrt{N}} \left(w(F^\prime) +   u d_{L_2}(F^\prime)\right)
 \end{equation*} where $d$ is the $L_2$ metric and $d_{L_2}(F^\prime)$ is the diameter of $(F^\prime, d)$.
\end{lemma}

\proof
To prove Lemma~\ref{lem:subgauss}, it is enough to show that $\left((P-P_N)\cL_f\right)_{f\in F^\prime}$ has $(L/\sqrt{N})$-subgaussian increments and then to apply \eqref{eq:dev_subgauss_process} where $\gamma_2\sim w(F^\prime)$ in this case. 

Let us prove that for some absolute constant $c_0$: for all $f,g\in F^\prime$,
\begin{equation*}
\norm{(P-P_N)(\cL_f - \cL_g)}_{\psi_2}\leq c_0(L/\sqrt{N}) \norm{f-g}_{L_2}
\end{equation*}
It follows from \eqref{eq:sum_psi_2} that
\begin{align*}
\norm{(P-P_N)(\cL_f - \cL_g)}_{\psi_2} & \leq 16 \left(\sum_{i=1}^N \frac{\norm{(\cL_f-\cL_g)(X_i,Y_i) - \bE(\cL_f-\cL_g)}_{\psi_2}^2}{N^2}\right)^{1/2}=\frac{16}{\sqrt{N}}\norm{\zeta_{f,g}}_{\psi_2}.
\end{align*}where $\zeta_{f,g} = (\cL_f-\cL_g)(X,Y) - \bE(\cL_f-\cL_g)$.Therefore, it only remains to show that $\norm{\zeta_{f,g}}_{\psi_2}\leq c_1 L \norm{f-g}_{L_2}$.

It follows from \eqref{eq:laplace_psi_2}, that the last inequality holds if one proves that for all $\lambda\geq c_1/(L\norm{f-g}_{L_2})$, 
\begin{equation}\label{eq:inter1}
\E\exp\left(\lambda |\zeta_{f,g}|\right)\leq \exp(c_2 \lambda^2 L^2\norm{f-g}_{L_2}^2)
\end{equation} for some absolute constants $c_1$ and $c_2$. To that end, it is enough to prove that, for some absolute constant $c_3$ -- depending only on $c_1$ and $c_2$ -- and all $\lambda>0$, 
\begin{equation*}
\E\exp\left(\lambda |\zeta_{f,g}|\right)\leq 2\exp(c_3 \lambda^2 L^2\norm{f-g}_{L_2}^2).
\end{equation*}

Note that if $Z$ is a real valued random variable and $\eps$ is a Rademacher variable independent of $Z$ then $\E\exp(|Z|)\leq 2 \exp(\eps Z)$. Hence, it follows from a symmetrization argument (cf. Lemma~6.3 in \cite{LT:91}), (a simple version of) the contraction principle (cf. Theorem~4.4 in \cite{LT:91}) and \eqref{eq:lemm_psi_2} that, for all $\lambda>0$,
\begin{align*}
&\E\exp\left(\lambda |\zeta_{f,g}|\right) \leq 2 \E \exp(\lambda \eps \zeta_{f,g}) \leq 2\E\exp\left(2 \lambda \eps(\cL_f-\cL_g)(X, Y)\right)\\
&\leq  2\E\exp\left(2\lambda \eps(f-g)(X)\right)\leq 2 \E \exp\left(c_4 \lambda^2 L^2 \norm{f-g}_{\psi_2}^2\right)
\end{align*}where $\eps$ is a Rademacher variable independent of $(X, Y)$ and where we used in the last but one  inequality that $|\cL_f(X,Y) - \cL_g(X, Y)|\leq |f(X) - g(X)|$ a.s..
\endproof

\begin{proposition}\label{prop:stochastic}
We assume that Assumption~\ref{ass:lipschitz},  \ref{ass:sub-gauss} and \ref{ass:margin} hold. Then the probability measure of $\Omega_0$ is at least as large as the one in \eqref{eq:proba}.
\end{proposition}

\beginproof The proof is based on a peeling argument (cf. \cite{MR1739079}) with respect to the two distances naturally associated with this problem: the regularization norm $\norm{\cdot}$ and the $L_2$-norm $\norm{\cdot}_{L_2}$ associated with the design $X$. The peeling according to $\norm{\cdot}$ is performed along the radii $\rho_j=  2^j \rho^*$ for $j\in\N$ and the peeling according to $\norm{\cdot}_{L_2}$ is performed within the class $\{f\in F : \norm{f-f^*}\leq \rho_j\}:=f^*+\rho_j B$ along the radii $2^i r(\rho_j)$ for all $i=0,1,2, \cdots$ up to a radius such that $2^i r(\rho_j)$ becomes larger than the radius of $f^*+\rho_j B$ in $L_2$, that is for all $i\in I_j$.

We introduce the following partition of the class $F$. We first introduce the "true model", i.e. the subset of $F$ where we want to show that $\hM$ belongs to with high probability:
\begin{equation*}
F_{0,0}=\big\{f\in F : \norm{f-f^*}\leq \rho_0 \mbox{ and } \norm{f-f^*}_{L_2}\leq r(\rho_0)\big\}
\end{equation*}(note that $\rho_0=\rho^*$). Then we peel the remaining set $F\backslash F_{0,0}$ according to the two norms: for every $i\in I_0$,
\begin{equation*}
F_{0,i}=\big\{f\in F : \norm{f-f^*}\leq \rho_0 \mbox{ and } 2^{i-1} r(\rho_0) < \norm{f-f^*}_{L_2}\leq 2^{i} r(\rho_0)\big\},
\end{equation*}for all $j\geq 1$ ,
\begin{equation*}
F_{j,0}=\big\{f\in F : \rho_{j-1}<\norm{f-f^*}\leq \rho_{j} \mbox{ and } \norm{f-f^*}_{L_2}\leq r(\rho_{j})\big\}
\end{equation*}and for every integer $i\in I_j$,
\begin{equation*}
F_{j,i}=\big\{f\in F : \rho_{j-1}<\norm{f-f^*}\leq \rho_j \mbox{ and } 2^{i-1} r(\rho_{j})< \norm{f-f^*}_{L_2}\leq 2^{i} r(\rho_{j})\big\}.
\end{equation*}
We also consider the sets $F_{j,i}^* = \rho_j B \cap (2^i r(\rho_j))B_{L_2}$ for all integers $i$ and $j$.

Let $j$ and $i\in I_j$ be two integers. It follows from Lemma~\ref{lem:subgauss} that for any $u>0$, with probability larger than $1-2\exp(-u^2)$,
\begin{equation}\label{eq:dev_F_j_i}
 \sup_{f\in F_{j,i}}\left|(P-P_N)\cL_f\right|\leq \sup_{f,g\in F_{j,i}^* + f^*}\left|(P-P_N)(\cL_f-\cL_g)\right|\leq \frac{c_0 L}{\sqrt{N}} \left(w(F_{j,i}^*) +   u d_{L_2}(F_{j,i}^*)\right)
 \end{equation} where $d_{L_2}(F^*_{j,i})\leq 2^{i+1} r(\rho_j)$.

Note that for any $\rho>0$, $h:r\to w(\rho B\cap rB_{L_2})/r$ is non-increasing (cf. Lemma~\ref{lem:h} in the Appendix) and note that, by definition of $r(\rho)$ (cf. Definition~\ref{def:function_r}), $h(r(\rho))\leq \theta r(\rho)^{{2\kappa}-1} \sqrt{N} /(CL)$. Since  $h(\cdot)$ is non-increasing, we have $w(F^*_{j,i})/(2^ir(\rho_j)) \leq   h(2^i r(\rho_{j}))\leq h(r(\rho_{j}))\leq \theta r(\rho_{j})^{{2\kappa}-1} \sqrt{N} /(CL)$ and so $w(F_{j,i}^*) \leq
\theta 2^i r(\rho_{j})^{{2\kappa}} \sqrt{N} / (CL)$. Therefore, it follows from \eqref{eq:dev_F_j_i} for $u = \theta \sqrt{N} (2^{(i-1)\vee0} r(\rho_j))^{2\kappa-1}/(2CL)$, if $C\geq 4c_0$ then, with probability at least
\begin{equation}\label{eq:proba_one_peel}
1- 2\exp\left(- \theta^2 N (2^{(i-1)\vee0} r(\rho_j))^{4\kappa-2} /(4C^2L^2) \right),
\end{equation} 
for every $f\in F_{j,i}$, 
\begin{equation*}
|(P-P_N)\cL_f|\leq\theta (2^{(i-1)\vee0} r(\rho_j))^{2\kappa}\leq \theta \max\left(r(2 \max(\norm{f-f^*}, \rho^*))^{2\kappa}, \norm{f-f^*}_{L_2}^{2\kappa}\right).
\end{equation*}
The result follows from a union bound.
\endproof

Now we turn to the proof of Theorem~\ref{thm:main_subgaussian} under the boundedness assumption. The proof follows the same strategy as in the "subgaussian case": we first use Proposition~\ref{prop:algebra} and then show (under the boundedness assumption) that event $\Omega_0$ holds with probability at least as large as the one in \eqref{eq:proba_bounded}.

Similar to Proposition~\ref{prop:stochastic}, we prove the following result under the boundedness assumption.

\begin{proposition}\label{prop:stochastic_boundedness}
We assume that Assumption~\ref{ass:lipschitz},  \ref{ass:bounded} and \ref{ass:margin} hold. Then the probability measure of $\Omega_0$ is at least as large as the one in \eqref{eq:proba_bounded}.
\end{proposition}
\proof
Using the same notation as in the proof of Proposition~\ref{prop:stochastic}, we have for any integer $j$ and $i$ such that $2^{i} r(\rho_j)\leq b$ that by Talagrand's concentration inequality: for any $x>0$, with probability larger than $1-2e^{-x}$,
\begin{equation}\label{eq:talagrand}
Z_{j,i} \leq 2 \E Z_{j,i} + \sigma(\cL_{F_{j,i}})\sqrt{\frac{8x}{N}} + \frac{69\norm{\cL_{F_{j,i}}}_\infty x}{2N}
\end{equation}
where
\begin{equation*}
Z_{j,i} = \sup_{f\in F_{j,i}}|(P-P_N)\cL_f|, \quad \sigma(\cL_{F_{j,i}}) = \sup_{f\in F_{j,i}} \sqrt{\E \cL_f^2} \mbox{ and } \norm{\cL_{F_{j,i}}}_\infty = \sup_{f\in F_{j,i}}\norm{\cL_f}_\infty.
\end{equation*}By the Lipschitz assumption, one has 
\begin{equation*}
\sigma(\cL_{F_{j,i}})\leq 2^{i+1} r(\rho_j) \mbox{ and }\norm{\cL_{F_{j,i}}}_\infty\leq 2b.
\end{equation*}Therefore, it only remains to upper bound the expectation $\E Z_{j,i}$. Let $\eps_1, \ldots, \eps_N$ be a $N$ i.i.d. Rademacher variables independent of the $(X_i, Y_i)$'s. For all function $f$, we set
\begin{equation*}
 P_{N, \eps}f = \frac{1}{N}\sum_{i=1}^N \eps_i f(X_i)
 \end{equation*} It follows from a  symmetrization and a contraction argument (cf. Chapter~4 in \cite{LT:91}) that
\begin{equation*}
\E Z_{j,i}\leq 4 \E \sup_{f\in F_{j,i}}|P_{N, \eps}(f-f^*)| \leq  \frac{4 {\rm Rad}(\rho_j B \cap (2^i r(\rho_j)) B_{L_2})}{\sqrt{N}}\leq (\theta/12)2^i r(\rho_j)^{2\kappa}.
\end{equation*}

Now, we take $x = c_2\theta^2 N (2^{i-1} r(\rho_j))^{4\kappa-2}$ in \eqref{eq:talagrand}  and note that $2^i r(\rho_j)\leq b$ and $\kappa\geq 1$: with probability larger than 
\begin{equation}
1-2\exp(-c_2\theta N(2^i r(\rho_j))^{4\kappa-2}),
\end{equation}
for any $f\in F_{j,i}$,
\begin{align*}
|(P-P_N)\cL_f|& \leq \theta 2^{i-1} r(\rho_j)^{2\kappa}/3 + 2\sqrt{8c_2}\theta\left(2^{i-1} r(\rho_j)\right)^{{2\kappa}} + 69c_2 \theta^2 b(2^{i-1} r(\rho_j))^{4\kappa-2}
\\
&
\leq \theta \left(2^{(i-1)\vee0} r(\rho_j)\right)^{{2\kappa}} \left[
\frac{1}{3} +2\sqrt{8 c_2} + 69 c_2 \theta b (2^i r(\rho_j))^{2\kappa - 2}
\right]
\\
& \leq \theta \left(2^{(i-1)\vee0} r(\rho_j)\right)^{{2\kappa}} \left[
\frac{1}{3} +2\sqrt{8 c_2} + 69 c_2 \theta b^{2\kappa-1}
\right]
\\
& \leq \theta \max\left(r(2 \max(\norm{f-f^*}, \rho^*))^{2\kappa}, \norm{f-f^*}_{L_2}^{2\kappa}\right)
\end{align*}
if $c_2$ is defined by
\begin{equation} 
c_2 = \min\left(\frac{1}{48},\frac{1}{207 \theta b^{2\kappa-1} }\right). 
\end{equation}
We conclude with a union bound.
\endproof

\section{Proof of Theorem~\ref{lower_bound_logistic_S1}}
\label{sec:proof_minimax_logistic_S1}
For the sake of simplicity, assume that $m\geq T$ so $\max(m,T)=m$.
Fix $r\in\{1,\dots,T\}$.
Fix $x>0$ such that $\exp(x)/[1+\exp(x)]\leq b$, we define the set of matrices
\begin{equation*}
\mathcal{C}_{x}
=
\left\{
A \in \mathbb{R}^{m\times r}: \forall (p,q),
A_{p,q} \in \{0,x\}
\right\}
\end{equation*}
and
\begin{equation*}
\mathcal{M}_x = \{A \in \mathbb{R}: A = (B|\dots|B|O), B\in\mathcal{C}_x \}
\end{equation*}
where the block $B$ is repeated $ \lfloor T/r \rfloor $ times
(this construction is taken from~\cite{koltchinskii2011nuclear}). Varshamov-Gilbert bound (Lemma 2.9 in~\cite{tsybakov2009introduction}) implies that there is a finite subset $\mathcal{M}_x^0\subset \mathcal{M}_x$ with ${\rm card}(\mathcal{M}_x^0) \geq 2^{rm/8}+1$ with $0\in\mathcal{M}_x^0$, and for any distinct $A,B\in\mathcal{M}_x^0$,
\begin{equation*}
  \|A-B\|^2_{S_2} \geq \frac{mr \lfloor T/r \rfloor }{8} x^2 \geq \frac{mT}{16} x^2
\end{equation*}
and so
\begin{equation*}
  \frac{1}{mT}\|A-B\|^2_{S_2} \geq \frac{x^2}{16}.
\end{equation*}
Then, for $A\in\mathcal{M}_x^0\setminus\{0\}$,
\begin{align*}
  \mathcal{K}(\mathbb{P}_0,\mathbb{P}_A)
   & = \frac{n}{ mT} \sum_{i=1}^m \sum_{j=1}^T
    \left[ \frac{1}{2}\log \left( \frac{1+\exp(M_{i,j})}{2\exp(M_{i,j})} \right)
    + \frac{1}{2}\log\left( \frac{1+\exp(M_{i,j})}{2} \right)
    \right]
    \\
   & = \frac{n}{mT} \sum_{i=1}^m \sum_{j=1}^T
   \left[
   \log\left(\frac{1+\exp(M_{i,j})}{2}\right) - \frac{1}{2} M_{i,j}
   \right]
   \\
   & \leq
   n \left[ \log\left(\frac{1+\exp(x)}{2}\right) - \frac{1}{2} x \right]
   \\
   & \leq c(b) n x^2 
\end{align*}
where $c(b)>0$ is a constant that depends only on $b$. So:
\begin{equation*}
 \frac{1}{{\rm card}(\mathcal{M}_x^0)-1}
 \sum_{A\in \mathcal{M}_x^0} \mathcal{K}(\mathbb{P}_0,\mathbb{P}_A)
 \leq c(b) n x^2
 \leq c(b) \log({\rm card}(\mathcal{M}_x^0)-1)
\end{equation*}
as soon as we choose
\begin{equation*}
  x
 \leq \sqrt{\frac{\log({\rm card}(\mathcal{M}_x^0)-1)}{n }}
 \leq \sqrt{ \frac{r m \log(2) }{8 n }}
\end{equation*}
(note that the condition $n\geq r m \log(2)/(8b^2)$ implies that
$\exp(x)/[1+\exp(x)] \leq b$).
Then, Theorem 2.5 in~\cite{tsybakov2009introduction} leads to the existence of $\beta,c>0$ such that
$$
\inf_{\widehat{M}} \sup_{A\in \mathcal{M}_x^0}
\mathbb{P}_A\left(
\frac{1}{mT}\|\widehat{M}-A\|^2_{S_2} \geq c \frac{m r}{N}
\right) \geq \beta.
$$
\endproof

\section{Proof of Theorem~\ref{theo:minimax_lower_bound_hinge_S1}}
\label{sec:proof_minimax_hinge_S1}
For the sake of simplicity, assume that $m\geq T$ so $\max(m,T)=m$.
Fix $r\in\{2,\dots,T\}$ and assume that $rT\leq N \leq mT$.

We recall that $\{E_{p,q}: 1\leq p\leq m, 1\leq p\leq T \}$ is the canonical basis of $\R^{m\times T}$. We consider the following ``blocks of coordinates'': for every $1\leq k\leq r-1$ and $1\leq l\leq T$,  
\begin{equation*}
B_{kl} = \left\{E_{p, l}: \frac{(k-1)mT}{N}+1\leq p < \frac{kmT}{N}+1 \right\}
\end{equation*}(note that $(r-1)mT/N+1\leq m$ when $rT\leq N\leq mT$). We also introduce the ``blocks'' of ``remaining'' coordinates:
\begin{equation*}
B_0 = \left\{E_{p, q}: \frac{(r-1)mT}{N}+1\leq p, 1\leq q\leq T \right\}
\end{equation*}

For every $\sigma = (\sigma_{kl})\in\{0,1\}^{(r-1)\times T}$, we denote by $\bP_\sigma$ the probability distribution of a pair $(X, Y)$ taking its values in $\R^{m\times T}\times \{-1, 1\}$ where $X$ is uniformly distributed over the basis $\{E_{p,q}: 1\leq p\leq m, 1\leq p\leq T \}$ and for every $(p,q)\in\{1, \ldots,m\}\times \{1, \ldots, T\}$, 
\begin{equation*}
 \bP_\sigma[Y=1|X=E_{p,q}] =  \left\{
\begin{array}{cc}
\sigma_{kl} & \mbox{ if } E_{p,q}\in B_{kl}\\
1 & \mbox{ otherwise.} 
\end{array}
 \right.
\end{equation*}We also introduce $\eta_\sigma(E_{p,q}) = \E[Y=1|X=E_{p,q}] = 2\bP_\sigma[Y=1|X=E_{p,q}] - 1$. It follows from \cite{zhang2004statistical} that the Bayes rules minimizes the Hinge risk, that is $f_\sigma^*\in\argmin_{f}\E_\sigma (Y-f(X))_+$, where the minimum runs over all measurable functions and $\E_\sigma$ denotes the expectation w.r.t. $(X, Y)$ when $(X, Y)\sim \bP_\sigma$, is achieved by $f^*_\sigma = {\rm sgn}(\eta_\sigma(\cdot))$. Therefore, $f_\sigma^*(\cdot) =\inr{M^*_\sigma, \cdot} $ where for every $(p,q)\in\{1, \ldots,m\}\times \{1, \ldots, T\}$,
\begin{equation*}
 (M^*_\sigma)_{pq} =  \left\{
\begin{array}{cc}
2\sigma_{kl}-1 & \mbox{ if } E_{p,q}\in B_{kl}\\
1 & \mbox{ otherwise.} 
\end{array}
 \right. = \eta_\sigma(E_{p,q}).
 \end{equation*} In particular, $M^*_\sigma$ has a rank at most equal to $r$.

Let $\sigma=(\sigma_{p,q}), \sigma^\prime=(\sigma^\prime_{pq})$ be in $\{0,1\}^{(r-1) T}$. We denote by $\rho(\sigma, \sigma^\prime)$ the Hamming distance between $\sigma$ and $\sigma^\prime$ (i.e. the number of times the coordinates of $\sigma$ and $\sigma^\prime$ are different). We denote by $H(\bP_{\sigma}, \bP_{\sigma^\prime})$ the Hellinger distance between the probability measures $\bP_{\sigma}$ and $\bP_{\sigma^\prime}$. We have
\begin{align*}
H(\bP_{\sigma}, \bP_{\sigma^\prime}) = \int \left( \sqrt{d\bP_{\sigma}} - \sqrt{d \bP_{\sigma^\prime}}\right)^2 = \frac{2\rho(\sigma, \sigma^\prime)}{N}.
\end{align*}
Then, if $\rho(\sigma, \sigma^\prime) = 1$,  it follows that (cf. Section~2.4 in \cite{tsybakov2009introduction}),
\begin{equation*}
H^2(\bP_{\sigma}^{\otimes N}, \bP_{\sigma^\prime}^{\otimes N}) = 2\left(1-\left(1-\frac{H^2(\bP_{\sigma}, \bP_{\sigma^\prime})}{2}\right)^N\right) = 2\left(1-\left(1- \frac{1}{N}\right)^N\right) \leq 2(1-e^{-2}):=\alpha.
\end{equation*}Now, it follows from Theorem~2.12 in \cite{tsybakov2009introduction}, that
\begin{equation}\label{eq:sacha_minimax_lower_bound}
\inf_{\hat\sigma}\max_{\sigma\in\{0,1\}^{(r-1)T}} \E_\sigma^{\otimes N}\norm{\hat\sigma- \sigma}_{l_1}\geq \frac{(r-1)T}{8}\left(1-\sqrt{\alpha(1-\alpha/4)}\right)
\end{equation}where the infimum $\inf_{\hat\sigma}$ runs over all measurable functions $\hat \sigma$ of the data $(X_i,Y_i)_{i=1}^N$ with values in $\R$ (note that Theorem~2.12 in \cite{tsybakov2009introduction} is stated for functions $\hat \sigma$ taking values in $\{0,1\}^{(r-1)T}$ but its is straightforward to extend this result to any $\hat \sigma$ valued in $\R$) and $\E_\sigma^{\otimes N}$ denotes the expectation w.r.t. those data distributed according to $\bP_{\sigma}^{\otimes N}$.

Now, we lower bound the excess risk of any estimator. Let $\hat f$ be an estimator with values in $\R$. Using a truncation argument it is not hard to see that one can restrict the values of $\hat f$ to $[-1,1]$. In that case,  We have
\begin{align*}
&\cE_{hinge}(\hat f) = \E\left[|2\eta_\sigma(X)-1||\hat f(X) - f^*_\sigma(X)|\right] = \E|\hat f(X) - f^*_\sigma(X)| \\
&= \sum_{p,q}|\hat f(E_{p,q}) - f^*_\sigma(E_{p,q})|\bP[X=E_{p,q}]\geq \sum_{kl} \frac{1}{mT}\sum_{E_{p,q}\in B_{kl}} |\hat f(E_{p,q}) - (2\sigma_{pq}-1)| \geq \frac{2}{N}\sum_{kl} |\hat \sigma_{kl} - \sigma_{pq}|
\end{align*}where $\hat\sigma_{kl}$ is the mean of $\{(\hat f(E_{p,q})+1)/2:E_{p,q}\in B_{kl}\}$. Then we obtain,
\begin{equation*}
\inf_{\hat f}\sup_{\sigma\in\{0,1\}^{(r-1)T}}\E_\sigma^{\otimes N} \cE_{hinge}(\hat f) \geq \frac{2}{N}\inf_{\hat\sigma}\max_{\sigma\in\{0,1\}^{(r-1)T}} \E_\sigma^{\otimes N} \norm{\hat \sigma - \sigma}_{l_1}
\end{equation*}and, using \eqref{eq:sacha_minimax_lower_bound}, we get
\begin{equation*}
\inf_{\hat f}\sup_{\sigma\in\{0,1\}^{(r-1)T}}\E_\sigma^{\otimes N} \cE_{hinge}(\hat f) \geq c_0 \frac{rT}{N}
\end{equation*}for $c_0 = \left(1-\sqrt{\alpha(1-\alpha/4)}\right)/4$.

\endproof

\section{Proofs of Section~\ref{sec:bernstein_margin_condition}}
\subsection{Proof of Section~\ref{sub:logistic_loss}}
The proof of Proposition~\ref{prop:margin_logistic_loss} may be found in several papers (cf., for instance, \cite{bartlett2003large}). Let us recall this argument since we will be using it at a starting point to prove the Bernstein condition in the subgaussian case.

{\bf Proof of Proposition~\ref{prop:margin_logistic_loss}:}
The  logistic risk of a function $f:\cX\to \R$ can be written as $ P \ell_f =\mathbb{E}[ g(X,f(X))]$ where for all $x,a\in\R$,  $g(x,a) := \left((1+\eta(x))/2\right) \log\left(1+e^{-a}\right) +  \left((1-\eta(x))/2\right) \log\left(1+e^{a}\right)$ and $\eta(x)=\E[Y|X=x]$ is the conditional expectation of $Y$ given $X=x$.

Since $f^*$ minimizes $f\to P\ell_f$ over the convex class $F$, one has by the first order condition that for every $f\in F$, $\E \partial_2 g(X,f^*(X))(f-f^*)(X)\geq0$. Therefore, it follows from a second order Taylor expansion that the excess logistic loss of every $f\in F$ is such that
\begin{equation}
\label{eq:taylor_logistic_BC}
\cE_{logistic}(f) = P\cL_f\geq \E \left[(f(X)-f^*(X))^2 \int_0^1 (1-u) \delta(f^*(X) + u (f-f^*)(X))du\right]
\end{equation}where $\delta(u)=\partial_2^2g(x,u) = e^u/(1+e^u)^2$ for every $u\in\R$.

Since $|f^*(X)|, |f(X)|\leq b$ a.s. then for every $u\in[0,1]$, $|f^*(X) + u (f-f^*)(X)|\leq 2b$, a.s. and since $\delta(v)\geq \delta(2b)\geq \exp(-2b)/4$ for every $|v|\leq 2b$, it follows from \eqref{eq:taylor_logistic_BC} that $P\cL_f\geq \delta(2b) \norm{f-f^*}_{L_2}^2$.
\endproof

{\bf Proof of Proposition~\ref{prop:bernstein_logistic_subgaussian}:}
Let $t^*\in RB_{l_2}$ be such that $f^*=\inr{\cdot, t^*}$, where $f^*$ is an oracle in $F=\{\inr{\cdot, t}:t\in R B_{l_2}\}$ w.r.t. the logistic loss risk.
Let $f=\inr{\cdot, t}\in F$ for some $t\in RB_{l_2}$. It follows from \eqref{eq:taylor_logistic_BC} that the excess logistic risk of $f$ satisfies
\begin{equation*}
P\cL_f\geq \int_0^1 \E \left[\inr{X, t^*-t}^2 \delta\left(\inr{X, t^*+ u(t-t^*)}\right)\right]du.
\end{equation*}The result will follow if one proves  that for every $t_0,t\in \R^d$, 
\begin{equation}\label{eq:lower_bound_gauss}
\E \left[\inr{X, t}^2 \delta\left(\inr{X, t_0}\right)\right]\geq \frac{\min\left(\pi, \pi^2\left(\norm{t_0}_2\sqrt{2\pi+\norm{t_0}_2^2}\right)^{-1}\right)}{\sqrt{2\pi+\norm{t_0}_2^2} + (\pi-1)\norm{t_0}_2} \frac{\norm{t}_2^2}{8\sqrt{2\pi}}.
\end{equation}

Let us now prove \eqref{eq:lower_bound_gauss}. We write $t = t_0^\perp + \lambda t_0$ where $t_0^\perp$ is a vector orthogonal to $t_0$ and $\lambda\in\R$. Since $\inr{X, t_0^\perp}$ and $\inr{X, t_0}$ are independent random variables, we have
\begin{align*}
\E \left[\inr{X, t}^2 \delta\left(\inr{X, t_0}\right)\right] & = \E \left[\inr{X, t_0^\perp}^2\right] \E \left[\delta\left(\inr{X, t_0}\right)\right] + \lambda^2 \E \left[\inr{X, t_0}^2 \delta\left(\inr{X, t_0}\right)\right],\\
& = \norm{t_0^\perp}_2^2 \E \delta(\norm{t_0}_2g) + \lambda^2\norm{t_0}_2^2 \E g^2 \delta(\norm{t_0}_2g) 
\end{align*} where $g\sim\cN(0,1)$ is standard Gaussian variable and we recall that $\delta(v) = e^v/(1+e^v)^2$ for all $v\in\R$.  Now, it remains to lower bound $\E \delta(\sigma g)$ and $\E g^2 \delta(\sigma g)$ for every $\sigma>0$.

Since $\delta(v)\geq \exp(-|v|)/4$ for all $v\in\R$, one has for all $\sigma>0$,
\begin{align*}
\E \delta(\sigma g)\geq \E \exp(-\sigma |g|)/4 = \exp(\sigma^2/2)\bP[g\geq \sigma]/2
\end{align*}and
\begin{align*}
\E g^2\delta(\sigma g)\geq \E g^2\exp(-\sigma |g|)/4 = (1/2) \exp(\sigma^2/2)\left[(1+\sigma^2)\bP[g\geq \sigma] - \frac{\sigma\exp(-\sigma^2/2)}{\sqrt{2\pi}}\right].
\end{align*}Therefore, for $\sigma=\norm{t_0}_2$,
\begin{align*}
\E \left[\inr{X, t}^2 \delta\left(\inr{X, t_0}\right)\right] \geq &\exp(\sigma^2/2)\bP[g\geq \sigma]\norm{t^\perp_0}_2^2\\ & + 2\lambda^2 \norm{t_0}_2^2 \exp(\sigma^2/2)\left[\left(1+\sigma^2\right)\bP[g\geq \sigma] - \frac{\sigma \exp(-\sigma^2/2)}{\sqrt{2\pi}}\right]
\end{align*}and since $\norm{t}_2^2 = \norm{t^\perp_0}_2^2
+\lambda^2 \norm{t_0}_2^2$, one has,
\begin{equation}\label{eq:excess_risk_mills_ratio}
\E \left[\inr{X, t}^2 \delta\left(\inr{X, t_0}\right)\right] \geq \frac{\norm{t}_2^2}{\sqrt{2\pi}}\min\left\{\left(\frac{1-\Phi(\sigma)}{\phi(\sigma)}\right), (1+\sigma^2)\left(\frac{1-\Phi(\sigma)}{\phi(\sigma)}\right) - \sigma \right\}
\end{equation}where $\phi$ and $\Phi$ denote the standard Gaussian density and distribution functions, respectively. 

We lower bound the right-hand side of \eqref{eq:excess_risk_mills_ratio} using estimates on the Mills ratio $(1-\Phi)/\phi$ that follows from Equation~(10) in \cite{mills_ratio}:  for every $\sigma>0$, 
\begin{equation*}
\frac{1-\Phi(\sigma)}{\phi(\sigma)} > \frac{\pi}{\sqrt{2\pi+\sigma^2} + (\pi-1)\sigma}.
\end{equation*}
\endproof 

\subsection{Proof of Section~\ref{sub:quantile_loss}}
 
{\bf Proof of Proposition~\ref{prop:quantile_loss}:}
 We globally follow a proof of~\cite{elsener2016robust}. We have
 \begin{align*}
 P \mathcal{L}_f
 & = \mathbb{E}[\rho_{\tau}(Y-f(X))-\rho_{\tau}(Y-f^*(X))]  = \mathbb{E}\Bigl\{ \mathbb{E}[\rho_{\tau}(Y-f(X))-\rho_{\tau}(Y-f^*(X))|X ] \Bigr\}.
 \end{align*}
 For all $x\in\cX$, denote by $F_x$ the c.d.f. associated with $f_x$. We have
 \begin{align*}
 & \mathbb{E}[\rho_{\tau}(Y-f(X))|X=x ]
   = (\tau-1) \int_{y < f(x)} (y-f(x))F_x({\rm d}y)
  + \tau \int_{y \geq f(x)} (y-f(x)) F_x({\rm d}y)
  \\
  & = \int_{y \geq f(x)} (y-f(x)) F_x({\rm d}y)
  +(\tau-1) \int_{\mathbb{R}} (y-f(x)) F_x({\rm d}y)
  \\
  & =  \int_{y \geq f(x)} (1-F_x(y)) {\rm d}y
  + (\tau-1) \left( \int_{\mathbb{R}} yF_x({\rm d}y) - f(x) \right)
   = g(x,f(x))
  +(\tau-1)  \int_{\mathbb{R}} yF_x({\rm d}y)
 \end{align*}
 where $g(x,a)= \int_{y \geq a} (1-F_x(y)) {\rm d}y + (1-\tau) a $.
 Note that $\partial_2 g (x,f^*(x))=0$ (can be checked by calculations but
 also obvious from the definition).
 So
 \begin{align*}
 &\mathbb{E}[\rho_{\tau}(Y-f(X))-\rho_{\tau}(Y-f^*(X))|X=x]
  = g(x,f(x))-g(x,f^*(x))
  = \int_{f^*(x)}^{f(x)} (f(x)-u) \partial_2^2 g (x,u) {\rm d}u \\
 & = \int_{f^*(x)}^{f(x)} (f(x)-u) f_x(u) {\rm d}u  \geq \frac{1}{C} \int_{f^*(x)}^{f(x)} (f(x)-u) {\rm d}u 
  = \frac{(f(x)-f^*(x))^2}{2C^2}.
 \end{align*}
It follows that
\begin{equation*}
\cE_{quantile}(f) = P \mathcal{L}_f
\geq
\mathbb{E}\left\{  \frac{(f(X)-f^*(X))^2}{2C} \right\} = \frac{1}{2C}
\|f-f^*\|_{L_2}^2.
\end{equation*}
\endproof

\appendix

\section{Technical lemmas} 

\label{sec:technical_lemmas}
\begin{lemma}\label{lem:Delta}
If $\rho \to r(2\rho)/\rho$ is non-increasing then $\rho\to \Delta(\rho)/\rho$ is non-decreasing.
\end{lemma}
\beginproof
We have for all $\rho>0$
\begin{equation*}
\frac{\Delta(\rho)}{\rho} = \inf_{H \in S \cap (r(2\rho)/\rho) B_{L_2}} \sup_{G \in \partial\norm{\cdot}(M^*)}\inr{H, G}.
\end{equation*}The result follows since $\rho \to S \cap (r(2\rho)/\rho)B_{L_2}$ is non-increasing.
\endproof

\begin{lemma}\label{lem:h}Let $\rho>0$. The function $h:r>0\to w(\rho B\cap rB_{L_2})/r$ is non-increasing. 
\end{lemma}
\beginproof
Let $r_1\geq r_2$. By convexity of $B$ and $B_{L_2}$, we have 
\begin{equation}(\rho B \cap r_1 B_{L_2})/r_1 = (\rho/r_1) B \cap B_{L_2} \subset (\rho/r_2)B \cap B_{L_2}= (\rho B \cap r_2 B_{L_2})/r_2.\end{equation} 
\endproof

\begin{footnotesize}
\bibliographystyle{plain}
\bibliography{biblio}
\end{footnotesize}

\end{document}

%% file: packages_notes.tex
\usepackage{graphicx}
\usepackage{color}
\usepackage{amsmath}
\usepackage{amssymb}
\usepackage{amscd}
\usepackage{bbm}
\usepackage{tikz}
\usepackage{hyperref}
\hypersetup{
    bookmarks=true,         
    colorlinks=true,       
    linkcolor=red,          
    citecolor=green,        
    filecolor=magenta,      
    urlcolor=cyan           
}

\usepackage[utf8]{inputenc}
\usepackage{amsthm}
\usepackage{bbm} 
\usepackage{algorithm}
\usepackage{algorithmic}
\usepackage[small,nohug,heads=vee]{diagrams} 

\usepackage{marginnote}

\usepackage{tcolorbox}
\newtcolorbox{mybox}[1]{colback=red!5!white, colframe=red!75!black,fonttitle=\bfseries, title=#1}

\newtheorem{theo}{Theorem}[section]
\newtheorem{assumption}{Assumption}[section]
\newtheorem{definition}{Definition}[section]
\newtheorem{lemma}{Lemma}[section]
\newtheorem{proposition}{Proposition}[section]

\newtheorem{rmk}{Remark}[section]

%% file: commandes_guillaume.tex
\newcommand{\inr}[1]{\bigl< #1 \bigr>}

\newcommand{\norm}[1]{\left\|#1\right\|}%
\newcommand{\normL}[1]{\left\|#1\right\|_{L_2}}%

\newcommand\eps{\epsilon}

\newcommand{\beginproof}{{\bf Proof. {\hspace{0.2cm}}}}
\def \endproof
{{\mbox{}\nolinebreak\hfill\rule{2mm}{2mm}\par\medbreak}}

\DeclareMathOperator*{\argmin}{argmin}

\DeclareMathOperator*{\rank}{rank}

\def\ds1{\textrm{1\kern-0.25emI}} 

\newcommand \E{\mathbb{E}}
\newcommand \R{\mathbb{R}}
\newcommand \N{\mathbb{N}}

\newcommand \cC{{\cal C}}

\newcommand \cE{{\cal E}}
\newcommand \cF{{\cal F}}

\newcommand \cH{{\cal H}}

\newcommand \cL{{\cal L}}

\newcommand \cN{{\cal N}}

\newcommand \cX{{\cal X}}
\newcommand \cY{{\cal Y}}

\newcommand \bE{{\mathbb E}}

\newcommand \bG{{\mathbb G}}

\newcommand \bP{{\mathbb P}}

\newcommand \bR{{\mathbb R}}

\newcommand{\hM}{\hat f}
\newcommand{\hrho}{\hat \rho}

%% file: 1bitMatrixCompletion_2017_01_27.bbl
\begin{thebibliography}{10}

\bibitem{alquier2015properties}
Pierre Alquier, James Ridgway, and Nicolas Chopin.
\newblock On the properties of variational approximations of gibbs posteriors.
\newblock {\em Journal Of Machine Learning Research}, 17(239):1--41, 2016.

\bibitem{MR2336861}
Jean-Yves Audibert and Alexandre~B. Tsybakov.
\newblock Fast learning rates for plug-in classifiers.
\newblock {\em Ann. Statist.}, 35(2):608--633, 2007.

\bibitem{MR2123199}
Franck Barthe, Olivier Gu{\'e}don, Shahar Mendelson, and Assaf Naor.
\newblock A probabilistic approach to the geometry of the {$l^n_p$}-ball.
\newblock {\em Ann. Probab.}, 33(2):480--513, 2005.

\bibitem{MR2166554}
Peter~L. Bartlett, Olivier Bousquet, and Shahar Mendelson.
\newblock Local {R}ademacher complexities.
\newblock {\em Ann. Statist.}, 33(4):1497--1537, 2005.

\bibitem{bartlett2003large}
Peter~L Bartlett, Michael~I Jordan, and Jon~D McAuliffe.
\newblock Large margin classifiers: Convex loss, low noise, and convergence
  rates.
\newblock In {\em NIPS}, pages 1173--1180, 2003.

\bibitem{bartlett2006convexity}
Peter~L Bartlett, Michael~I Jordan, and Jon~D McAuliffe.
\newblock Convexity, classification, and risk bounds.
\newblock {\em Journal of the American Statistical Association},
  101(473):138--156, 2006.

\bibitem{MR2240689}
Peter~L. Bartlett and Shahar Mendelson.
\newblock Empirical minimization.
\newblock {\em Probab. Theory Related Fields}, 135(3):311--334, 2006.

\bibitem{belloni2011l1}
A.~Belloni and V.~Chernozhukov.
\newblock $\ell$-1-penalized quantile regression in high-dimensional sparse
  models.
\newblock {\em The Annals of Statistics}, 39(1):82--130, 2011.

\bibitem{MR3418717}
Ma{\l}gorzata Bogdan, Ewout van~den Berg, Chiara Sabatti, Weijie Su, and
  Emmanuel~J. Cand{\`e}s.
\newblock S{LOPE}---adaptive variable selection via convex optimization.
\newblock {\em Ann. Appl. Stat.}, 9(3):1103--1140, 2015.

\bibitem{MR2182250}
St{\'e}phane Boucheron, Olivier Bousquet, and G{\'a}bor Lugosi.
\newblock Theory of classification: a survey of some recent advances.
\newblock {\em ESAIM Probab. Stat.}, 9:323--375, 2005.

\bibitem{boyd2011distributed}
Stephen Boyd, Neal Parikh, Eric Chu, Borja Peleato, and Jonathan Eckstein.
\newblock Distributed optimization and statistical learning via the alternating
  direction method of multipliers.
\newblock {\em Foundations and Trends{\textregistered} in Machine Learning},
  3(1):1--122, 2011.

\bibitem{cai2010singular}
Jian-Feng Cai, Emmanuel~J Cand{\`e}s, and Zuowei Shen.
\newblock A singular value thresholding algorithm for matrix completion.
\newblock {\em SIAM Journal on Optimization}, 20(4):1956--1982, 2010.

\bibitem{CandesP10}
E.~J. Cand{\`e}s and Y.~Plan.
\newblock {Matrix Completion With Noise}.
\newblock {\em Proceedings of the IEEE}, 98(6):925--936, 2010.

\bibitem{catoni2004statistical}
Olivier Catoni.
\newblock {\em Statistical learning theory and stochastic optimization: Ecole
  d'Et{\'e} de Probabilit{\'e}s de Saint-Flour, XXXI-2001}, volume~31.
\newblock Springer, 2004.

\bibitem{catoni2007pac}
Olivier Catoni.
\newblock Pac-bayesian supervised classification: the thermodynamics of
  statistical learning.
\newblock {\em arXiv preprint arXiv:0712.0248}, 2007.

\bibitem{MR3113826}
Djalil Chafa{\"{\i}}, Olivier Gu{\'e}don, Guillaume Lecu{\'e}, and Alain Pajor.
\newblock {\em Interactions between compressed sensing random matrices and high
  dimensional geometry}, volume~37 of {\em Panoramas et Synth\`eses [Panoramas
  and Syntheses]}.
\newblock Soci\'et\'e Math\'ematique de France, Paris, 2012.

\bibitem{MR2989474}
Venkat Chandrasekaran, Benjamin Recht, Pablo~A. Parrilo, and Alan~S. Willsky.
\newblock The convex geometry of linear inverse problems.
\newblock {\em Found. Comput. Math.}, 12(6):805--849, 2012.

\bibitem{Cottet2016}
V.~{Cottet} and P.~{Alquier}.
\newblock {1-bit Matrix Completion: PAC-Bayesian Analysis of a Variational
  Approximation}.
\newblock {\em ArXiv e-prints}, April 2016.

\bibitem{MR1864085}
Felipe Cucker and Steve Smale.
\newblock On the mathematical foundations of learning.
\newblock {\em Bull. Amer. Math. Soc. (N.S.)}, 39(1):1--49 (electronic), 2002.

\bibitem{Davenport14}
Mark~A Davenport, Yaniv Plan, Ewout van~den Berg, and Mary Wootters.
\newblock 1-bit matrix completion.
\newblock {\em Information and Inference}, 3(3):189--223, 2014.

\bibitem{MR1932358}
R.~M. Dudley.
\newblock {\em Real analysis and probability}, volume~74 of {\em Cambridge
  Studies in Advanced Mathematics}.
\newblock Cambridge University Press, Cambridge, 2002.
\newblock Revised reprint of the 1989 original.

\bibitem{mills_ratio}
Lutz D{\"u}mbgen.
\newblock Bounding standard gaussian tail probabilities.
\newblock Technical report, University of Bern, 2010.

\bibitem{elsener2016robust}
Andreas Elsener and Sara van~de Geer.
\newblock Robust low-rank matrix estimation.
\newblock {\em arXiv preprint arXiv:1603.09071}, 2016.

\bibitem{MR2721710}
Manuel Garcia-Magari{\~n}os, Anestis Antoniadis, Ricardo Cao, and Wenceslao
  Gonz{\'a}lez-Manteiga.
\newblock Lasso logistic regression, {GS}oft and the cyclic coordinate descent
  algorithm: application to gene expression data.
\newblock {\em Stat. Appl. Genet. Mol. Biol.}, 9:Art. 30, 30, 2010.

\bibitem{halko2011finding}
Nathan Halko, Per-Gunnar Martinsson, and Joel~A Tropp.
\newblock Finding structure with randomness: Probabilistic algorithms for
  constructing approximate matrix decompositions.
\newblock {\em SIAM review}, 53(2):217--288, 2011.

\bibitem{hsieh2014nuclear}
Cho-Jui Hsieh and Peder~A Olsen.
\newblock Nuclear norm minimization via active subspace selection.
\newblock In {\em ICML}, pages 575--583, 2014.

\bibitem{MR0161415}
Peter~J. Huber.
\newblock Robust estimation of a location parameter.
\newblock {\em Ann. Math. Statist.}, 35:73--101, 1964.

\bibitem{klopp2014noisy}
O.~Klopp.
\newblock Noisy low-rank matrix completion with general sampling distribution.
\newblock {\em Bernoulli}, 20(1):282--303, 2014.

\bibitem{MR1892654}
V.~Koltchinskii and D.~Panchenko.
\newblock Empirical margin distributions and bounding the generalization error
  of combined classifiers.
\newblock {\em Ann. Statist.}, 30(1):1--50, 2002.

\bibitem{MR2329442}
Vladimir Koltchinskii.
\newblock Local {R}ademacher complexities and oracle inequalities in risk
  minimization.
\newblock {\em Ann. Statist.}, 34(6):2593--2656, 2006.

\bibitem{MR2829871}
Vladimir Koltchinskii.
\newblock {\em Oracle inequalities in empirical risk minimization and sparse
  recovery problems}, volume 2033 of {\em Lecture Notes in Mathematics}.
\newblock Springer, Heidelberg, 2011.
\newblock Lectures from the 38th Probability Summer School held in Saint-Flour,
  2008, {\'E}cole d'{\'E}t{\'e} de Probabilit{\'e}s de Saint-Flour.
  [Saint-Flour Probability Summer School].

\bibitem{MR2906869}
Vladimir Koltchinskii, Karim Lounici, and Alexandre~B. Tsybakov.
\newblock Nuclear-norm penalization and optimal rates for noisy low-rank matrix
  completion.
\newblock {\em Ann. Statist.}, 39(5):2302--2329, 2011.

\bibitem{koltchinskii2011nuclear}
Vladimir Koltchinskii, Karim Lounici, Alexandre~B Tsybakov, et~al.
\newblock Nuclear-norm penalization and optimal rates for noisy low-rank matrix
  completion.
\newblock {\em The Annals of Statistics}, 39(5):2302--2329, 2011.

\bibitem{lafond2014probabilistic}
Jean Lafond, Olga Klopp, Eric Moulines, and Joseph Salmon.
\newblock Probabilistic low-rank matrix completion on finite alphabets.
\newblock In {\em Advances in Neural Information Processing Systems}, pages
  1727--1735, 2014.

\bibitem{MR2364224}
Guillaume Lecu{\'e}.
\newblock Optimal rates of aggregation in classification under low noise
  assumption.
\newblock {\em Bernoulli}, 13(4):1000--1022, 2007.

\bibitem{HDR-lecue}
Guillaume Lecu{\'e}.
\newblock {\em Interplay between concentration, complexity and geometry in
  learning theory with applications to high dimensional data analysis}.
\newblock Habilitation {\`a} Diriger des Recherches Universit{\'e}. Paris-Est
  Marne-la-vall{\'e}e, December 2011.

\bibitem{MR2933668}
Guillaume Lecu{\'e} and Shahar Mendelson.
\newblock General nonexact oracle inequalities for classes with a
  subexponential envelope.
\newblock {\em Ann. Statist.}, 40(2):832--860, 2012.

\bibitem{LM13}
Guillaume Lecu{\'e} and Shahar Mendelson.
\newblock Learning subgaussian classes: Upper and minimax bounds.
\newblock Technical report, CNRS, Ecole polytechnique and Technion, 2013.

\bibitem{LM_sparsity}
Guillaume Lecu{\'e} and Shahar Mendelson.
\newblock Regularization and the small-ball method {I}: sparse recovery.
\newblock Technical report, CNRS, Ecole Polytechnique and Technion, 2015.

\bibitem{LM_comp}
Guillaume Lecu{\'e} and Shahar Mendelson.
\newblock Regularization and the small-ball method {II}: complexity dependent
  error rates.
\newblock Technical report, CNRS, Ecole Polytechnique and Technion, 2015.

\bibitem{LT:91}
Michel Ledoux and Michel Talagrand.
\newblock {\em Probability in {B}anach spaces}, volume~23 of {\em Ergebnisse
  der Mathematik und ihrer Grenzgebiete (3) [Results in Mathematics and Related
  Areas (3)]}.
\newblock Springer-Verlag, Berlin, 1991.
\newblock Isoperimetry and processes.

\bibitem{mai2015Bayesian}
T.~T. Mai and P.~Alquier.
\newblock A bayesian approach for noisy matrix completion: Optimal rate under
  general sampling distribution.
\newblock {\em Electronic Journal of Statistics}, 9:823--841, 2015.

\bibitem{MR2699823}
Carmen Mak.
\newblock {\em Polychotomous logistic regression via the {L}asso}.
\newblock ProQuest LLC, Ann Arbor, MI, 1999.
\newblock Thesis (Ph.D.)--University of Toronto (Canada).

\bibitem{Mammen1999}
E.~Mammen and A.~Tsybakov.
\newblock Smooth discrimination analysis.
\newblock {\em The Annals of Statistics}, 27(6):1808--1829, 1999.

\bibitem{mazumder2010spectral}
Rahul Mazumder, Trevor Hastie, and Robert Tibshirani.
\newblock Spectral regularization algorithms for learning large incomplete
  matrices.
\newblock {\em Journal of machine learning research}, 11(Aug):2287--2322, 2010.

\bibitem{MR2412631}
Lukas Meier, Sara van~de Geer, and Peter B{\"u}hlmann.
\newblock The group {L}asso for logistic regression.
\newblock {\em J. R. Stat. Soc. Ser. B Stat. Methodol.}, 70(1):53--71, 2008.

\bibitem{mendelson2002improving}
Shahar Mendelson.
\newblock Improving the sample complexity using global data.
\newblock {\em IEEE transactions on Information Theory}, 48(7):1977--1991,
  2002.

\bibitem{MR2075996}
Shahar Mendelson.
\newblock On the performance of kernel classes.
\newblock {\em J. Mach. Learn. Res.}, 4(5):759--771, 2004.

\bibitem{MR2426759}
Shahar Mendelson.
\newblock Obtaining fast error rates in nonconvex situations.
\newblock {\em J. Complexity}, 24(3):380--397, 2008.

\bibitem{MR2590050}
Shahar Mendelson and Joseph Neeman.
\newblock Regularization in kernel learning.
\newblock {\em Ann. Statist.}, 38(1):526--565, 2010.

\bibitem{MR1113700}
M.~M. Rao and Z.~D. Ren.
\newblock {\em Theory of {O}rlicz spaces}, volume 146 of {\em Monographs and
  Textbooks in Pure and Applied Mathematics}.
\newblock Marcel Dekker, Inc., New York, 1991.

\bibitem{MR1890178}
M.~M. Rao and Z.~D. Ren.
\newblock {\em Applications of {O}rlicz spaces}, volume 250 of {\em Monographs
  and Textbooks in Pure and Applied Mathematics}.
\newblock Marcel Dekker, Inc., New York, 2002.

\bibitem{rohde2011estimation}
Angelika Rohde and Alexandre~B Tsybakov.
\newblock Estimation of high-dimensional low-rank matrices.
\newblock {\em The Annals of Statistics}, 39(2):887--930, 2011.

\bibitem{MR3073790}
N.~Sabbe, O.~Thas, and J.-P. Ottoy.
\newblock E{ML}asso: logistic lasso with missing data.
\newblock {\em Stat. Med.}, 32(18):3143--3157, 2013.

\bibitem{srebro2004maximum}
Nathan Srebro, Jason Rennie, and Tommi~S Jaakkola.
\newblock Maximum-margin matrix factorization.
\newblock In {\em Advances in neural information processing systems}, pages
  1329--1336, 2004.

\bibitem{MR2450103}
Ingo Steinwart and Andreas Christmann.
\newblock {\em Support vector machines}.
\newblock Information Science and Statistics. Springer, New York, 2008.

\bibitem{MR3485953}
Weijie Su and Emmanuel Cand{\`e}s.
\newblock S{LOPE} is adaptive to unknown sparsity and asymptotically minimax.
\newblock {\em Ann. Statist.}, 44(3):1038--1068, 2016.

\bibitem{Talagrand:05}
Michel Talagrand.
\newblock {\em The generic chaining}.
\newblock Springer Monographs in Mathematics. Springer-Verlag, Berlin, 2005.
\newblock Upper and lower bounds of stochastic processes.

\bibitem{MR2427362}
Guo-Liang Tian, Man-Lai Tang, Hong-Bin Fang, and Ming Tan.
\newblock Efficient methods for estimating constrained parameters with
  applications to regularized (lasso) logistic regression.
\newblock {\em Comput. Statist. Data Anal.}, 52(7):3528--3542, 2008.

\bibitem{Tsy04}
Alexandre~B. Tsybakov.
\newblock Optimal aggregation of classifiers in statistical learning.
\newblock {\em Ann. Statist.}, 32(1):135--166, 2004.

\bibitem{tsybakov2009introduction}
Alexandre~B Tsybakov.
\newblock {\em Introduction to nonparametric estimation}.
\newblock Springer Series in Statistics, 2009.

\bibitem{MR3526202}
Sara van~de Geer.
\newblock {\em Estimation and testing under sparsity}, volume 2159 of {\em
  Lecture Notes in Mathematics}.
\newblock Springer, [Cham], 2016.
\newblock Lecture notes from the 45th Probability Summer School held in
  Saint-Four, 2015, \'Ecole d'\'Et\'e de Probabilit\'es de Saint-Flour.
  [Saint-Flour Probability Summer School].

\bibitem{MR1739079}
Sara~A. van~de Geer.
\newblock {\em Applications of empirical process theory}, volume~6 of {\em
  Cambridge Series in Statistical and Probabilistic Mathematics}.
\newblock Cambridge University Press, Cambridge, 2000.

\bibitem{van2008high}
Sara~A Van~de Geer.
\newblock High-dimensional generalized linear models and the lasso.
\newblock {\em The Annals of Statistics}, pages 614--645, 2008.

\bibitem{MR1641250}
Vladimir~N. Vapnik.
\newblock {\em Statistical learning theory}.
\newblock Adaptive and Learning Systems for Signal Processing, Communications,
  and Control. John Wiley \& Sons Inc., New York, 1998.
\newblock A Wiley-Interscience Publication.

\bibitem{zhang2004statistical}
Tong Zhang.
\newblock Statistical behavior and consistency of classification methods based
  on convex risk minimization.
\newblock {\em Annals of Statistics}, pages 56--85, 2004.

\end{thebibliography}
